%% file: arxiv.tex
\documentclass{article}

\usepackage{geometry}
\AtBeginDocument{%
  \providecommand\BibTeX{{%
    \normalfont B\kern-0.5em{\scshape i\kern-0.25em b}\kern-0.8em\TeX}}}

\usepackage{setspace}
\usepackage{amsmath}

\newcommand{\mymagenta}{black}
\newcommand{\myred}{black}
\newcommand{\mygreen}{black}
\newcommand{\myblue}{black}

\usepackage{bm}
\usepackage{multirow}
\usepackage{booktabs}

\usepackage{amsmath}
\usepackage{scalerel,amssymb}
\usepackage{mathtools}

\usepackage{svg}
\usepackage{wrapfig}
\usepackage{tikz}
\usepackage{pgfplots}
\usepackage{booktabs}
\usepackage{float}
\usepackage{subfig}
\usepackage{hyperref}

\usepackage{enumitem}

\setenumerate[1]{label={(\arabic*)}}

\usepackage{tikz-uml}
\usepackage{pgf-pie}
\usetikzlibrary{backgrounds}

\usepackage{setspace}
\usepackage[normalem]{ulem}
\usepgfplotslibrary{fillbetween}

\usepackage[framemethod=tikz]{mdframed}
\usetikzlibrary{scopes, backgrounds}

\usepackage{listings}
\usepackage[ruled,vlined,linesnumbered]{algorithm2e}

\newenvironment{acks}[1]{
\section{Acknowledgements}

#1
}

\usepackage{textcomp}
\definecolor{lbcolor}{rgb}{0.93,0.93,0.93}
\definecolor{gnuplot@orange}{RGB}{229,158,0}
\definecolor{gnuplot@purple}{RGB}{148,0,212}
\definecolor{gnuplot@red}{RGB}{200,0,0}
\definecolor{gnuplot@lightblue}{RGB}{87,181,232}
\definecolor{gnuplot@green}{RGB}{0,158,65}
\definecolor{gnuplot@darkblue}{RGB}{0,115,179}
\definecolor{gnuplot@yellow}{RGB}{240,227,66}
\pgfplotscreateplotcyclelist{colorGPL}{%
gnuplot@darkblue,every mark/.append style={fill=gnuplot@darkblue!80!black},mark=o\\%
gnuplot@red,every mark/.append style={fill=gnuplot@red!50!black},mark=square\\%
gnuplot@green,every mark/.append style={fill=gnuplot@green!50!black},mark=triangle*\\%
gnuplot@orange,every mark/.append style={fill=gnuplot@orange!80!black,mark size=2.5pt},mark=x\\%
gnuplot@purple,mark=diamond\\%
black,densely dashed,every mark/.append style={solid,fill=gnuplot@darkblue!80!black},mark=square*\\%
gnuplot@lightblue,densely dashed,every mark/.append style={solid,fill=gnuplot@lightblue!30!black},mark=otimes*\\%
red!80!white,densely dashed,every mark/.append style={solid,fill=red!40!white},mark=oplus*\\%
}
\makeatother

\graphicspath{{svg_with_latex/}}

\usepackage{cite}
\usepackage[numbers,sort&compress]{natbib}
\setcitestyle{numbers}
\setcitestyle{square}

\newif\ifwithacks

\begin{document}

\title{Efficient distributed matrix-free multigrid methods\\ on locally refined meshes for FEM computations}

\author{Peter Munch\thanks{Institute of Material Systems Modeling, Helmholtz-Zentrum Hereon, Max-Planck-Str. 1, 21502 Geesthacht, Germany (\texttt{peter.muench@hereon.de}).}~\thanks{Institute for Computational Mechanics, Department of Mechanical Engineering, Technical University of Munich, Boltzmannstr.~15, 85748 Garching b. M\"unchen, Germany (\texttt{peter.muench@tum.de}).}
\and
Timo Heister\thanks{Clemson University, South Carolina, USA (\texttt{heister@clemson.edu}).}
\and
Laura Prieto Saavedra\thanks{Department of Chemical Engineering, Ecole Polytechnique de Montreal, PO Box 6079, Stn Centre-Ville, H3C 3A7, Montreal, QC, Canada (\texttt{laura.prieto-saavedra@polymtl.ca}).}
\and
Martin Kronbichler\thanks{Department of Information Technology, Uppsala University, Box 337, 75105 Uppsala, Sweden (\texttt{martin.kronbichler@it.uu.se}).}}


\maketitle
  
\input{chapters/abstract}

\noindent \textbf{Key words.}
  \input{chapters/keywords}

\input{chapters/introduction}

\input{chapters/theory}

\input{chapters/implementation}
\input{chapters/modelling}
\input{chapters/performance}
\input{chapters/stokes}
\input{chapters/conclusions}


\input{chapters/appendix}


\bibliographystyle{unsrtnat}
\bibliography{document}

\end{document}

%% file: chapters/abstract.tex
\begin{abstract}
This work studies three multigrid variants for 
matrix-free finite-element computations on locally refined meshes:
geometric local smoothing, geometric global coarsening, and polynomial global coarsening.
We have integrated the algorithms into the same framework---the open-source finite-element library \texttt{deal.II}---, which allows us to make fair comparisons regarding their implementation complexity,
computational efficiency, and parallel scalability as well as
to compare the measurements with theoretically derived performance models. Serial simulations and parallel weak and strong scaling 
on up to 147,456 CPU cores {on 3,072 compute nodes} are presented. 
The results obtained indicate that global coarsening algorithms {\color{\myred}show}
a better parallel behavior for comparable smoothers due to
{\color{\myred}the better load balance particularly on the expensive fine levels}.
In the serial case, the costs of applying hanging-node constraints might be significant, leading to advantages of local smoothing, even though the number of solver iterations {\color{\myred}needed is} slightly higher.
\end{abstract}

%% file: chapters/keywords.tex
multigrid, finite element computations, linear solvers, matrix-free method

%% file: chapters/introduction.tex
\section{Introduction}\label{sec:introduction}

Many solvers {for finite element methods (FEM) rely on efficient solution methods for} second-order partial differential equations (PDEs), e.g., for the Poisson equation:
\begin{align*}
-\Delta u = f,
\end{align*}
{\color{\mymagenta}where $u$ is the solution variable and $f$ is the source term}.
Poisson-like problems {also} {\color{\myred}frequently} occur as subproblems, e.g., in
computational fluid dynamics~\cite{Deville02, KronbichlerDiagneHolmgren2016, ExaDG2020} or in computational plasma physics~\cite{munch2020hyperdeal}.
Efficient {realizations} often rely on adaptively refined meshes to resolve
{geometries} or features in the solution itself and on robust iterative
solvers for such meshes.

Multigrid methods are among the most competitive solvers for such problems \citep{Gholami2016}. 
The three basic steps {of a two-level algorithm} are 
1)~presmoothing, in which the high-frequency error components in the initial guess are removed with a \textit{presmoother},
2)~coarse-grid correction, in which the {\color{\myred}given} problem is solved on a coarse grid, requiring
\textit{intergrid transfer operators} and a \textit{coarse-grid solver}, and 3)~postsmoothing, in which
the high-frequency error components introduced during interpolation are removed with a
\textit{postsmoother}. {Nesting two-level algorithms recursively gives a multigrid algorithm.}
In library implementations, these steps are generally hidden behind operators. {\color{\myred}The latter} can be generally chosen and/or configured 
by the user and strongly depend on the multigrid variant selected. This publication discusses massively parallel multigrid
{variants} for locally refined meshes and the efficient implementation of their operators.


\subsection{Multigrid variants}

Which multigrid approach to choose depends on the way the mesh is generated and on the {underlying}
finite-element space. 
If the mesh is generated by globally refining each cell on a coarse grid recursively, it is a 
natural choice to {\color{\myred}apply} geometric multigrid (abbreviated here as $h$-multigrid), which uses 
the levels of the resulting mesh 
hierarchy as multigrid levels. Alternatively, in the context of high-order finite elements, it is  possible to create levels by reducing the polynomial order of the shape functions $p$ of the elements, while keeping the mesh the same, as done by polynomial multigrid (abbreviated as $p$-multigrid). 
For a very fine, unstructured mesh with low-order elements, it is not as trivial to explicitly construct enough multigrid levels and one might need to 
fall back to non-nested multilevel algorithms~\cite{bittencourt2001nonnested, bramble1991analysis} or algebraic multigrid (AMG; see the review by \citet{Stueben2001}).
These basic multigrid strategies can be nested
in hybrid multigrid solvers \cite{fehn2020hybrid, Rudi2015, OMalley2017b, Lu2014, OMalley2017, Stiller2016a, sundar2012parallel}
and, in doing so, one can exploit the advantages of all of them regarding robustness.
Most {common} are $hp$-multigrid, which combines $h$- and $p$-multigrid, and AMG as 
black-box coarse-grid solver of geometric or polynomial multigrid solvers.

All the {\color{\myred}above-mentioned} multigrid variants are applicable to locally refined meshes. However, local refinement
comes with additional options (local vs. global definition of the multigrid levels) and  
with additional challenges, which are connected, e.g., with the presence of hanging-node constraints.

\subsection{Related work}

Some authors of this study have been involved 
in various publications in the field of multigrid methods in the past. 
Their implementations will be used and
extended in this work.
In \cite{fehn2020hybrid}, 
an efficient hybrid multigrid solver for
discontinuous Galerkin methods (DG) for globally refined 
meshes was presented. It relies on auxiliary-space approximation~\cite{antonietti2017uniform},
i.e., the transfer into a continuous space, as well as on sequential
execution of $p$-multigrid, $h$-multigrid, and AMG.
In \cite{Kronbichler2021}, that solver was extended to leverage locally refined meshes.
In \cite{clevenger2020flexible,kronbichler2019multigrid,clevenger2021comparison, kronbichler2018performance}, 
matrix-free implementations of parallel geometric local-smoothing algorithms for CPU and GPU were investigated and comparisons with AMG were conducted.

\subsection{Our contribution}

In this publication, we will consider three well-known multigrid algorithms
for locally refined meshes for continuous higher-order matrix-free 
FEM:
geometric local smoothing,
geometric global coarsening, and polynomial global coarsening.
We have implemented them into the same framework, which allows us to {\color{\myred}compare their implementation complexity and} performance for a
large variety of problem sizes. This has not been
done in an extensive way in the literature, often using only one of them~\cite{clevenger2020flexible, Sundar2015}. {Furthermore, we rely on matrix-free operator evaluations, which are {\color{\myred}optimal, state-of-the-art} implementations 
in terms of node-level performance on modern hardware~\cite{kronbichler2018performance}, and hence embed the methods in a {\color{\myblue}challenging context in terms of communication costs}
where differences are most pronounced.}


The algorithms presented in this publication have been 
integrated into the open-source finite-element library
\texttt{deal.II}~\cite{dealii2019design} and are mostly part of its 9.3 release~\cite{dealII93}. 
Their implementation has been 
used in \cite{Kronbichler2021} to simulate the flow through a lung geometry and is applied in the 
\texttt{ExaDG} incompressible Navier--Stokes solver~\cite{ExaDG2020}. 
All results of
this publication have been obtained with small benchmark programs
leveraging on the infrastructure of \texttt{deal.II}.
The programs are 
available on GitHub under 
\url{https://github.com/peterrum/dealii-multigrid}.

{\color{\myblue}The results obtained in this publication for continuous FEM are transferable
to the DG case, where one does not have to consider hanging-node constraints but fluxes between
differently refined cells. In the case of auxiliary-space approximation~\cite{Kronbichler2021}, this difference
only involves the finest level and the rest of the multigrid algorithm could be as described in this 
publication.}

The remainder of this work is organized as follows. In Section~\ref{sec:multigrid}, we give a short overview
of multigrid variants applicable to locally refined meshes. Section~\ref{sec:impl} presents
implementation details of our solver, and Section~\ref{sec:modelling} discusses relevant performance models.
{\color{\myred}Sections~\ref{sec:performance_h} and \ref{sec:performance_hp} demonstrate performance results for geometric 
multigrid and polynomial multigrid, and {\color{\mymagenta}in Section~\ref{sec:stokes} the solver is applied to a challenging Stokes problem}.} Finally, Section~\ref{sec:outlook} 
summarizes our conclusions and points to further research directions.

%% file: chapters/theory.tex
\section{Multigrid methods for locally refined meshes}\label{sec:multigrid}

\begin{figure}
\footnotesize
\begin{algorithm}[H]
 \caption{Multigrid V-cycle ${\bm x}\gets\texttt{MultigridVCycle}(\bm b)$ including the copy of $\bm b$ to and of $\bm x$ from the multigrid level(s).}\label{algo:multigrid:copy}
 $[{\color{gray}\bm b^{(0)}, \dots,} \bm b^{(L)}] \gets \bm b$\tcc*{copy to multigrid level(s)}
 $\texttt{VCycleLevel}(L) $\tcc*{Algorithm~\ref{algo:multigrid:cycle} w. input/output $[\bm b^{(0)}/\bm x^{(0)}, \dots, \bm b^{(L)}/\bm x^{(L)}]$}
 \Return $\bm x \gets [{\color{gray}\bm x^{(0)}, \dots,} \bm x^{(L)}] $ \tcc*{copy from multigrid level(s)}
\end{algorithm}
\begin{algorithm}[H]
 \caption{\color{\myred}Actual multigrid V-cycle $\texttt{VCycleLevel}(l)$ called recursively on each level~$l$. 
 It operates on vectors of vectors $[\bm b^{(0)}, \dots, \bm b^{(L)}]$ and $[\bm x^{(0)}, \dots, \bm x^{(L)}]$, which are filled/read in Algorithm~\ref{algo:multigrid:copy},
 and uses the level operators $\bm A^{(l)}$, smoothers, intergrid operators, and coarse-grid solvers, which are 
 set up on or between the multigrid levels. 
 In the case of local smoothing, we distinguish between interior DoFs (not labeled specially)
 and DoFs on the internal boundaries ($\bm x_E^{(l)}$) as well as decompose the level operator $\bm A^{(l)}$ into $\bm A_{SS}^{(l)}$,
  $\bm A_{SE}^{(l)}$, $\bm A_{ES}^{(l)}$, and $\bm A_{EE}^{(l)}$ (see the explanation in Subsection~\ref{sec:multigrid:ls}).
  For global coarsening, $\bm A^{(l)} = \bm A_{SS}^{(l)}$. 
 }\label{algo:multigrid:cycle}
 \uIf{l = 0}{$\bm{x}^{(0)} \gets \texttt{CoarseGridSolver}(\bm A_{SS}^{(0)}, \bm b^{(0)}) $  \tcc*{coarse-grid solver: $\bm A^{(0)}_{SS} \overset{!}{=} \bm A^{(0)}$}}
 \Else{
    $\bm x^{(l)} \gets \texttt{Smoother}(\bm A^{(l)}_{SS}, \bm 0, \bm b^{(l)}) $ \tcc*{presmoothing}
    $\left(\bm r^{(l)}, {\color{gray}\bm r^{(l)}_E}\right) \gets \left( \bm b^{(l)} - \bm A^{(l)}_{SS} \bm x^{(l)}, {\color{gray}- \bm A^{(l)}_{ES} \bm x^{(l)}} \right) $ \tcc*{compute residual}
    $\bm b^{(l-1)} \gets \bm b^{(l-1)} + \texttt{Restrictor} \left(\bm r^{(l)}, {\color{gray}\bm r^{(l)}_E}\right) $ \tcc*{restrict residual}
    $\texttt{VCycleLevel}(l-1) $ \tcc*{recursion}
    $ \left(\bm x^{(l)}, {\color{gray}\bm x^{(l)}_E}\right) \gets \bm x^{(l)} + \texttt{Prolongator}\left(\bm x^{(l-1)}\right) $  \tcc*{prolongation}
    \uIf{local smoothing}{
    {\color{gray}$\bm b^{(l)} \gets \bm b^{(l)} - \bm A^{(l)}_{SE} \bm x^{(l)}_E $}  \tcc*{edge}
    }
    $\bm x^{(l)} \gets \texttt{Smoother}(\bm A_{SS}^{(l)}, \bm x^{(l)}, \bm b^{(l)}) $  \tcc*{postsmoothing}
  }
\end{algorithm}
\end{figure}

Algorithms~\ref{algo:multigrid:copy} and~\ref{algo:multigrid:cycle} present
the basic multigrid algorithm {to solve an equation system of the form $\bm A\bm x=\bm b$ ($\bm A$ is the system matrix, $\bm b$ is the right-hand side vector containing the source term $f$ and the boundary conditions, and $\bm x$ is the solution vector).} It is general enough for meshes obtained by global as well
as by local refinement. In 
the first step, data is transferred to the multigrid levels, after which a
multigrid cycle ({\color{\myred}in this study,} a V-cycle with the steps: presmoothing, computation of the residual, restriction,
solution on the coarser grid, prolongation, and postsmoothing) is {\color{\myblue}performed. Then, the} result
is copied back from the multigrid levels. 
The steps to copy data from and to the multigrid levels are {\color{\myblue}not strictly} needed in all cases, however, are
required by local smoothing and can be used to switch from double to single precision {\color{\myred}to reduce the costs of multigrid
if it is used as a preconditioner~\cite{ljungkvist2014matrix}}.
{The multigrid algorithm is complemented with the algorithms 
of a pre-/postsmoother (e.g., Chebyshev smoother~\cite{adams2003parallel})
and of a coarse-grid solver. 
{\color{\myblue}We use the term ``coarse-grid solver'' also in the 
case that it is applied to a grid that could be coarsened further. This term is rather an indication that the recursion
is terminated.}
Instead of using multigrid as a solver, 
we choose to precondition a conjugate-gradient
solver~\cite{hestenes1952methods} with one multigrid cycle per iteration as this
is often more robust. This algorithm is not presented here.}

\begin{figure}[t]
    \centering
    \def\svgwidth{0.94\columnwidth}
    {\footnotesize 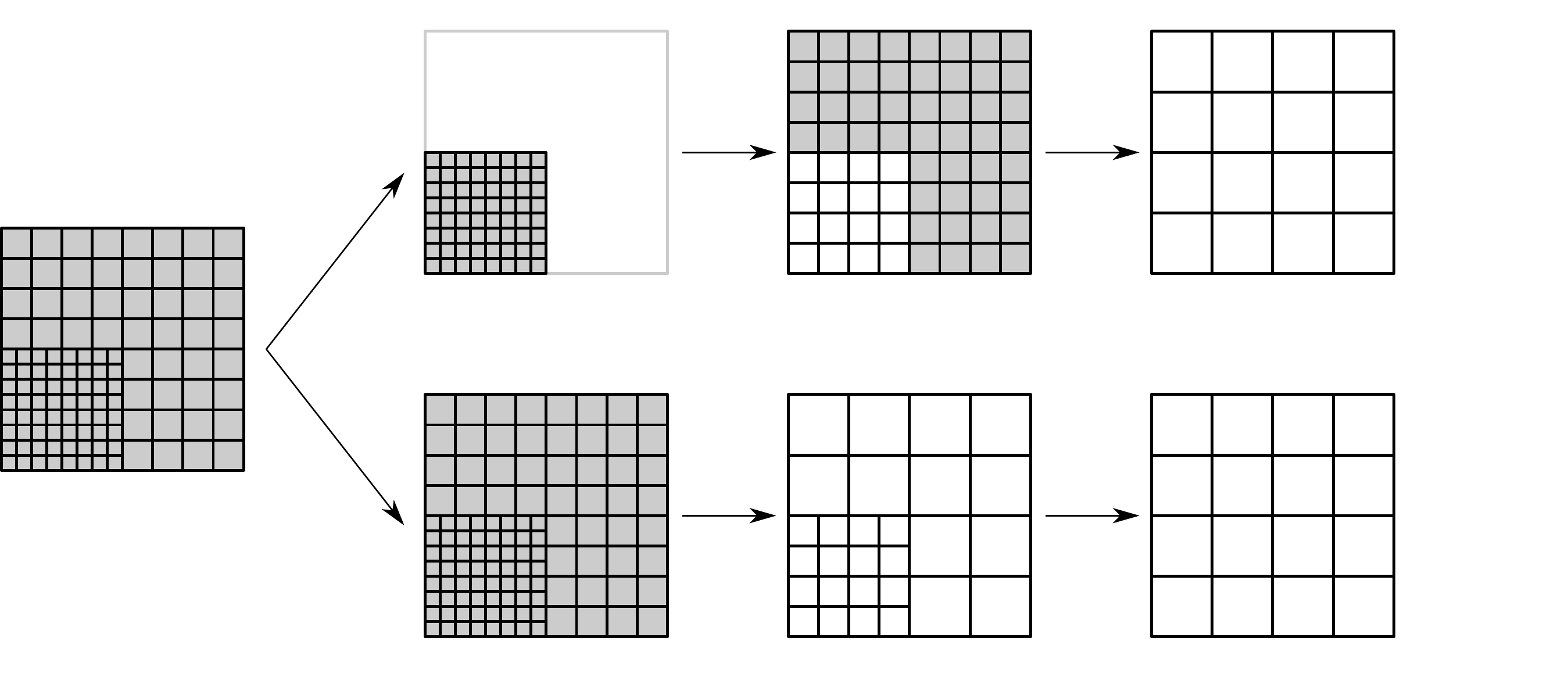}
    
    \vspace{0.3cm}

\caption{Visual comparison of geometric multigrid methods {\color{\myred}for
locally refined meshes}. Top: (geometric) local smoothing; bottom:
  (geometric) global coarsening.
  Local smoothing only considers cells strictly on the same refinement level. This typically introduces an internal boundary (at the refinement edge) when the cells do not cover the
  whole computational domain. Only if they do (here, for level 1 and 0, not for level 2),
  one can switch to a coarse-grid solver.
  Instead, global coarsening {considers the whole domain and} typically introduces hanging nodes on the
multigrid levels. Global coarsening tends to have more cells in total compared to local-smoothing algorithms, but often reduces the number of cells {\color{\myred}per multigrid level} quicker on the finer levels. {\color{\myred}The gray shading indicates active cells.}
}\label{fig:impl:comparsions:geom}

\end{figure}

The various multigrid algorithms for locally refined meshes differ in the construction
of the levels and the concrete details in the implementations of the multigrid steps.
We will consider two types of geometric multigrid methods: geometric local smoothing in Subsection~\ref{sec:multigrid:ls} and geometric global coarsening in Subsection~\ref{sec:multigrid:gc}. Figure~\ref{fig:impl:comparsions:geom} gives a visual comparison of them 
and points out the issues resulting from the local or global definition of the levels, 
which will be discussed 
extensively in the following. Furthermore, we will detail polynomial global coarsening in Subsection~\ref{sec:multigrid:p}.

AMG can be used for the solution on locally refined meshes as well. Since the levels are constructed recursively via the Galerkin
operator $\bm A^{(c)} := \bm R^{( c, f)} \bm A^{(f)} \bm P^{(f, c)}$
with {\color{\myred}the restriction matrix $\bm R^{(c, f)}$ and the prolongation
matrix  $\bm P^{(f, c)}$} constructed algebraically, no distinction regarding the 
local or global {\color{\myred}definition of the levels is possible}. Since we will use AMG in the following only as a
coarse-grid solver, {we refer to the literature for more details}: {\color{\myblue}Clevenger et al.~\cite{clevenger2020flexible} present} 
a scaling comparison between AMG and a matrix-free version of local smoothing for 
a Laplace problem with $\mathcal{Q}_2$, showing the advantages
of matrix-free multigrid methods on modern computing systems.

\subsection{Geometric local smoothing}\label{sec:multigrid:ls}

Geometric local-smoothing algorithms~\cite{brandt1977multi, bastian2006multigrid, mccormick1989multilevel, 
storti1991multigrid, lopez1997algorithmic, schuller2013portable, kanschat2004multilevel,
jouhaud2005multigrid, clevenger2020flexible, janssen2011adaptive} use the refinement hierarchy also for multigrid levels and perform
smoothing refinement level by refinement level: cells of less refined parts of the mesh on the finest level
are skipped (see Figure~\ref{fig:impl:comparsions:geom}) so that hanging-node constraints
do not need to be considered during smoothing. {\color{\myblue}Authors in \cite{mccormick1989multilevel, becker2000multigrid, iliev2001multigrid, wu2006uniform} also have investigated a version of local smoothing in which smoothing is also performed on a halo of a single coarse cell so that hanging-node constraints
need to be applied here as well. We will not consider this
form of local smoothing in the following.}

The fact that domains on each level might not cover the whole computational domain results in multiple issues.
Data needs to be transferred in Algorithm~\ref{algo:multigrid:copy} to and from all multigrid levels 
that are active, i.e., have cells that are not refined.
Furthermore, internal interfaces ({\color{\myred}also known as refinement edges, abbreviated as \textit{edges}}) might result; they need special 
treatment. For details, interested readers are referred to~\cite{kronbichler2019multigrid, janssen2011adaptive}. In the following, we will only summarize key aspects relevant for our investigations.

For the purpose of explanation, let us split the degrees of freedom (DoFs) associated with the cells on an arbitrary level $l$ into 
the interior ones $\bm x^{(l)}_S$ and the ones
at the refinement edges $\bm x^{(l)}_E$   
so that the associated matrix system $\bm A^{(l)} \bm x^{(l)} = \bm b^{(l)}$ has the following block
structure:
\begin{align*}
\left(
\begin{array}{cc}
\bm A^{(l)}_{SS} & \bm A^{(l)}_{SE} \\
\bm A^{(l)}_{ES} & \bm A^{(l)}_{EE} 
\end{array}
\right)
\left(
\begin{array}{c}
\bm x^{(l)}_{S} \\
\bm x^{(l)}_{E} 
\end{array}
\right)
=
\left(
\begin{array}{c}
\bm b^{(l)}_{S} \\
\bm b^{(l)}_{E} 
\end{array}
\right)
\end{align*}
For presmoothing on level $l$, only the contributions from $\bm A^{(l)}_{SS}$ are considered. This is equivalent to applying homogeneous Dirichlet boundary conditions at the 
refinement edges, i.e., $\bm A^{(l)}_{SE}=\bm 0$, $\bm A^{(l)}_{ES} = \bm 0$, and $\bm A^{(l)}_{EE}=\bm I$.
However, when switching to a finer or coarser level, the coupling matrices need to be considered.
The residual to be restricted becomes:
\begin{align}\label{eq:multigrid:ls:edge}
\left(
\begin{array}{c}
\bm r^{(l)}_{S} \\
\bm r^{(l)}_{E} 
\end{array}
\right)
=
\left(
\begin{array}{c}
\bm b^{(l)}_{S} \\
\bm b^{(l)}_{E} 
\end{array}
\right)
-
\left(
\begin{array}{cc}
\bm A^{(l)}_{SS} & \bm A^{(l)}_{SE} \\
\bm A^{(l)}_{ES} & \bm A^{(l)}_{EE} 
\end{array}
\right)
\left(
\begin{array}{c}
\bm x^{(l)}_{S} \\
\bm x^{(l)}_{E} 
\end{array}
\right)
\stackrel{\bm x^{(l)}_{E} {=} \bm 0}{=}
\left(
\begin{array}{c}
\bm 0 \\
\bm b^{(l)}_{E} 
\end{array}
\right)+
\underbrace{
\left(
\begin{array}{c}\bm b^{(l)}_{S} \\
\bm 0 
\end{array}
\right)
-
\left(
\begin{array}{c}
\bm A^{(l)}_{SS} \\
\bm A^{(l)}_{ES} 
\end{array}
\right)
\bm x^{(l)}_{S} 
}_*.
\end{align}
Since $\bm b^{(l)}_{E} $ has been transferred
to the coarser level by Algorithm~\ref{algo:multigrid:copy} already, one only has to restrict the result of term 
$*$ to the coarser level.
During postsmoothing, a modified right-hand side
$
\bm b^{(l)}_{S} := \bm b^{(l)}_{S} - \bm A^{(l)}_{SE} \bm x^{(l)}_{E}
$ 
needs to be considered,
which is equivalent to applying an inhomogeneous Dirichlet boundary condition {\color{\myred}with 
boundary values prescribed by the coarser level}.

A natural choice to partition the multigrid levels for local-smoothing algorithms is to partition the
active level and assign cells on the
lower refinement levels recursively  to the process of their children. 
A simple variant of it is the ``first-child policy''~\cite{clevenger2020flexible}: it recursively assigns the parent cell the
rank of its first child cell.
Since generally in adaptive FEM codes the parents of locally owned and ghost cells are
already available on processes due to tree-like data-structure storage of adaptively refined 
meshes~\cite{Bangerth2011, Burstedde2011}, no additional data structures need to be constructed and saved, but the storage of an 
additional flag (multigrid {rank of the cell}) is enough, leading to 
low memory consumption. Furthermore,
intergrid transfer operations are potentially cheap as data is mainly transferred locally. A 
disadvantage---besides of having to consider the edge constraints---is the potential load imbalance
on the levels{\color{\myred}, as discussed in \cite{clevenger2020flexible}}. This load imbalance could be alleviated by partitioning each level for itself. 
Such alternative partitioning algorithm would lead to similar problems as in the case of global-coarsening 
algorithms (discussed next) and, as a result, the  needed data structures would become more {complex. 
This would object} to the claimed simplicity of the data structures of this {\color{\mygreen}method; hence,}
we will consider geometric-local smoothing only with ``first-child policy'' 
in this publication.

Furthermore, the fact that the transfer to and from the multigrid levels involves all active levels
prevents an early switch to a coarse-grid solver, i.e., 
one can only switch to this solver on levels that are indeed not
locally refined anymore{\color{\myred}, i.e., $\bm A^{(l)}_{SS} \overset{!}{=} \bm A^{(l)}$}.
{On the other hand, the lack of hanging nodes allows the usage of such types of smoothers
that have been developed for uniformly refined meshes, e.g., patch smoothers \cite{arnold1997preconditioning, kanschat2015multigrid}.}

Since geometric local smoothing is the only local-smoothing approach we will consider here, we will
call it simply---{\color{\myblue}as common in the literature}---\textit{local smoothing} in the following.

\subsection{Geometric global coarsening}\label{sec:multigrid:gc}

Geometric global-coarsening algorithms~\cite{becker2000multigrid, becker2007parallel} coarsen all cells simultaneously, 
translating to meshes with hanging nodes also on coarser levels of the multigrid hierarchy
(see Figure~\ref{fig:impl:comparsions:geom}). 
The computational complexity---i.e., the total {\color{\myred}number of cells to be processed}---is slightly higher than in the case of local smoothing
and might be non-optimal for some extreme examples of meshes~\cite{bastian2006multigrid, kanschat2004multilevel}.

The fact that all levels cover the whole computational domain has the advantages that no internal interfaces
have to be considered and the transfer to/from the multigrid levels becomes a simple copy operation
to/from the finest level ($\bm b^{(L)} \gets \bm b$, $\bm x \gets \bm x^{(L)}$). However, hanging
nodes have to be considered during the application of the smoothers on the levels. This is normally not a problem, since codes supporting adaptive meshes
will already have the right infrastructure for this (at least for the active level). {\color{\myred}However, the operator evaluation and the {\color{\myred}applicable smoothers} might be more expensive per cell than in the case of uniformly refined meshes, since the application of hanging-node constraints is not free~\cite{munch2021hn}.
On the other hand, global-coarsening approaches show---for comparable smoothers---a better convergence behavior, which improves
with the number of smoothing iterations~\cite{aulisa2019construction, aulisa2018improved}.}

As the work on the levels generally increases compared to local smoothing, 
it is a valid option to repartition each level separately. 
On the one hand, this implies higher pressure on the transfer operators, since
they need to transfer data between independent meshes\footnote{In \texttt{deal.II}, 
one needs to create a sequence of grids for global coarsening (each with its own hierarchical description). This is generally acceptable, since repartitioning of each level often leads to 
non-overlapping trees so that a single data structure containing all geometric multigrid levels 
would have little benefit for reducing memory consumption.
}, requiring potentially complex
internal data structures, which describe the connectivities, and involved setup routines\footnote{\citet{Sundar2015} present a two-step setup routine: the original fine mesh is coarsened and the resulting ``surrogate mesh'' is repartitioned. For space-filling-curve-based partitioning, this 
approach turns out to be highly efficient.}. On the other hand, it opens the possibility
to control the load balance between processes and the minimal granularity of work per process
(by removing processes on the coarse level in a controlled way, allowing to switch to subcommunicators~\cite{Sundar2015}) and also to apply a coarse-grid solver on any level.
{Furthermore, the construction of full multigrid solvers, which visit the finest level only
a few times, is easier.}

\subsection{Polynomial global coarsening}\label{sec:multigrid:p}

Polynomial global-coarsening algorithms~\cite{Ronquist1987, Maday1988, Hu1995, Hu1997, Guo1998, Guo2000, Atkins2005, Helenbrook2003, Helenbrook2006, Helenbrook2008, Helenbrook2016, Mascarenhas2009, Mascarenhas2010, Sundar2015, Stiller2016b, Stiller2017, OMalley2017, Bassi2003, Hillewaert2006, Liang2009, Luo2006a, Luo2006b, Darmofal2004, Bassi2009, Nastase2006, Premasuthan2009, Fidkowski2005, Ghidoni2014, Shahbazi2009, Jiang2015, Brown2010} are based on keeping the mesh size $h$ constant on all levels,
but reducing the polynomial {\color{\myred}degree} $p$ of shape functions, e.g., to $p=1$. Hence, the multigrid levels
in this case have the same mesh but different polynomial orders. There are various strategies to reduce
the order of the polynomial degree~\cite{helenbrook2008solving, fehn2020hybrid}: the most common is the bisection strategy, which repeatedly halves the degree 
$p^{(c)}= \lfloor  p^{(f)} / 2 \rfloor$. This strategy reduces the number of DoFs in 
the case of a globally refined mesh 
similarly to the geometric multigrid strategies and is a good compromise between the number of iterations
and the cost of a V-cycle~\cite{fehn2020hybrid}.

The statements made in Section~\ref{sec:multigrid:gc} about geometric global coarsening are also valid for polynomial global coarsening. 
However, in contrast to geometric
global coarsening, the levels here do not need to be partitioned for themselves in order to obtain a good load balance, since the
number of unknowns is reduced uniformly in each cell. This leads to a transfer operation that mainly
works on locally owned DoFs.

In the following, we call geometric global coarsening simply \textit{global coarsening}
and polynomial global coarsening \textit{polynomial coarsening}.

%% file: svg_with_latex/local_vs_global.pdf_tex
\begingroup%
  \makeatletter%
  \providecommand\color[2][]{%
    \errmessage{(Inkscape) Color is used for the text in Inkscape, but the package 'color.sty' is not loaded}%
    \renewcommand\color[2][]{}%
  }%
  \providecommand\transparent[1]{%
    \errmessage{(Inkscape) Transparency is used (non-zero) for the text in Inkscape, but the package 'transparent.sty' is not loaded}%
    \renewcommand\transparent[1]{}%
  }%
  \providecommand\rotatebox[2]{#2}%
  \newcommand*\fsize{\dimexpr\f@size pt\relax}%
  \newcommand*\lineheight[1]{\fontsize{\fsize}{#1\fsize}\selectfont}%
  \ifx\svgwidth\undefined%
    \setlength{\unitlength}{776.93407003bp}%
    \ifx\svgscale\undefined%
      \relax%
    \else%
      \setlength{\unitlength}{\unitlength * \real{\svgscale}}%
    \fi%
  \else%
    \setlength{\unitlength}{\svgwidth}%
  \fi%
  \global\let\svgwidth\undefined%
  \global\let\svgscale\undefined%
  \makeatother%
  \begin{picture}(1,0.43605009)%
    \lineheight{1}%
    \setlength\tabcolsep{0pt}%
    \put(0,0){\includegraphics[width=\unitlength,page=1]{local_vs_global.pdf}}%
    \put(0.34884739,0.43018267){\color[rgb]{0,0,0}\makebox(0,0)[t]{\lineheight{1.25}\smash{\begin{tabular}[t]{c}\underline{\smash{\textbf{mg level 2}}}\end{tabular}}}}%
    \put(0.58044158,0.43018267){\color[rgb]{0,0,0}\makebox(0,0)[t]{\lineheight{1.25}\smash{\begin{tabular}[t]{c}\underline{\smash{\textbf{mg level 1}}}\end{tabular}}}}%
    \put(0.81301091,0.43018267){\color[rgb]{0,0,0}\makebox(0,0)[t]{\lineheight{1.25}\smash{\begin{tabular}[t]{c}\underline{\smash{\textbf{mg level 0}}}\end{tabular}}}}%
    \put(0.34884738,0.23515599){\color[rgb]{0,0,0}\makebox(0,0)[t]{\lineheight{1.25}\smash{\begin{tabular}[t]{c}\textit{64 cells}\end{tabular}}}}%
    \put(0.58092915,0.23515599){\color[rgb]{0,0,0}\makebox(0,0)[t]{\lineheight{1.25}\smash{\begin{tabular}[t]{c}\textit{64 cells}\end{tabular}}}}%
    \put(0.81301091,0.23515599){\color[rgb]{0,0,0}\makebox(0,0)[t]{\lineheight{1.25}\smash{\begin{tabular}[t]{c}\textit{16 cells}\end{tabular}}}}%
    \put(0.34884738,0.0020991){\color[rgb]{0,0,0}\makebox(0,0)[t]{\lineheight{1.25}\smash{\begin{tabular}[t]{c}\textit{112 cells}\end{tabular}}}}%
    \put(0.58092915,0.0020991){\color[rgb]{0,0,0}\makebox(0,0)[t]{\lineheight{1.25}\smash{\begin{tabular}[t]{c}\textit{28 cells}\end{tabular}}}}%
    \put(0.81301091,0.0020991){\color[rgb]{0,0,0}\makebox(0,0)[t]{\lineheight{1.25}\smash{\begin{tabular}[t]{c}\textit{16 cells}\end{tabular}}}}%
    \put(0.07971056,0.30731585){\color[rgb]{0,0,0}\makebox(0,0)[t]{\lineheight{1.25}\smash{\begin{tabular}[t]{c}\underline{\smash{\textbf{``active'' level}}}\end{tabular}}}}%
    \put(0.93218291,0.23515599){\color[rgb]{0,0,0}\makebox(0,0)[lt]{\lineheight{1.25}\smash{\begin{tabular}[t]{l}\textit{$\sum=$144 cells}\end{tabular}}}}%
    \put(0.93045399,0.0020991){\color[rgb]{0,0,0}\makebox(0,0)[lt]{\lineheight{1.25}\smash{\begin{tabular}[t]{l}\textit{$\sum$=156 cells}\end{tabular}}}}%
    \put(0.25262273,0.04110443){\color[rgb]{0,0,0}\makebox(0,0)[rt]{\lineheight{1.25}\smash{\begin{tabular}[t]{r}\underline{\smash{\textbf{global coarsening:}}}\end{tabular}}}}%
    \put(0.25262273,0.38922708){\color[rgb]{0,0,0}\makebox(0,0)[rt]{\lineheight{1.25}\smash{\begin{tabular}[t]{r}\underline{\smash{\textbf{local smoothing:}}}\end{tabular}}}}%
    \put(0,0){\includegraphics[width=\unitlength,page=2]{local_vs_global.pdf}}%
    \put(0.35513019,0.35385889){\color[rgb]{0,0,0}\makebox(0,0)[t]{\lineheight{1.25}\smash{\begin{tabular}[t]{c}\textit{\color{red}refinement edge}\end{tabular}}}}%
    \put(0,0){\includegraphics[width=\unitlength,page=3]{local_vs_global.pdf}}%
    \put(0.56171162,0.12797078){\color[rgb]{0,0,0}\makebox(0,0)[t]{\lineheight{1.25}\smash{\begin{tabular}[t]{c}\textit{\color{red}hanging node}\end{tabular}}}}%
    \put(0,0){\includegraphics[width=\unitlength,page=4]{local_vs_global.pdf}}%
  \end{picture}%
\endgroup%

%% file: chapters/implementation.tex
\section{Implementation details}\label{sec:impl}

In this section, we {\color{\myred}will detail efficient implementations of} the multigrid ingredients listed in 
Algorithm~\ref{algo:multigrid:cycle} and needed for the multigrid variants 
for locally refined meshes{, which 
have been considered in Section~\ref{sec:multigrid}}. We start with the handling of constraints. Then, we {\color{\myred}proceed with 
matrix-free evaluations of operator $A$, which is needed
on the active and the multigrid levels, as well as with smoothers} and  coarse-grid
solvers. The discussion of matrix-free transfer operators concludes this section.

\subsection{Handling constraints}\label{sec:imp:constraints}

Constraints need to be considered---with slight differences---in the case both of {local-smoothing and global-}coarsening algorithms. 
First, we impose Dirichlet boundary conditions in a strong form and express them as constraints. Secondly, hanging-node constraints,
which force the solution representation of the refined side to be matching the polynomial
representation of the coarse side, need to be considered to maintain $H^1$ regularity 
of the tentative solution~\cite{shephard1984linear}. 
In a general way, these constraints can be expressed  as
$x_i = \sum_j c_{ij} x_j + b_i$,
where $x_i$ is a constrained DoF, $x_j$ a {\color{\myred}constraining} DoF, $c_{ij}$
the coefficient relating the DoFs, and $b_i$ a real value, which can be 
used to consider inhomogeneities. 
We do not 
eliminate constraints, but use a condensation approach~\cite{bangerth2009data, zave79design}.
We will continue to talk about constraints and their
efficient application in the context of matrix-free
loops in Subsections~\ref{sec:impl:mf} and \ref{sec:impl:transfer}, where we
discuss operators that indeed need to apply constraints.

\subsection{Matrix-free operator evaluation}\label{sec:impl:mf}

\begin{figure}[t]
    \centering
    \def\svgwidth{1.0\columnwidth}
    {\footnotesize 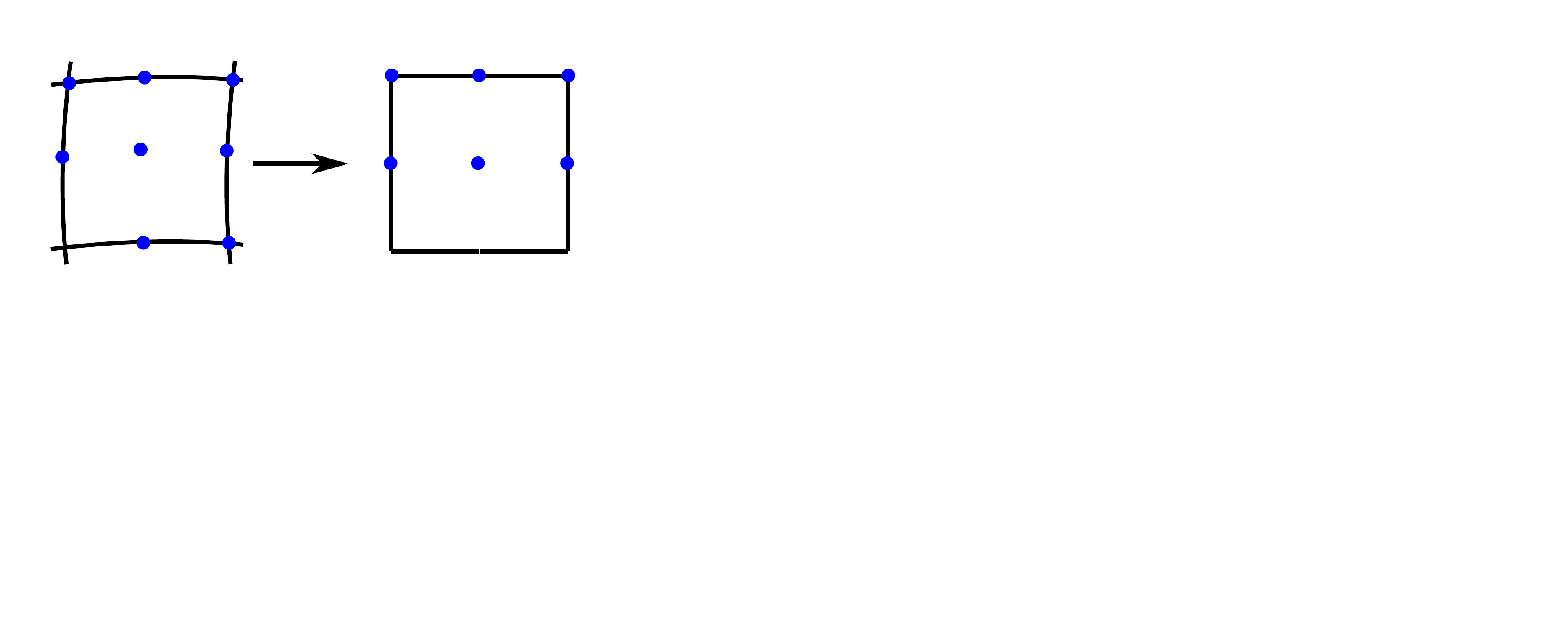}

\caption{Basic steps of a matrix-free operator evaluation according to~\eqref{impl:mf} and of a matrix-free
polynomial prolongation according to~\eqref{eq:impl:loop:categorized} for a single cell $e$.}\label{fig:impl:mf}

\end{figure}

Instead of assembling the system matrix $\bm A$ and performing matrix-vector multiplications of the form $\bm A \bm x$, {\color{\mygreen}the matrix-free operator evaluation computes the underlying finite-element integrals to represent} $\mathcal{A}(\bm x)$. The 
structure of a matrix-free operator evaluation in the context of continuous finite elements generally is:
\begin{align}\label{impl:mf}
\bm v=\mathcal{A}(\bm x)=\sum_e\mathcal{G}^T_e \circ \mathcal{C}^T_e \circ \tilde{\mathcal{S}}^T_e \circ \mathcal{Q}_e \circ \mathcal{S}_e \circ \mathcal{C}_e  \circ\mathcal{G}_e \circ \bm x.
\end{align}
This structure is depicted in Figure~\ref{fig:impl:mf}. For each cell $e$, cell-relevant values are gathered with operator $\mathcal{G}_e$, constraints are---as
discussed in Subsection~\ref{sec:imp:constraints}---applied with $\mathcal{C}_e$, and quantities---like values, gradients, or Hessians---are computed with $\mathcal{S}_e$ at the quadrature points. These quantities are processed by a quadrature-point operation $\mathcal{Q}_e$; the result is integrated and summed into the result vector $\bm v$ by {\color{\myred}applying $\tilde{\mathcal{S}}_e^T$, $\mathcal{C}_e^T$, and $\mathcal{G}_e^T$. In this publication, we consider symmetric
(self-adjoint) PDE operators with $\tilde{\mathcal{S}}_e = \mathcal{S}_e$.}

In the literature, there are both GPU~\citep{kronbichler2019multigrid, anderson2021mfem, ljungkvist2014matrix, ljungkvist2017matrix, kolev2021efficient} and CPU~\cite{Muthing2017, Kronbichler2012, Kronbichler2017a, anderson2021mfem, munch2020hyperdeal} implementations for operations as expressed in~\eqref{impl:mf}, 
which require their own hardware-specific optimizations. For tensor-product 
(quadrilateral and hexahedral) elements,  a technique known as sum factorization~\cite{Melenk99fullydiscrete, Orszag1980} is often applied, which allows to 
replace full interpolations from the local solution values to the quadrature points by a sequence of 1D ones. 
In the context of CPUs, it is an option to 
vectorize over multiple elements~\cite{Kronbichler2012, Kronbichler2017a}, i.e., perform $(\mathcal{S}^T \circ \mathcal{Q} \circ \mathcal{S})_e$ for
multiple elements at the same time. However, in order to be able to do this, the data already needs to be 
laid out in a struct-of-arrays fashion. 
The reshuffling of the data from array-of-structs format to struct-of-arrays format 
and back can be done, e.g., by $\mathcal{G}_e$,
while looping through all elements~\cite{Kronbichler2012}.
We would like to point out that we do not perform $\mathcal{C}_e  \circ\mathcal{G}_e \circ x$
in two steps, but apply the constraints right away, while we are setting up the value of a DoF. 
For the application of hanging-node constraints, we use the special-purpose algorithm introduced in~\cite{munch2021hn},
which is based on the update of the DoF map $\mathcal{G}_e$ and applies in-place sum factorization for the interpolation
during the application of edge and face constraints as well.
Even though the application of hanging-node constraints uses state-of-the-art algorithms {\color{\myred}with small
overhead for high-order finite elements ($<20\%$)}, the additional
steps are {\color{\mygreen}not free particularly} for
linear elements and might lead to load imbalances in a parallel setting if there is a process with
a disproportionately high number of cells with hanging nodes. For a quantitative analysis of this
problem, see the discussion in~\cite{munch2021hn}.


Even though we are using matrix-free algorithms, some of the multigrid ingredients (e.g., smoother and coarse-grid solver) need an explicit representation of the linear operator in the form of a matrix or a part of it (e.g., of its diagonal). In a naive approach {\color{\myblue}relying only} on the matrix-free kernels~\eqref{impl:mf}, one could compute the matrix by simply applying a {unit base vector $e_i$ globally to the operator
$
A(:, i)=\mathcal{A}(e_i)
$
and} hereby reconstructing the matrix column by column.
One can use the same approach also on the cell level:
\begin{align}\label{impl:mf:matrix}
A_e(:, i)=\mathcal{G}^T_e \circ \mathcal{Q}_e \circ \mathcal{S}_e \circ e_i
\end{align}
and assemble the resulting element matrix{\color{\myred}, as usual, also applying the constraints}. Computing the diagonal is slightly more complex
if one does not want to store the complete element matrix and to apply the
constraints during assembly, as described above. {\color{\myred}Instead,} we choose to compute {\color{\myred}the j-th entry
of the locally relevant diagonal contribution via}
\begin{align}\label{impl:mf:diagonal}
d_e(j) = \sum_i \left[ \bm C^T_e \bm A_e(:, i) \right] \bm C_e (i, j) \quad \forall j \in \{ j \;|\; \bm C_{e,ji} \neq 0 \} ,
\end{align}
i.e., we apply the local constraint matrix $\bm C_e$ from the left to the i-th column of the element matrix (computed via~\eqref{impl:mf:matrix}) and apply $\bm C_e$ again from the right. This approach needs 
as many 
basis vector applications {as there are shape functions per cell}. 
The local result can be simply added to the global diagonal via $\displaystyle d = \sum_e \mathcal{G}^T_e \circ d_e$. 
For cells without constrained DoFs, \eqref{impl:mf:diagonal} simplifies to $d_e(i)=\bm A_e(i,i)$ so that one
can use~\eqref{impl:mf:matrix} to compute the element matrix 
column by column and only store the diagonal entries. 

%

\subsection{Smoother}

{For demonstration purposes, we} use Chebyshev iterations around a point-Jacobi method~\cite{adams2003parallel}{, since the construction of efficient and 
robust smoothers for locally refined meshes is out of the scope of this publication. 
The smoother we use only needs an operator {evaluation} and 
its diagonal representation}{\color{\myred}, for which we can use the 
algorithms described in Subsection~\ref{sec:impl:mf}. It is run on the levels either globally or locally, depending on
whether we use a global-coarsening or a local-smoothing multigrid algorithm. 
Interior boundaries
separating the current level from a coarser one are treated as homogeneous Dirichlet boundaries. As a consequence, the
constraint matrix does not have to consider hanging-node
constraints anymore, and the algorithms} significantly simplify---with the result that smoother applications have a higher throughput in the local-smoothing case than in the case
of global coarsening. 

\subsection{Coarse-grid solver}

The algorithms described in Subsection~\ref{sec:impl:mf} allow to set up traditional coarse-grid solvers, e.g., {a Jacobi solver,
a Chebyshev solver}, direct solvers, but also AMG. In this publication, we mostly apply AMG as a coarse-grid 
solver, since we use
it either on very coarse meshes (in this case, it falls back to a direct solver) or for
problems discretized with linear elements, for which AMG solvers are very competitive.

\subsection{Transfer operator}\label{sec:impl:transfer}

The prolongation operator $\mathcal{P}^{(f,c)}$ prolongates the result $\bm x$ from a coarse 
space to a fine space (this includes prolongation from a coarse grid to a fine grid and from a coarser {polynomial} degree to a finer degree):
\begin{align*}
\bm x^{(f)} = \mathcal{P}^{(f,c)} \circ \bm x^{(c)}
\end{align*}
According to the literature~\cite{aulisa2019construction,Sundar2015}, this can be done in three steps:
\begin{align}\label{equ:impl:prolongation:global:dg}
\bm x^{(f)} = \mathcal{W}^{(f)} \circ \tilde{\mathcal{P}}^{(f,c)} \circ \mathcal{C}^{(c)} \circ \bm x^{(c)} 
\end{align}
with $\mathcal{C}^{(c)}$ setting the values of constrained DoFs
on the coarse mesh, particularly
resolving the hanging-node constraints, $\tilde{\mathcal{P}}^{(f,c)}$ performing the
{prolongation} on the discontinuous space as if no hanging nodes were existing, and 
the weighting operator $\mathcal{W}^{(f)}$
zeroing out the DoFs constrained on the fine mesh.

\begin{figure}[t]
    \centering
    \def\svgwidth{1.0\textwidth}
    {\footnotesize 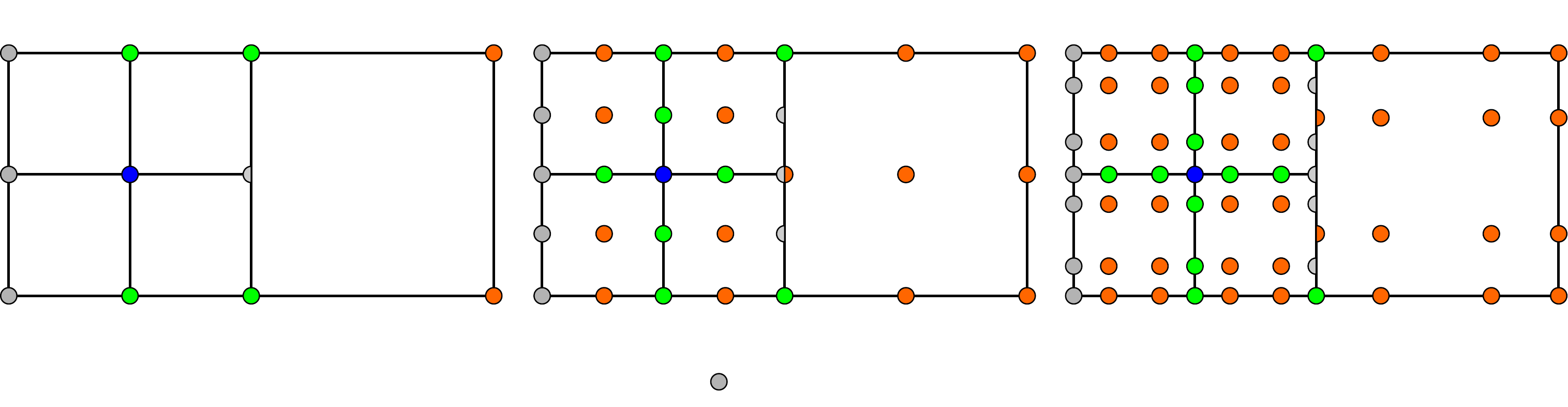}
    

\caption{Example for entries of $\mathcal{W}^{(f)}$ for a mesh with two coarse cells, {\color{\myred}of which one is refined}, and 
Dirichlet boundary at the left face for a scalar continuous Lagrange element of degree $1\le p \le 3$. {\color{\myred}In our
implementation, constrained DoFs do not contribute to the 
valence of constraining DoFs, which results in valences of
one for DoFs inside constraining (coarse) edges/faces.} {In the case that $p > 2$ and all DoFs} of a geometric entity are constrained in the same way, it is enough to store the valence per geometric entity. For efficient access to the
information during the cell loop, one would store the information for all entities of a
``fine cell'' (9 integers in 2D and 27 in 3D, just like in the $p=2$ case).
}\label{fig:impl:mfinterpolation}

\end{figure}

In order to derive a matrix-free implementation, one can express~\eqref{equ:impl:prolongation:global:dg}
for nested meshes as loops over all cells (see also Figure~\ref{fig:impl:mf}):
\begin{align}\label{eq:impl:loop:notcategorized}
\bm x^{(f)} = \sum_{e \in \{ cells\}} \mathcal{S}^{(f)}_e \circ \mathcal{W}^{(f)}_e \circ \mathcal{P}^{(f,c)}_e \circ \mathcal{C}^{(c)}_e \circ \mathcal{G}^{(c)}_e \circ \bm x^{(c)}
\end{align}
Here, $\mathcal{C}^{(c)}_e \circ \mathcal{G}^{(c)}_e$ gathers the cell-relevant coarse DoFs and applies the constraints just as in the case of matrix-free operator evaluations~\eqref{impl:mf}.
$\mathcal{P}^{(f,c)}_e$ performs the {prolongation} onto the fine space for the given (coarse) cell and  $\mathcal{S}^{(f)}_e$ sums the result back to a global vector.
Since multiple elements could add to the same global entry of the vector $\bm x^{(f)}$ during the cell loop, the values to be added have to be weighted
with the inverse of the valence of the corresponding DoF. This is done by $\mathcal{W}^{(f)}_e$, which also ignores constrained DoFs (zero valence) in order to be consistent with~\eqref{equ:impl:prolongation:global:dg}.
Figure~\ref{fig:impl:mfinterpolation} shows, as an example, the values of $\mathcal{W}^{(f)}$ 
for a simple mesh for a scalar Lagrange element of degree $1\le p \le 3$.

\begin{figure}[t]
    \centering
    \def\svgwidth{1.0\columnwidth}
    {\footnotesize 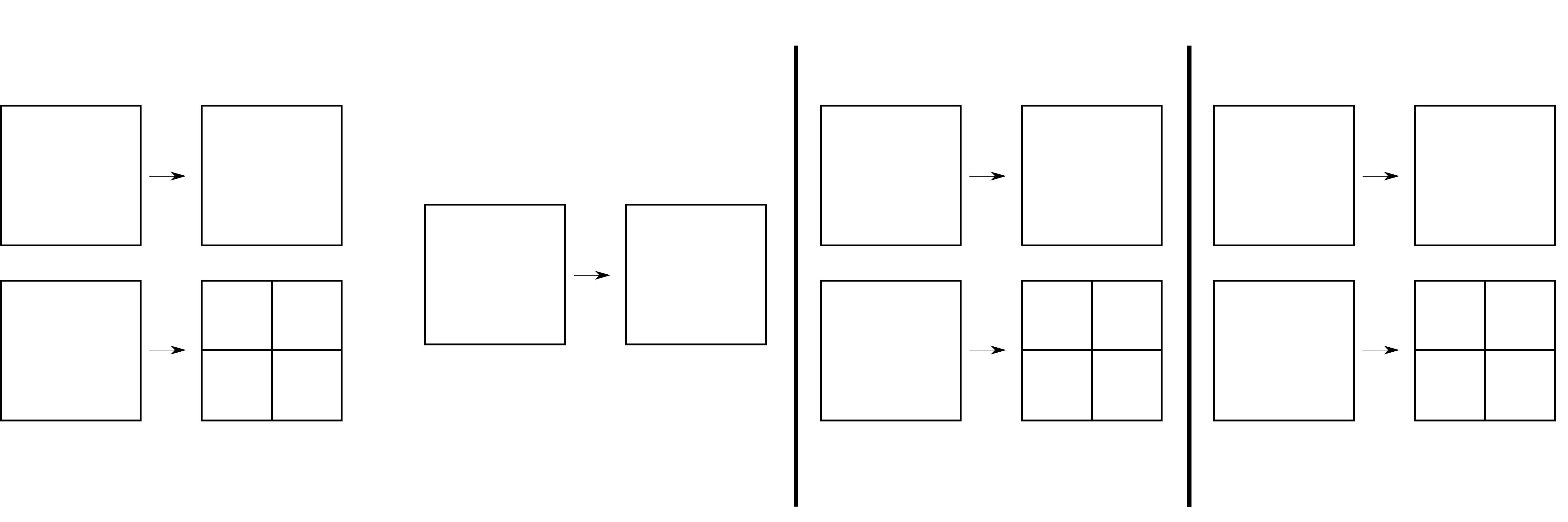}
    
    \vspace{0.3cm}

\caption{Left: Construction of the element prolongation $\mathcal{P}_{e}^{(f,c)}$, based on the \textbf{refinement-type} 
(with/without refinement) and \textbf{element-types pair} (coarse and fine FE). 
Right: Examples for prolongation with and without refinement for equal-degree and 
different-degree finite elements ($p^{(c)}=p^{(f)}=1$ vs. $p^{(c)}=1$, $p^{(f)}=2$). Relevant prolongation types for local
smoothing, global coarsening, and polynomial coarsening are highlighted.
Note that, in the case of global coarsening, two types of prolongation (categories) are needed.}\label{fig:impl:prolongation}

\end{figure}

We construct the element prolongation matrices via:
\begin{align*}
\left(\bm P_{e}^{(f,c)}\right)_{ij}=\left(\bm M^{(f)}_e \right)^{-1} \left( \phi^{(f)}_i,\, \phi^{(c)}_j \right)_{\Omega_e}
\quad\text{with}\quad
\left(\bm M_{e}^{(f)}\right)_{ij} =\left( \phi^{(f)}_i,\, \phi^{(f)}_j \right)_{\Omega_e}
\end{align*}
with $\phi^{(c)}$ being the shape functions defined on the coarse cell and  $\phi^{(f)}$ being the ones 
defined on the ``fine cell''. Instead of treating each ``fine cell'' {on} its own, we group direct children
of the coarse cells together
and define the fine shape functions on the appropriate subregions. As a consequence,
the ``finite element'' on the fine mesh depends both on the actual finite-element type (like on the 
polynomial degree and continuity) and on the refinement type, as indicated in Figure~\ref{fig:impl:prolongation}.
In the case of local smoothing, the finite-element type remains the same and cells are only refined by definition, since cells without refinements are ignored during smoothing. In the case of polynomial global coarsening, the mesh stays the same
and only the polynomial degree is uniformly increased for all cells. In the case of geometric
global coarsening, cells are either refined or not, but the element and its polynomial degree stay
the same. This means that---while in the case of local smoothing and polynomial global coarsening a single
$\bm P^{(f,c)}_e$ is enough---one needs two {variants} in the case of geometric global coarsening {\color{\myred}(note: for non-refined cells, $\bm P^{(f,c)}_e$ is an identity matrix)}.
We define the set of all coarse-fine-cell pairs connected via the same element prolongation matrix
as category $\mathcal{C}$.


Since $\mathcal{P}^{(f,c)}_e=\mathcal{P}^{(f,c)}_{\mathcal{C}(e)}$, i.e., 
all cells of the same category $\mathcal{C}{(e)}$ have the same element prolongation matrix,
and in order to be able to apply them for multiple elements in one go in a 
vectorization-over-elements fashion~\cite{Kronbichler2012} as in the case of matrix-free loops~\eqref{impl:mf}, we loop over the cells type by type so that \eqref{eq:impl:loop:notcategorized} becomes:
\begin{align}\label{eq:impl:loop:categorized}
\bm x^{(f)} = 
\sum_{c}
\sum_{e \in \{ e | \mathcal{C}(e) = c \}}
 \mathcal{S}^{(f)}_e \circ \mathcal{W}^{(f)}_e \circ \mathcal{P}^{(f,c)}_c \circ \mathcal{C}^{(c)}_e \circ \mathcal{G}^{(c)}_e \circ \bm x^{(c)}.
\end{align}

We choose the restriction operator as the transpose of the prolongation
operator
\begin{align*}
\mathcal{R}^{(c, f)} = \left(\mathcal{P}^{(f,c)}\right)^T 
\quad \leftrightarrow \quad 
\mathcal{R}^{(c, f)}_e = \left(\mathcal{P}^{(f,c)}_e\right)^T,
\end{align*}
which implies that the element restriction matrix is the transpose of the 
cell {prolongation} matrix as well.

We conclude this subsection with discussing the appropriate data structures for a transfer operator suitable
both for global geometric and for polynomial coarsening. For details on local smoothing, see~\cite{clevenger2020flexible, kronbichler2019multigrid} and 
the {\color{\mygreen}documentation of} \texttt{MGTransferMatrixFree} class in \texttt{deal.II}~\cite{arndt2017deal}. Since global-coarsening algorithms smoothen
on the complete computational domain, data structures only need to be able to perform two-level 
transfers~\eqref{eq:impl:loop:categorized} independently 
between arbitrary fine (f) and coarse (c) grids. $\mathcal{C}^{(c)}_e \circ \mathcal{G}^{(c)}_e $ 
is identical to $\mathcal{C}_e \circ \mathcal{G}_e $ in the matrix-free loop~\eqref{impl:mf} so that {\color{\myred}specialized
algorithms and data structures~\cite{munch2021hn} can be applied and reused}. 
$\mathcal{S}^{(f)}_e$ needs the indices of DoFs for the given element $e$ in order to be able to scatter the values, and 
$\mathcal{W}^{(f)}_e$
stores the weights of DoFs {\color{\myred}(or of geometric entities - see also the argumentation in Figure~\ref{fig:impl:mfinterpolation})} for the given element $e$.
$\mathcal{P}^{(f,c)}_{\mathcal{C}(e)}$ (and $\mathcal{R}^{(c,f)}_{\mathcal{C}(e)}$) need to be available for each category. 
We regard them 
as general operators and choose the most efficient way of evaluation based on the element
types: simple dense-matrix representation vs. some substructure with sum factorization based on 
smaller 1D {prolongation} matrices. 
{\color{\myred}The category of each cell has to be known for each cell.}


{Alongside these process-local data structures, one needs access to all constraining DoFs on the coarse level 
required during $\mathcal{C}^{(c)}_e \circ \mathcal{G}^{(c)}_e$ 
and to the DoFs of all child cells on the fine level during the process-local part of $\mathcal{S}^{(f)}_e$, which is 
concluded by a data exchange.  
{\color{\myred}If external vectors do not allow access to the required DoFs, we copy data to/from internal} temporal global vectors
with appropriate ghosting~\cite{Kronbichler2012}.}
We choose the {\color{\myred}coarse-side identification of cells} due to its implementation simplicity and structured data access at the price of more ghost transfer.\footnote{\citet{Sundar2015} showed that, by assigning all children
of a cell to the same process, one can easily derive an algorithm that allows to perform the cell-local
prolongation/restriction on the fine side, potentially reducing the amount of data to be communicated
during the transfer. Since we allow levels to be partitioned arbitrarily in our implementation, we
do not use this approach. Furthermore, one should note that the algorithm proposed there does not allow
to apply the constraints $\mathcal{C}_e^{(c)}$ during a single cell loop as in~\eqref{eq:impl:loop:categorized}, but needs a global preprocessing step as 
in~\eqref{equ:impl:prolongation:global:dg}, potentially requiring additional sweeps through
the whole data {with access to the slow main memory}.}
For setup of the communication pattern, we use consensus-based sparse dynamic algorithms~\cite{hoefler2010scalable, arndt2020deal}.

 
For the sake of separation of concerns, one might create three classes {\color{\myred}to implement a global-coarsening transfer operator} as we have done {\color{\myred}in \texttt{deal.II}. The relation of these classes is} shown in
Figure~\ref{fig:transfer:uml}: the multigrid transfer class (\texttt{MGTransferGlobalCoarsening}) delegates the actual transfer tasks to the right two-level implementation (\texttt{MGTwoLevelTransfer}), which performs 
communications needed as well as evaluates \eqref{eq:impl:loop:categorized} for each category and cell by using 
category-specific information from the third class (\texttt{MGTransferSchemes}).

\begin{figure}[!t]

\begin{center}
\begin{tikzpicture}[thick,scale=0.55, every node/.style={transform shape}]

\umlclass[x=-1.4,y=0.7,fill=blue!20, opacity=0.5]{MGTransferMatrixFree}{
}{
prolongate\_and\_add(level, dst, src)\\
restrict\_and\_add(level, dst, src)\\ \\
copy\_to\_mg(dst, src)\\ 
copy\_from\_mg(dst, src)\\ 
interpolate\_to\_mg(dst, src)
}

\umlinterface[x=-9,y=0,fill=blue!20]{MGTransferBase}{}{
prolongate\_and\_add(level, dst, src) = 0\\
restrict\_and\_add(level, dst, src) = 0 \\ \\
copy\_to\_mg(dst, src) = 0\\ 
copy\_from\_mg(dst, src) = 0\\ 
interpolate\_to\_mg(dst, src) = 0
}

\umlclass[x=-1,y=-0.7,fill=blue!20]{MGTransferGlobalCoarsening}{
}{
prolongate\_and\_add(level, dst, src)\\
restrict\_and\_add(level, dst, src)\\ \\
copy\_to\_mg(dst, src)\\ 
copy\_from\_mg(dst, src)\\ 
interpolate\_to\_mg(dst, src)
}

\umlclass[x=7.5,y=-0.7,fill=blue!20]{MGTwoLevelTransfer}{
category\_ptrs\\
level\_dof\_indices\_coarse ($\mathcal{G}$)\\
constraints\_coarse ($\mathcal{C}$)\\
weights\_fine ($\mathcal{W}$)\\
level\_dof\_indices\_fine ($\mathcal{S}$)\\\\
ghosted\_vector\_coarse ($\mathcal{G}$ - optional)\\
ghosted\_vector\_fine ($\mathcal{S}$ - optional)\\
}{
prolongate\_and\_add(vec\_fine, vec\_coarse)\\
restrict\_and\_add(vec\_coarse, vec\_fine)\\ 
interpolate(vec\_coarse, vec\_fine)
}

\umlclass[x=14.5,y=-0.7,fill=blue!20]{MGTransferScheme}{
prolongation\_matrix \\
prolongation\_matrix\_1D \\ \\ 
interpolation\_matrix\\
interpolation\_matrix\_1D \\
}{}

\umlaggreg [mult={$L+1$ }] {MGTransferGlobalCoarsening}{MGTwoLevelTransfer}

\umlaggreg [mult={$\#C$ }] {MGTwoLevelTransfer}{MGTransferScheme}

\umlinherit []{MGTransferGlobalCoarsening}{MGTransferBase}
\umlinherit []{MGTransferMatrixFree}{MGTransferBase}

\end{tikzpicture}
\end{center}
\caption{UML diagram of the global-coarsening transfer operator \texttt{MGTransferGlobalCoarsening} in
\texttt{deal.II}. It implements the base class \texttt{MGTransferBase}, which is also the base class of
\texttt{MGTransferMatrixFree} (local smoothing), and delegates its prolongation/restriction/interpolation tasks to the right \texttt{MGTwoLevelTransfer} instance.
Each of these instances is defined between two levels and is responsible for looping over categories/cells and for evaluating \eqref{eq:impl:loop:categorized} by using the prolongation matrices from the correct \texttt{MGTransferScheme} object. Furthermore, it is responsible for the communication, for which it has 
{\color{\myred}two optional internal vectors with appropriate ghosting}. 
\texttt{MGTwoLevelTransfer} objects can be initialized for geometric or polynomial global coarsening---with 
the consequence that \texttt{MGTransferGlobalCoarsening} {\color{\myred}can handle (global) $h$-, $p$-, and $hp$-multigrid}.
} \label{fig:transfer:uml}
\end{figure}

%% file: svg_with_latex/vmult_vs.pdf_tex
\begingroup%
  \makeatletter%
  \providecommand\color[2][]{%
    \errmessage{(Inkscape) Color is used for the text in Inkscape, but the package 'color.sty' is not loaded}%
    \renewcommand\color[2][]{}%
  }%
  \providecommand\transparent[1]{%
    \errmessage{(Inkscape) Transparency is used (non-zero) for the text in Inkscape, but the package 'transparent.sty' is not loaded}%
    \renewcommand\transparent[1]{}%
  }%
  \providecommand\rotatebox[2]{#2}%
  \newcommand*\fsize{\dimexpr\f@size pt\relax}%
  \newcommand*\lineheight[1]{\fontsize{\fsize}{#1\fsize}\selectfont}%
  \ifx\svgwidth\undefined%
    \setlength{\unitlength}{1074.41603866bp}%
    \ifx\svgscale\undefined%
      \relax%
    \else%
      \setlength{\unitlength}{\unitlength * \real{\svgscale}}%
    \fi%
  \else%
    \setlength{\unitlength}{\svgwidth}%
  \fi%
  \global\let\svgwidth\undefined%
  \global\let\svgscale\undefined%
  \makeatother%
  \begin{picture}(1,0.40153768)%
    \lineheight{1}%
    \setlength\tabcolsep{0pt}%
    \put(0,0){\includegraphics[width=\unitlength,page=1]{vmult_vs.pdf}}%
    \put(0.20979869,0.31340034){\color[rgb]{0,0,0}\makebox(0,0)[t]{\lineheight{0}\smash{\begin{tabular}[t]{c}$\mathcal{C}_e\circ \mathcal{G}_e$\end{tabular}}}}%
    \put(0,0){\includegraphics[width=\unitlength,page=2]{vmult_vs.pdf}}%
    \put(0.40238245,0.31185248){\color[rgb]{0,0,0}\makebox(0,0)[t]{\lineheight{0}\smash{\begin{tabular}[t]{c}$\mathcal{S}_e$\end{tabular}}}}%
    \put(0.58400545,0.38346273){\color[rgb]{0,0,0}\makebox(0,0)[t]{\lineheight{0}\smash{\begin{tabular}[t]{c}{$\mathcal{Q}_e$}\end{tabular}}}}%
    \put(0,0){\includegraphics[width=\unitlength,page=3]{vmult_vs.pdf}}%
    \put(0.60632078,0.3110739){\color[rgb]{0,0,0}\makebox(0,0)[t]{\lineheight{0}\smash{\begin{tabular}[t]{c}$\tilde{\mathcal{S}}^T_e$\end{tabular}}}}%
    \put(0,0){\includegraphics[width=\unitlength,page=4]{vmult_vs.pdf}}%
    \put(0.82640222,0.31305501){\color[rgb]{0,0,0}\makebox(0,0)[t]{\lineheight{0}\smash{\begin{tabular}[t]{c}$\mathcal{G}^T_e\circ \mathcal{C}^T_e$\end{tabular}}}}%
    \put(0,0){\includegraphics[width=\unitlength,page=5]{vmult_vs.pdf}}%
    \put(0.20979869,0.11794574){\color[rgb]{0,0,0}\makebox(0,0)[t]{\lineheight{0}\smash{\begin{tabular}[t]{c}$\mathcal{C}^{(c)}_e\circ \mathcal{G}^{(c)}_e$\end{tabular}}}}%
    \put(0,0){\includegraphics[width=\unitlength,page=6]{vmult_vs.pdf}}%
    \put(0.50709007,0.11639784){\color[rgb]{0,0,0}\makebox(0,0)[t]{\lineheight{0}\smash{\begin{tabular}[t]{c}$\mathcal{P}^{(f,c)}_e$\end{tabular}}}}%
    \put(0,0){\includegraphics[width=\unitlength,page=7]{vmult_vs.pdf}}%
    \put(0.78591512,0.11760041){\color[rgb]{0,0,0}\makebox(0,0)[t]{\lineheight{0}\smash{\begin{tabular}[t]{c}$\mathcal{S}^{(f)}_e$\end{tabular}}}}%
    \put(0,0){\includegraphics[width=\unitlength,page=8]{vmult_vs.pdf}}%
    \put(0.7302087,0.18448198){\color[rgb]{0,0,0}\makebox(0,0)[t]{\lineheight{0}\smash{\begin{tabular}[t]{c}{$\mathcal{W}_e^{(f)}$}\end{tabular}}}}%
    \put(0,0){\includegraphics[width=\unitlength,page=9]{vmult_vs.pdf}}%
    \put(0.19994652,0.39588055){\color[rgb]{0,0,0}\makebox(0,0)[t]{\lineheight{0}\smash{\begin{tabular}[t]{c}evaluation\end{tabular}}}}%
    \put(0,0){\includegraphics[width=\unitlength,page=10]{vmult_vs.pdf}}%
    \put(0.82400564,0.39588055){\color[rgb]{0,0,0}\makebox(0,0)[t]{\lineheight{0}\smash{\begin{tabular}[t]{c}integration\end{tabular}}}}%
    \put(0,0){\includegraphics[width=\unitlength,page=11]{vmult_vs.pdf}}%
    \put(0.21839507,0.00098164){\color[rgb]{0,0,0}\makebox(0,0)[t]{\lineheight{0}\smash{\begin{tabular}[t]{c}coarse (c)\end{tabular}}}}%
    \put(0,0){\includegraphics[width=\unitlength,page=12]{vmult_vs.pdf}}%
    \put(0.79608641,0.00098164){\color[rgb]{0,0,0}\makebox(0,0)[t]{\lineheight{0}\smash{\begin{tabular}[t]{c}fine (f)\end{tabular}}}}%
    \put(0,0){\includegraphics[width=\unitlength,page=13]{vmult_vs.pdf}}%
    \put(0.50018589,0.06489285){\color[rgb]{0,0,0}\makebox(0,0)[t]{\lineheight{0}\smash{\begin{tabular}[t]{c}Fig.~\ref{fig:impl:prolongation}\end{tabular}}}}%
    \put(0,0){\includegraphics[width=\unitlength,page=14]{vmult_vs.pdf}}%
    \put(0.79457598,0.18226557){\color[rgb]{0,0,0}\makebox(0,0)[lt]{\lineheight{0}\smash{\begin{tabular}[t]{l}Fig.~\ref{fig:impl:mfinterpolation}\end{tabular}}}}%
    \put(0,0){\includegraphics[width=\unitlength,page=15]{vmult_vs.pdf}}%
    \put(0.00565714,0.30204061){\color[rgb]{0,0,0}\rotatebox{90}{\makebox(0,0)[t]{\lineheight{0}\smash{\begin{tabular}[t]{c}\underline{\smash{\textbf{op. evaluation}}}\end{tabular}}}}}%
    \put(0.00565714,0.10239726){\color[rgb]{0,0,0}\rotatebox{90}{\makebox(0,0)[t]{\lineheight{0}\smash{\begin{tabular}[t]{c}\underline{\smash{\textbf{poly. prolongation}}}\end{tabular}}}}}%
    \put(0,0){\includegraphics[width=\unitlength,page=16]{vmult_vs.pdf}}%
    \put(0.21947526,0.06340365){\color[rgb]{0,0,0}\makebox(0,0)[t]{\lineheight{0}\smash{\begin{tabular}[t]{c}\cite{munch2021hn}\end{tabular}}}}%
    \put(0,0){\includegraphics[width=\unitlength,page=17]{vmult_vs.pdf}}%
    \put(0.18447949,0.26183705){\color[rgb]{0,0,0}\makebox(0,0)[t]{\lineheight{0}\smash{\begin{tabular}[t]{c}\cite{munch2021hn}\end{tabular}}}}%
    \put(0,0){\includegraphics[width=\unitlength,page=18]{vmult_vs.pdf}}%
  \end{picture}%
\endgroup%

%% file: svg_with_latex/prolongation.pdf_tex
\begingroup%
  \makeatletter%
  \providecommand\color[2][]{%
    \errmessage{(Inkscape) Color is used for the text in Inkscape, but the package 'color.sty' is not loaded}%
    \renewcommand\color[2][]{}%
  }%
  \providecommand\transparent[1]{%
    \errmessage{(Inkscape) Transparency is used (non-zero) for the text in Inkscape, but the package 'transparent.sty' is not loaded}%
    \renewcommand\transparent[1]{}%
  }%
  \providecommand\rotatebox[2]{#2}%
  \newcommand*\fsize{\dimexpr\f@size pt\relax}%
  \newcommand*\lineheight[1]{\fontsize{\fsize}{#1\fsize}\selectfont}%
  \ifx\svgwidth\undefined%
    \setlength{\unitlength}{1435.63866812bp}%
    \ifx\svgscale\undefined%
      \relax%
    \else%
      \setlength{\unitlength}{\unitlength * \real{\svgscale}}%
    \fi%
  \else%
    \setlength{\unitlength}{\svgwidth}%
  \fi%
  \global\let\svgwidth\undefined%
  \global\let\svgscale\undefined%
  \makeatother%
  \begin{picture}(1,0.33305215)%
    \lineheight{1}%
    \setlength\tabcolsep{0pt}%
    \put(0,0){\includegraphics[width=\unitlength,page=1]{prolongation.pdf}}%
    \put(0.10932415,0.27342711){\color[rgb]{0,0,0}\makebox(0,0)[t]{\lineheight{0}\smash{\begin{tabular}[t]{c}\textit{without refinement}\end{tabular}}}}%
    \put(0.10932415,0.04463895){\color[rgb]{0,0,0}\makebox(0,0)[t]{\lineheight{0}\smash{\begin{tabular}[t]{c}\textit{with refinement}\end{tabular}}}}%
    \put(0,0){\includegraphics[width=\unitlength,page=2]{prolongation.pdf}}%
    \put(0.57056647,0.03790919){\color[rgb]{0,0,0}\makebox(0,0)[t]{\lineheight{0}\smash{\begin{tabular}[t]{c}$\mathcal{Q}^1$\end{tabular}}}}%
    \put(0.69496753,0.03790919){\color[rgb]{0,0,0}\makebox(0,0)[t]{\lineheight{0}\smash{\begin{tabular}[t]{c}$\mathcal{Q}^1$\end{tabular}}}}%
    \put(0.82162018,0.03790919){\color[rgb]{0,0,0}\makebox(0,0)[t]{\lineheight{0}\smash{\begin{tabular}[t]{c}$\mathcal{Q}^1$\end{tabular}}}}%
    \put(0.94602121,0.03790919){\color[rgb]{0,0,0}\makebox(0,0)[t]{\lineheight{0}\smash{\begin{tabular}[t]{c}$\mathcal{Q}^2$\end{tabular}}}}%
    \put(0.31500955,0.14992643){\color[rgb]{0,0,0}\makebox(0,0)[t]{\lineheight{0}\smash{\begin{tabular}[t]{c}FE$_c$\end{tabular}}}}%
    \put(0.44447671,0.14992643){\color[rgb]{0,0,0}\makebox(0,0)[t]{\lineheight{0}\smash{\begin{tabular}[t]{c}FE$_f$\end{tabular}}}}%
    \put(0,0){\includegraphics[width=\unitlength,page=3]{prolongation.pdf}}%
    \put(0.58236759,0.00188355){\color[rgb]{0,0,0}\makebox(0,0)[t]{\lineheight{0}\smash{\begin{tabular}[t]{c}\textit{local smoothing}\end{tabular}}}}%
    \put(0.63586294,0.28896292){\color[rgb]{0,0,0}\makebox(0,0)[t]{\lineheight{0}\smash{\begin{tabular}[t]{c}\textit{global coarsening}\end{tabular}}}}%
    \put(0.88297633,0.28896292){\color[rgb]{0,0,0}\makebox(0,0)[t]{\lineheight{0}\smash{\begin{tabular}[t]{c}\textit{polynomial coarsening}\end{tabular}}}}%
    \put(0,0){\includegraphics[width=\unitlength,page=4]{prolongation.pdf}}%
    \put(0.10932415,0.32517113){\color[rgb]{0,0,0}\makebox(0,0)[t]{\lineheight{0}\smash{\begin{tabular}[t]{c}\underline{\smash{\textbf{refinement type}}}\end{tabular}}}}%
    \put(0.37998516,0.32517113){\color[rgb]{0,0,0}\makebox(0,0)[t]{\lineheight{0}\smash{\begin{tabular}[t]{c}\underline{\smash{\textbf{element types}}}\end{tabular}}}}%
    \put(0.63586294,0.3247857){\color[rgb]{0,0,0}\makebox(0,0)[t]{\lineheight{0}\smash{\begin{tabular}[t]{c}\underline{\smash{\textbf{example 1:}}}\end{tabular}}}}%
    \put(0.8770697,0.3247857){\color[rgb]{0,0,0}\makebox(0,0)[t]{\lineheight{0}\smash{\begin{tabular}[t]{c}\underline{\smash{\textbf{example 2:}}}\end{tabular}}}}%
    \put(0,0){\includegraphics[width=\unitlength,page=5]{prolongation.pdf}}%
  \end{picture}%
\endgroup%

%% file: svg_with_latex/hn_weights_alternative_horizontal.pdf_tex
\begingroup%
  \makeatletter%
  \providecommand\color[2][]{%
    \errmessage{(Inkscape) Color is used for the text in Inkscape, but the package 'color.sty' is not loaded}%
    \renewcommand\color[2][]{}%
  }%
  \providecommand\transparent[1]{%
    \errmessage{(Inkscape) Transparency is used (non-zero) for the text in Inkscape, but the package 'transparent.sty' is not loaded}%
    \renewcommand\transparent[1]{}%
  }%
  \providecommand\rotatebox[2]{#2}%
  \newcommand*\fsize{\dimexpr\f@size pt\relax}%
  \newcommand*\lineheight[1]{\fontsize{\fsize}{#1\fsize}\selectfont}%
  \ifx\svgwidth\undefined%
    \setlength{\unitlength}{872.65497668bp}%
    \ifx\svgscale\undefined%
      \relax%
    \else%
      \setlength{\unitlength}{\unitlength * \real{\svgscale}}%
    \fi%
  \else%
    \setlength{\unitlength}{\svgwidth}%
  \fi%
  \global\let\svgwidth\undefined%
  \global\let\svgscale\undefined%
  \makeatother%
  \begin{picture}(1,0.2593368)%
    \lineheight{1}%
    \setlength\tabcolsep{0pt}%
    \put(0,0){\includegraphics[width=\unitlength,page=1]{hn_weights_alternative_horizontal.pdf}}%
    \put(0.15951577,0.24970438){\color[rgb]{0,0,0}\makebox(0,0)[t]{\lineheight{0}\smash{\begin{tabular}[t]{c}$p=1$\end{tabular}}}}%
    \put(0.50257841,0.24970438){\color[rgb]{0,0,0}\makebox(0,0)[t]{\lineheight{0}\smash{\begin{tabular}[t]{c}$p=2$\end{tabular}}}}%
    \put(0.83889853,0.25063043){\color[rgb]{0,0,0}\makebox(0,0)[t]{\lineheight{0}\smash{\begin{tabular}[t]{c}$p=3$\end{tabular}}}}%
    \put(0.47583051,0.0097895){\color[rgb]{0,0,0}\makebox(0,0)[lt]{\lineheight{0}\smash{\begin{tabular}[t]{l}0\end{tabular}}}}%
    \put(0,0){\includegraphics[width=\unitlength,page=2]{hn_weights_alternative_horizontal.pdf}}%
    \put(0.52052171,0.0097895){\color[rgb]{0,0,0}\makebox(0,0)[lt]{\lineheight{0}\smash{\begin{tabular}[t]{l}1\end{tabular}}}}%
    \put(0,0){\includegraphics[width=\unitlength,page=3]{hn_weights_alternative_horizontal.pdf}}%
    \put(0.56521281,0.0097895){\color[rgb]{0,0,0}\makebox(0,0)[lt]{\lineheight{0}\smash{\begin{tabular}[t]{l}1/2\end{tabular}}}}%
    \put(0,0){\includegraphics[width=\unitlength,page=4]{hn_weights_alternative_horizontal.pdf}}%
    \put(0.6202175,0.0097895){\color[rgb]{0,0,0}\makebox(0,0)[lt]{\lineheight{0}\smash{\begin{tabular}[t]{l}1/4\end{tabular}}}}%
    \put(0.43726619,0.0097895){\color[rgb]{0,0,0}\makebox(0,0)[rt]{\lineheight{0}\smash{\begin{tabular}[t]{r}$\mathcal{W}_i^{(f)}$:\end{tabular}}}}%
    \put(0,0){\includegraphics[width=\unitlength,page=5]{hn_weights_alternative_horizontal.pdf}}%
  \end{picture}%
\endgroup%

%% file: chapters/modelling.tex
\section{Performance modeling}\label{sec:modelling}

{\color{\myred}In Section~\ref{sec:impl}, we have presented efficient implementations
of the operators in Algorithms~\ref{algo:multigrid:copy} and \ref{algo:multigrid:cycle} for 
local smoothing, global coarsening, and polynomial coarsening.
Most of the discussion was independent of the multigrid
variant chosen, highlighting the similarities from the
implementation point of view. The main differences arise 
naturally from the local or global definition of levels. E.g.,
one might need to consider---possibly expensive---hanging-node constraints during 
matrix-free loops when doing
global coarsening or polynomial coarsening. 
Local smoothing, on the other hand, has the disadvantage
of performing additional
steps: 1) global transfer to/from multigrid levels and
2) special treatment of edges during smoothing, 
computation of the residual, and modification of the right-hand side vector for postsmoothing. 
In the following sections, we will quantify the influence of the costs of the 
potentially more expensive operator evaluations 
and of the additional operator evaluations, related to
the choice of the multigrid level definition.}

{\color{\myred}Our primary goal is
to minimize the \textbf{time to solution}. It consists of
setup costs and the actual solve time, which is the product of the \textbf{number of iterations} times \textbf{the time per iteration}. We will disregard the setup costs, since 
they normally amortize in time-dependent simulations, where
one does not remesh every time step, and ways to optimize
the setup of global-coarsening algorithms are known in the literature~\cite{Sundar2015}. The time of the solution process strongly depends
on the choice of the smoother, which influences the number of iterations
and, as a result, also the time to solution. 
Since different iteration numbers might distort the view
on the performance of the actual multigrid algorithm, we will also consider the value of the time per iteration as an important indicator of the computational performance of multigrid algorithms particularly to quantify the additional
costs in an iteration.
}

In order to get a first estimate of the benefits of an algorithm compared to another, one can derive following
metrics purely from geometrical information:
\begin{itemize}
\item The \textbf{serial workload} can be estimated as the sum of the number of cells on all levels
$
\mathcal{W}_s = \sum_l \mathcal{C}_l.
$
This metric is based on the assumption that all cells have the same costs, which {is not necessarily true in the context of hanging nodes~\cite{munch2021hn}}.
\item The \textbf{parallel workload} can be estimated as the sum of the maximum number of cells owned by any
process on each level:
$
\mathcal{W}_p = \sum_l \max_p \mathcal{C}_l^p,
$
i.e., the critical path of the cells.
In the ideal case, one would expect that $\mathcal{C}_l^p = \mathcal{C}_l/p$ and therefore 
$\mathcal{W}_p = \mathcal{W}_s/p$. However, due to potential load imbalances, the work {might not be} well
distributed on the levels, i.e.,  
$\max_p{(\mathcal{C}_l^{p})} \ge \mathcal{C}_s/p $. 
Since one can theoretically only proceed to process the next level once all
processes have finished a level, load imbalances will result in some processes waiting at some imaginary
barriers. We say ``imaginary barriers'' as level operators generally do not have any barriers, but only rely on
point-to-point communication between neighboring processes. Nevertheless,
this simplified point of view is acceptable, since multigrid algorithms are
multiplicative Schwarz methods between levels, inherently leading to a serial execution of the levels. 
We define \textbf{parallel workload efficiency} as
$
\mathcal{W}_s / (\mathcal{W}_p \cdot p) 
$, as has been also done in~\citep{clevenger2020flexible}.
\item We define \textbf{horizontal communication efficiency} as 50\% of the number of ghost cells accumulated over all ranks and divided by the total number of cells. {\color{\myblue}The division by two is necessary to take into account that only one neighbor is updating the ghost values.} As such, this ratio can be seen as a proxy of how much information needs to be communicated, when computing residuals and updating ghost values.
As this number counts cells, it is independent of the polynomial degree of the element chosen. The element degree used determines the absolute amount of communication necessary.
{\color{\myred}Note that in reality, only degrees of freedom located at the interface have to be exchanged such that the fraction of the solution that needs to be communicated is less than the fraction of those cells.}
\item \textbf{Vertical communication efficiency} is the share of fine cells that have the same
owning process as their corresponding coarse cell (parent). This quantity gives an indication on the efficiency of
the transfer operator and on how much data has to be sent around. A small number indicates that most of the
data has to be completely permuted, involving a large volume of communication.
This metric has been considered in ~\citep{clevenger2020flexible} as well.
\item Increasing \textbf{number of (multigrid) levels} leads to additional synchronization points and
communication steps and, as a result, might lead to increased latency.
\end{itemize}

Furthermore, \textbf{memory consumption} {of the grid class} is a metric we will consider.\footnote{We use the memory-consumption output provided by \texttt{deal.II}. No particular efforts have been put in reducing the memory consumption of the triangulations in the case of 
global coarsening.} A common argument supporting the 
usage of local-smoothing algorithms is that {\color{\myred}no space is needed for potentially differently partitioned meshes and complex data
structures providing the connectivity between them, since} the multigrid algorithm can simply reuse the 
already existing mesh hierarchy also for the multigrid levels~\citep{clevenger2020flexible}.

\subsection*{Examples}

\begin{figure}
\begin{minipage}{0.33\textwidth}
%
%
\begin{minipage}{0.49\textwidth}
\centering

{\scriptsize $L=4$}

\includegraphics[width=1.0\textwidth]{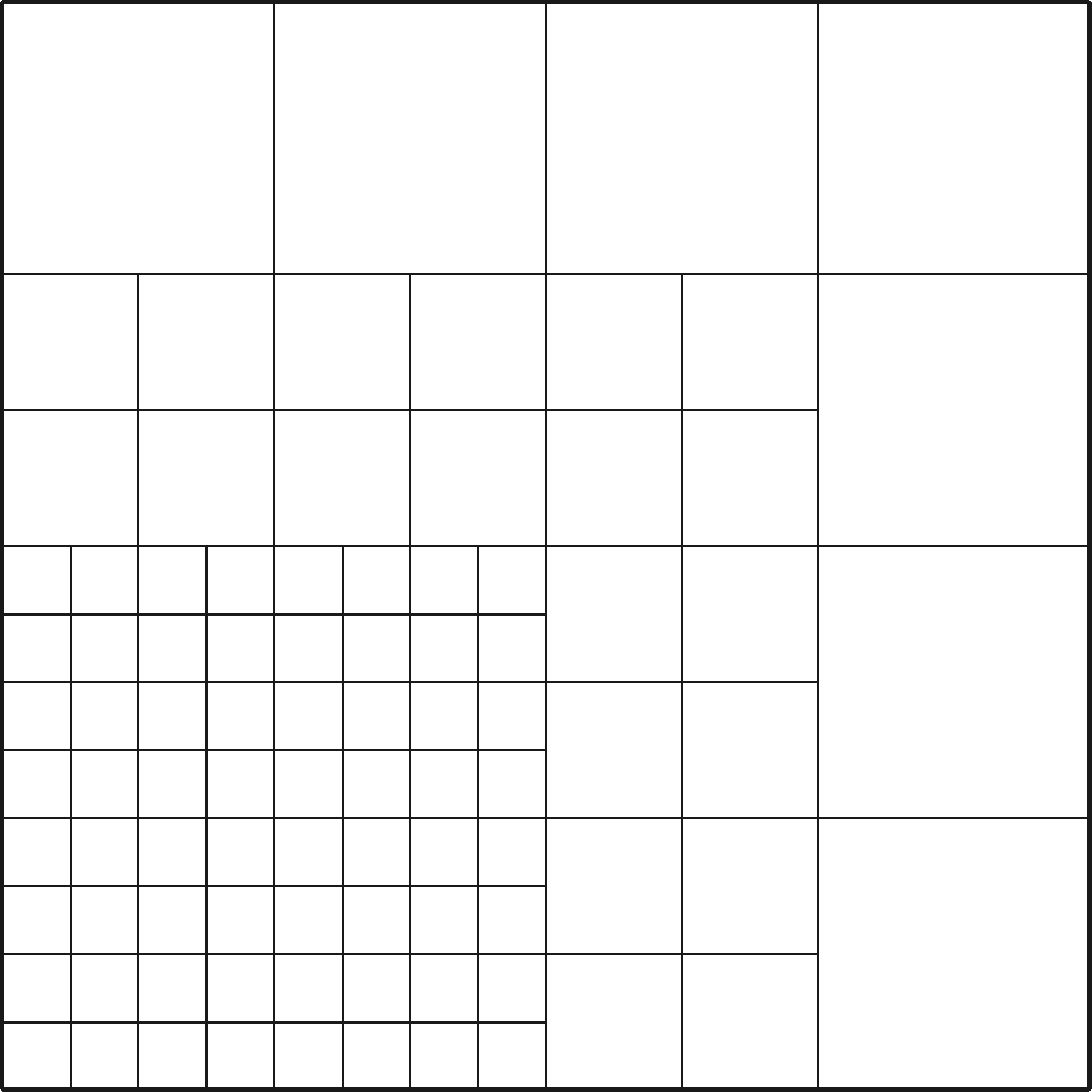}
\end{minipage}
\begin{minipage}{0.49\textwidth}
\centering

{\scriptsize $L=5$}

\includegraphics[width=1.0\textwidth]{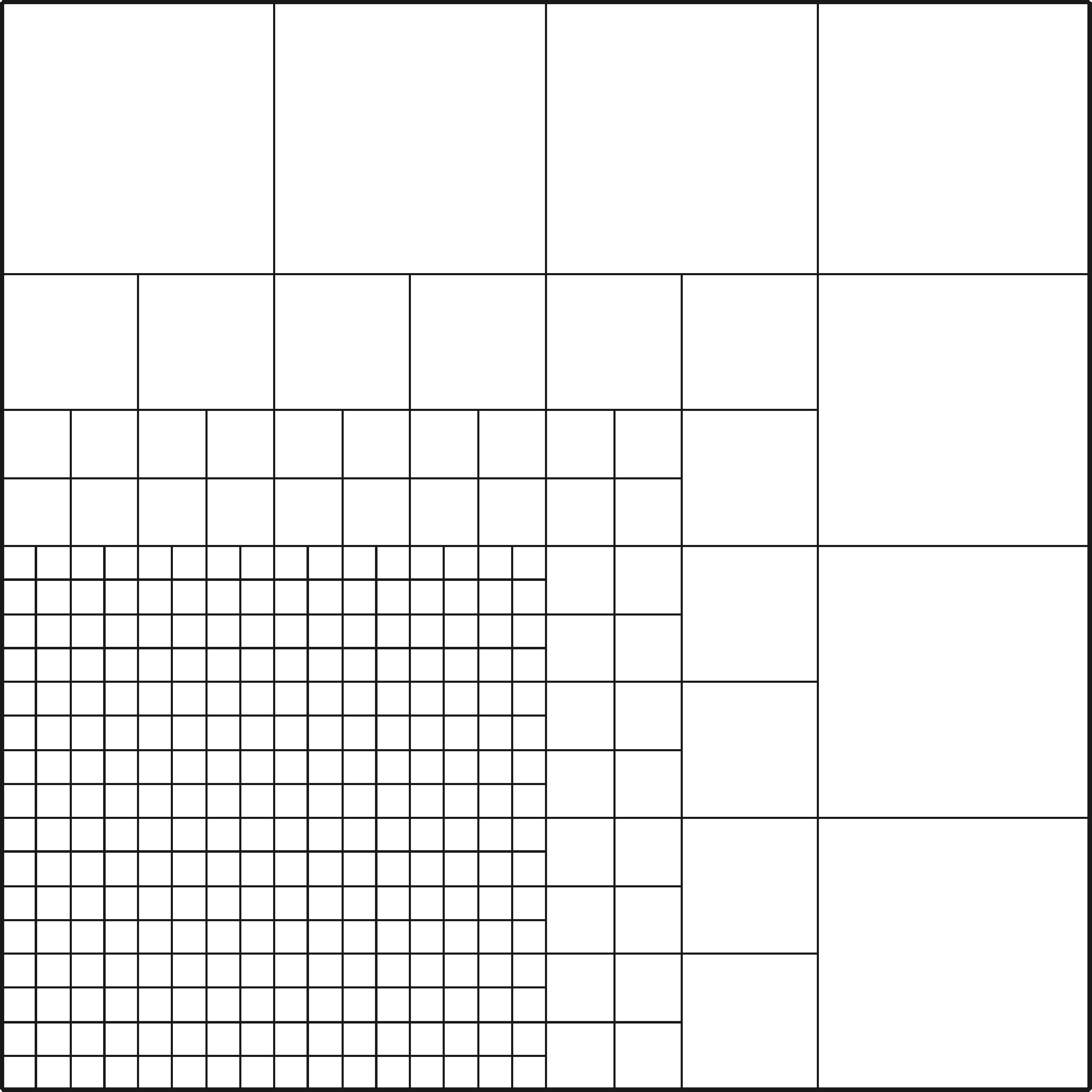}
\end{minipage}

%
%
\begin{minipage}{0.49\textwidth}
\centering

\includegraphics[width=1.0\textwidth]{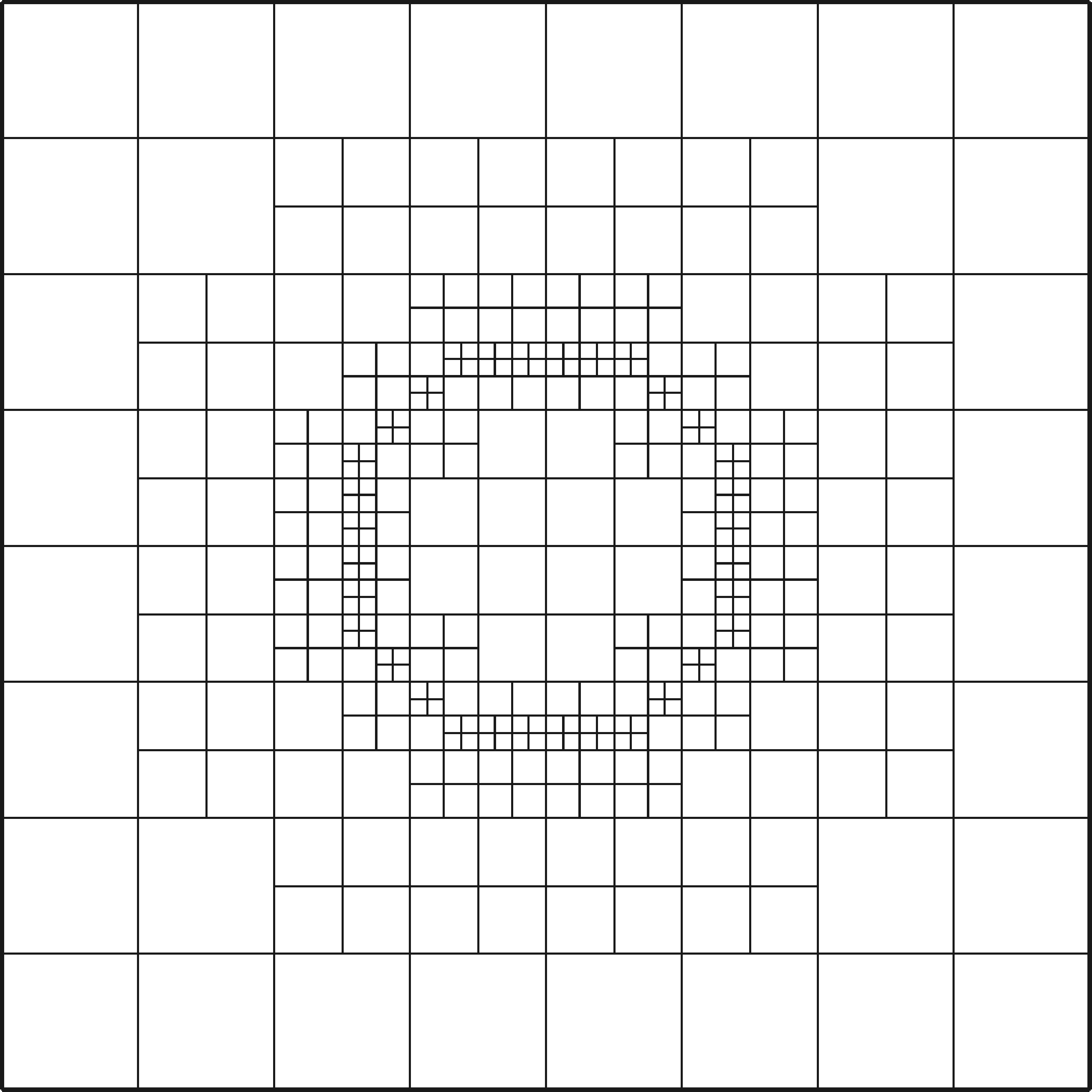}

\vspace{-0.2cm}

{\scriptsize $L=6$}
\end{minipage}
\begin{minipage}{0.49\textwidth}
\centering

\includegraphics[width=1.0\textwidth]{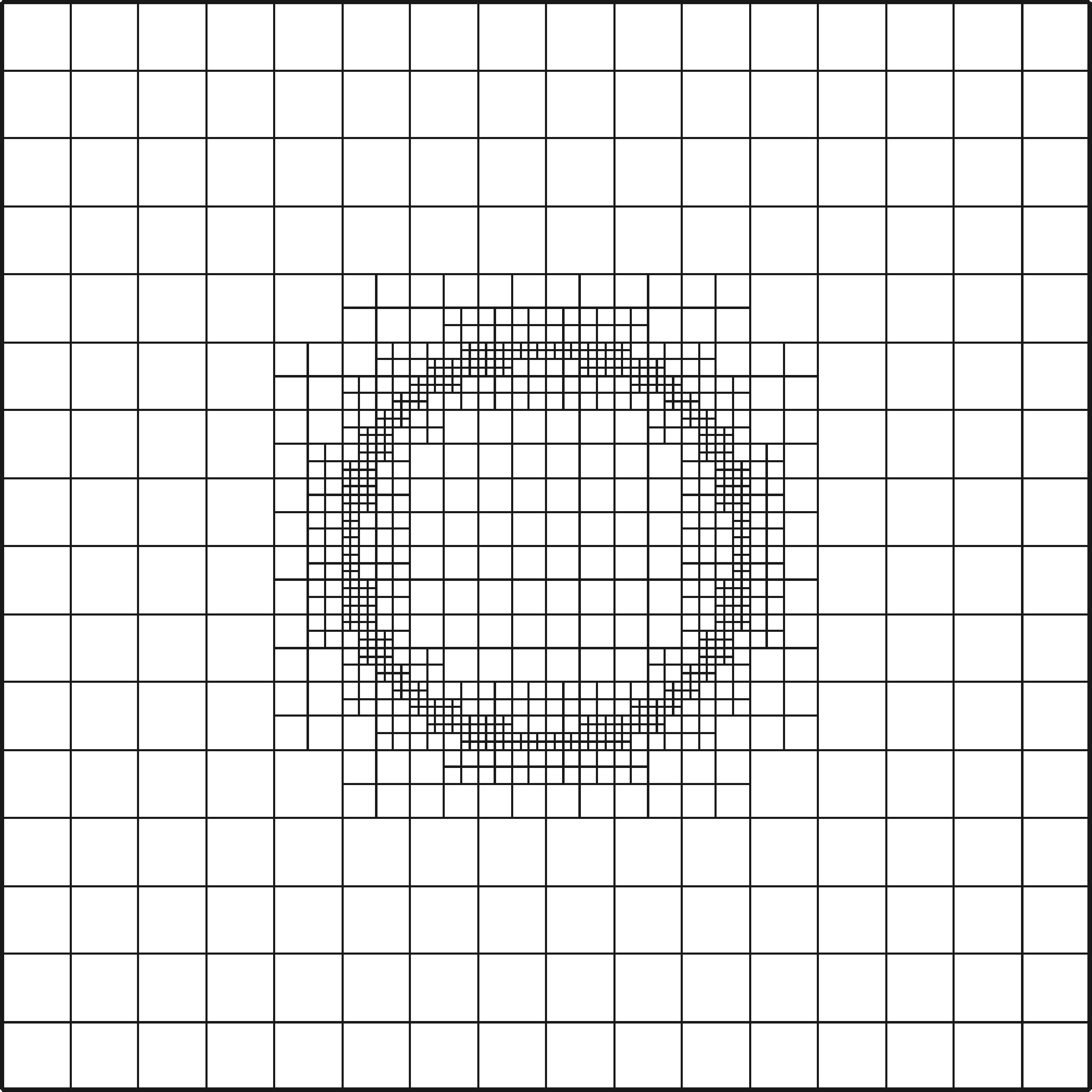}

\vspace{-0.2cm}

{\scriptsize $L=7$}
\end{minipage}

\end{minipage}
\quad
\begin{minipage}{0.59\textwidth}
\begin{scriptsize}
\begin{tabular}{c|cc|cc|cc|cc}
\toprule 
 & \multicolumn{4}{c|}{octant} & \multicolumn{4}{c}{shell} \\
 & \multirow{2}{*}{\#cells} & \multirow{2}{*}{\%HN} & \multicolumn{2}{c|}{\#dofs} & \multirow{2}{*}{\#cells} & \multirow{2}{*}{\%HN} & \multicolumn{2}{c}{\#dofs} \\
$L$ & &  & $p=1$ & $p=4$ &  &  & $p=1$ & $p=4$ \\
\midrule
3&	1.2e+2&31\%&	2.2e+2&	9.3e+3&		-&	-&	-&	- \\
4&	7.0e+2&37\%&	1.0e+3&	5.1e+4&		-&	-&	-&	- \\
5&	4.7e+3&23\%&	5.7e+3&	3.2e+5&		1.2e+3& 69\% &	2.0e+3&	9.3e+4 \\
6&	3.5e+4&12\%&	3.8e+4&	2.3e+6&		6.8e+3& 78\% &	9.8e+3&	5.1e+5 \\
7&	2.7e+5&6.2\%&	2.8e+5&	1.8e+7&		3.7e+4& 70\% &	4.8e+4&	2.6e+6 \\
8&	2.1e+6&3.1\%&	2.2e+6&	1.4e+8&		2.7e+5& 38\% &	3.2e+5&	1.9e+7 \\
9&	1.7e+7&1.6\%&	1.7e+7&	1.1e+9 &		2.2e+6	& 19\% &2.3e+6&	1.4e+8 \\
10&	1.3e+8&0.8\%&1.4e+8&8.6e+9&				1.7e+7&	10\% & 1.8e+7 &1.1e+9 \\	
11&	1.1e+9&0.4\%&1.1e+9&6.9e+10&				1.4e+8&	4.8\% & 1.4e+8 &8.9e+9 \\	
12&	8.6e+9&0.2\%&8.9e+9&-&				1.1e+9&	2.4\% & 1.1e+9 &7.1e+10 \\	
13&	-&-&-&-&				8.9e+9&	1.2\% & 8.9e+9 &5.7e+11 \\	
\bottomrule
\end{tabular}
\end{scriptsize}
\end{minipage}

\caption{Cross section at the center of the geometries of the \texttt{octant} (top) and the \texttt{shell} (bottom) simulation. Additionally, the number of cells,
{\color{\mymagenta}the share of cells with hanging-node constraints,} and the number of DoFs (for a scalar Lagrange element with $p=1$ and $p=4$) are given for each refinement case. Note that we only consider 3D geometries in this publication.}\label{fig:metric:geom}

\end{figure}

In the experimental sections~\ref{sec:performance_h} and~\ref{sec:performance_hp}, we will 
consider two types of static 3D meshes, as has been also done in~\citep{clevenger2020flexible}. They are obtained by
refining a coarse mesh consisting of a single cell defined by $[-1,1]^3$
according to one of  the following two solution criteria:
\begin{itemize}
\item \textbf{octant}: refine all mesh cells in the first octant $[-1,0]^3$ $L$ times and
\item \textbf{shell}: after $L-3$ uniform refinements, perform three local refinement steps 
with all cells whose center $\bm c$ is
$ |\bm c| \le0.55$, $0.3 \le |\bm c| \le0.43$, and $0.335 \le |\bm c| \le0.39$.

\end{itemize}
These two meshes are relevant in practice, since similar meshes occur in simulations of flows
with far fields and of multi-phase flows with bubbles~\cite{KronbichlerDiagneHolmgren2016} or any kind of interfaces.
All the refinement procedures are completed by a closure after each step, ensuring one-irregularity in the
sense that two leaf cells may only differ by one level if they share a vertex.  
Figure~\ref{fig:metric:geom} shows the considered meshes and provides numbers regarding
the cell count for $3\le  L \le 13$. {\color{\myred}All meshes
are partitioned along space-filling curves~\cite{Bangerth2011,Burstedde2011} with the
option to assign cells weights.}

\begin{figure}

\centering

\captionsetup{type=table}
\caption{Geometrical multigrid statistics for the \texttt{octant} test case for different numbers of refinements (wl: serial/parallel workload, wl-eff: parallel workload efficiency, v-eff: vertical communication efficiency; h-eff: horizontal communication efficiency; mem: memory consumption in bytes).}\label{fig:tab:geom:stat:quadrant}

\begin{scriptsize}
\begin{tabular}{c|cc|cc|ccccc|ccccc}
\toprule
& \multicolumn{4}{c|}{1 process} &  \multicolumn{10}{c}{192 processes} \\
 &  \multicolumn{2}{c|}{local smoothing} & \multicolumn{2}{c|}{global coarsening} & \multicolumn{5}{c|}{local smoothing} & \multicolumn{5}{c}{global coarsening} \\ \midrule
$L$ & wl &  mem & wl &  mem &  wl & wl-eff &  v-eff & h-eff &  mem & wl & wl-eff &  v-eff & h-eff &  mem \\
\midrule 
\input{./global-coarsening-paper-data/small-scaling/quadrant_summary.tex}
\bottomrule
\end{tabular}
\end{scriptsize}

\vspace{-0.1cm}

\captionsetup{type=table}
\caption{Geometrical multigrid statistics for the \texttt{shell} test case for different numbers of refinements.}\label{fig:tab:geom:stat:annulus}

\begin{scriptsize}
\begin{tabular}{c|cc|cc|ccccc|ccccc}
\toprule
& \multicolumn{4}{c|}{1 process} &  \multicolumn{10}{c}{192 processes} \\
 &  \multicolumn{2}{c|}{local smoothing} & \multicolumn{2}{c|}{global coarsening} & \multicolumn{5}{c|}{local smoothing} & \multicolumn{5}{c}{global coarsening} \\ \midrule
$L$ & wl &  mem & wl &  mem &  wl & wl-eff &  v-eff & h-eff &  mem & wl & wl-eff &  v-eff & h-eff &  mem \\
\midrule 
\input{./global-coarsening-paper-data/small-scaling/annulus_summary.tex}
\bottomrule
\end{tabular}
\end{scriptsize}

\vspace{0.2cm}

\captionsetup{type=figure}
\input{./global-coarsening-paper-data/small-scaling/workload-0001.tex}\label{fig:geom:stat:level:quadrant}


\input{./global-coarsening-paper-data/small-scaling/workload-0192.tex}\label{fig:geom:stat:level:annulus}

\end{figure}

Tables~\ref{fig:tab:geom:stat:quadrant} and~\ref{fig:tab:geom:stat:annulus} give{---as examples---}evaluated numbers
for geometrical metrics
of the two considered meshes for a single process and for 192 processes {\color{\myred}with
cells constrained by hanging nodes} weighted {with the factor of} 2 for partitioning, compared to the rest of the cells. 
For a single process, only 
workload and memory consumption are shown. 

Starting with the \texttt{octant} case, one can see that the serial workload in the case of
local smoothing and global coarsening is similar, with local smoothing having consistently less work.
The workload of each level is depicted in Figure~\ref{fig:geom:stat:level:quadrant}.
The behavior of the memory consumption is similar to the one of the workload: global coarsening has a slightly higher memory consumption,
since it explicitly  needs to store the coarser meshes as well; however, the second finest mesh has already 
{approximately} one eighth of the size of the coarsest triangulation so that the overhead is {small}.
In the parallel case, the memory consumption differences are higher: 
this is related to the fact that overlapping ghosted forests of trees need to be saved. In contrast,
the workload in the parallel case is much lower in the case of global coarsening: while global coarsening is
able to reach efficiencies higher than $90\%$, the efficiency is only approximately 50-60\% in the case
of local smoothing. The high value in the global-coarsening case was to be expected, since we 
repartition each level during construction. The low value in the case of local smoothing can be
explained by taking a look at Figure~\ref{fig:geom:stat:level:annulus}, where the minimum, maximum, and average workload are shown
for each level. {\color{\myblue}The workload on the finest level is optimally distributed between
processes with cells, however, there are a few processes without any cells, i.e., any work, on that
level. On the second level, most processes can reduce the number of cells nearly optimally by a factor of 8, but processes idle on the finest level start to 
participate in the smoothing process: since they have not participated in the coarsening yet,
the number of cells for those processes is much higher and approximately the same as other
processes had on the finest level. This pattern of discrepancy of minimum and maximum number of cells on the lower levels continues  and, as a consequence, the maximum 
workload per level is higher, increasing the overall parallel workload.}
{Figure~\ref{fig:geom:stat:level:annulus} might also
give the impression that there is a load imbalance in the case of global coarsening, 
since the minimum workload on the finest level
is the half of the maximum value. This is related to the way how we partition---penalizing of
cells that have hanging nodes with a weight of 2---and to the fact that a lot of cells with hanging nodes are
clustered locally. 
 The actually resulting load imbalance is small,
as shown in Figure~\ref{fig:performance_h:profile:parallel}.}
The difference
between local smoothing and global coarsening 
in the parallel workload comes at a price. While the vertical communication {efficiency} is---by construction---high
in the local-smoothing case, this is not true in the case of global coarsening: 20\% and less is not 
uncommon, requiring the permutation of data during transfer. The horizontal communication costs are similar 
in the case of local smoothing and global coarsening. {\color{\myblue}Assuming
that 1) pre- and postsmoothing are the most time-consuming steps, 2)} the workload is the
relevant metric, and 3) the transfer between levels is not dominant  in this case, the values indicate that global coarsening might be twice as fast as local smoothing. 
{\color{\myred}However, if the transfer is the bottleneck,
the picture might look differently so that a conclusive
statement only based on geometrical metrics is not possible
in this case.}

The observations made for the \texttt{octant} case are also valid, but  {\color{\myblue}even more pronounced in the \texttt{shell} case}. Here, the 
vertical communication is more favorable in the case of global coarsening and the workload efficiency is significantly 
worse in the local-smoothing case so that one can expect
a noticeable speedup when using global coarsening. 

In summary, one can state the following: in the serial case, global coarsening has to process at least as many cells as 
local smoothing so that one can expect that the latter has---with the assumption that the number of iterations is the same---an advantage regarding
throughput{\color{\myred}, particularly since no hanging-node
constraints have to be applied}. 
With increasing number of processes, the workload {\color{\myred}might not be well distributed anymore
in the case of local smoothing with the result that parallel efficiency drops} 
strongly. In contrast, this is---by construction---not an issue in the global-coarsening case, 
since the work is simply redistributed between the levels. The price is that one might need
to send around a lot of data during restriction and prolongation. The number of {levels
in the case of local smoothing and of global coarsening is the same}, leading to potentially same
scaling limits $\mathcal{O}(L)$. However, with appropriate partitioning of the levels one could decrease the number of participating processes on the coarser levels and switch to 
a coarse-grid solver at an earlier stage in the global-coarsening case, leading to a better scaling
limit.  In the following, we will
show experimentally that the above statements can be verified and will take a more detailed look
at the significance of a good load balance and of a cheap transfer, which are mutually contradicting requirements.

Making definite general conclusions is difficult as they depend on the 
number of processes as well as on the type of the coarse mesh and of the
refinement. While the two meshes considered here are prototypical for many problems we have encountered in practice, it is clear that the statements {\color{\myred} made in this
publication} cannot hold in all cases, since it is easy to construct examples that favor one over the other.
However, the two meshes demonstrate clearly, which aspects are dominating at which problem sizes; the actual crossover point,
however, might be problem- and hardware-specific.

%% file: global-coarsening-paper-data/small-scaling/quadrant_summary.tex
3  & 1.4e+2 & 6.6e+4 & 1.4e+2 & 8.7e+4 & 1.8e+1 & 3\% & 89\% & 60\% & 1.3e+6 & 2.5e+1 & 3\% & 17\% & 62\% & 2.7e+6 \\
4  & 8.0e+2 & 3.3e+5 & 8.4e+2 & 4.2e+5 & 2.2e+1 & 18\% & 85\% & 56\% & 1.2e+7 & 3.3e+1 & 13\% & 4\% & 58\% & 1.5e+7 \\
5  & 5.4e+3 & 2.0e+6 & 5.6e+3 & 2.5e+6 & 7.2e+1 & 38\% & 88\% & 58\% & 5.3e+7 & 6.7e+1 & 43\% & 1\% & 59\% & 6.8e+7 \\
6  & 4.0e+4 & 1.4e+7 & 4.0e+4 & 1.7e+7 & 4.0e+2 & 51\% & 96\% & 65\% & 1.3e+8 & 2.8e+2 & 76\% & 6\% & 66\% & 2.0e+8 \\
7  & 3.1e+5 & 1.1e+8 & 3.1e+5 & 1.3e+8 & 2.7e+3 & 58\% & 99\% & 75\% & 4.2e+8 & 1.7e+3 & 92\% & 13\% & 75\% & 6.2e+8 \\
8  & 2.4e+6 & 8.6e+8 & 2.4e+6 & 9.9e+8 & 2.0e+4 & 62\% & 99\% & 84\% & 1.8e+9 & 1.3e+4 & 96\% & 23\% & 84\% & 2.4e+9 \\
9  & 1.9e+7 & 6.8e+9 & 1.9e+7 & 7.8e+9 & 1.6e+5 & 64\% & 99\% & 91\% & 1.0e+10 & 1.0e+5 & 98\% & 38\% & 91\% & 1.3e+10 \\

%% file: global-coarsening-paper-data/small-scaling/annulus_summary.tex
5  & 1.4e+3 & 6.0e+5 & 1.8e+3 & 9.5e+5 & 3.1e+1 & 22\% & 80\% & 55\% & 3.2e+7 & 4.5e+1 & 21\% & 2\% & 56\% & 4.7e+7 \\
6  & 7.8e+3 & 3.2e+6 & 9.2e+3 & 4.4e+6 & 1.6e+2 & 26\% & 89\% & 59\% & 7.7e+7 & 1.0e+2 & 47\% & 12\% & 59\% & 1.3e+8 \\
7  & 4.2e+4 & 1.7e+7 & 4.9e+4 & 2.2e+7 & 7.6e+2 & 28\% & 96\% & 66\% & 1.5e+8 & 4.3e+2 & 58\% & 36\% & 65\% & 3.0e+8 \\
8  & 3.1e+5 & 1.2e+8 & 3.4e+5 & 1.4e+8 & 4.7e+3 & 34\% & 99\% & 76\% & 4.4e+8 & 2.3e+3 & 75\% & 78\% & 75\% & 7.7e+8 \\
9  & 2.5e+6 & 8.9e+8 & 2.5e+6 & 1.1e+9 & 3.5e+4 & 36\% & 99\% & 85\% & 1.8e+9 & 1.5e+4 & 86\% & 93\% & 84\% & 2.7e+9 \\
10  & - & - & - & - & 2.7e+5 & 38\% & 99\% & 91\% & 1.0e+10 & 1.1e+5 & 92\% & 97\% & 91\% & 1.3e+10 \\

%% file: global-coarsening-paper-data/small-scaling/workload-0001.tex
\begin{tikzpicture}[thick,scale=0.7, every node/.style={scale=0.75}]
\begin{axis}[
   title={Local smoothing},
   bar width=10pt,
   width=0.5\textwidth,
   height=0.35\textwidth,
   y tick label style={/pgf/number format/fixed, /pgf/number format/fixed zerofill,/pgf/number format/precision=1},
   ylabel={Workload},
   xlabel={Level}, 
   xmajorgrids,xminorgrids,
   ymajorgrids,yminorgrids,
   xmin=-0.5,xmax=8.500000,
   ymin=0.0, ymax=2339907.900000,
   xtick distance=1,legend columns=3, legend pos={north west}
]
\addplot [ybar stacked, blue, fill=blue!40!white] coordinates {
(0, 1.000000)
(1, 8.000000)
(2, 64.000000)
(3, 216.000000)
(4, 1000.000000)
(5, 5832.000000)
(6, 39304.000000)
(7, 287496.000000)
(8, 2097152.000000)
};
\end{axis}
\end{tikzpicture}
\begin{tikzpicture}[thick,scale=0.7, every node/.style={scale=0.75}]
\begin{axis}[
   title={Global coarsening},
   bar width=10pt,
   width=0.5\textwidth,
   height=0.35\textwidth,
   y tick label style={/pgf/number format/fixed, /pgf/number format/fixed zerofill,/pgf/number format/precision=1},
   xlabel={Level}, 
   xmajorgrids,xminorgrids,
   ymajorgrids,yminorgrids,
   xmin=-0.5,xmax=8.500000,
   ymin=0.0, ymax=2339907.900000,
   xtick distance=1,legend columns=3, legend pos={north west}
]
\addplot [ybar stacked, blue, fill=blue!40!white] coordinates {
(0, 1.000000)
(1, 8.000000)
(2, 15.000000)
(3, 120.000000)
(4, 701.000000)
(5, 4712.000000)
(6, 34903.000000)
(7, 269998.000000)
(8, 2127189.000000)
};
\end{axis}
\end{tikzpicture}

\caption{Workload of an octant simulation with a single process.}

%% file: global-coarsening-paper-data/small-scaling/workload-0192.tex
\begin{tikzpicture}[thick,scale=0.7, every node/.style={scale=0.75}]
\begin{axis}[
   title={Local smoothing},
   bar width=10pt,
   width=0.5\textwidth,
   height=0.35\textwidth,
   y tick label style={/pgf/number format/fixed, /pgf/number format/fixed zerofill,/pgf/number format/precision=1},
   ylabel={Workload},
   xlabel={Level}, 
   xmajorgrids,xminorgrids,
   ymajorgrids,yminorgrids,
   xmin=-0.5,xmax=8.500000,
   ymin=0.0, ymax=12540.000000,
   xtick distance=1,legend columns=3, legend pos={north west}
]
\addlegendimage{ybar,ybar legend,blue, fill=blue!40!white}
\addlegendentry{min};
\addlegendimage{ybar,ybar legend,red, fill=red!40!white}
\addlegendentry{max};
\addlegendimage{}
\addlegendentry{avg};
\addplot [ybar stacked, blue, fill=blue!40!white] coordinates {
(0, 0.000000)
(1, 0.000000)
(2, 0.000000)
(3, 0.000000)
(4, 2.000000)
(5, 20.000000)
(6, 160.000000)
(7, 1283.000000)
(8, 0.000000)
};
\addplot [ybar stacked, red, fill=red!40!white]  coordinates {
(0, 1.000000)
(1, 3.000000)
(2, 30.000000)
(3, 62.000000)
(4, 157.000000)
(5, 460.000000)
(6, 1480.000000)
(7, 5149.000000)
(8, 11400.000000)
};
\draw (axis cs:-0.400000,0.000000) -- (axis cs:0.400000,0.000000);
\draw (axis cs:0.600000,0.000000) -- (axis cs:1.400000,0.000000);
\draw (axis cs:1.600000,0.000000) -- (axis cs:2.400000,0.000000);
\draw (axis cs:2.600000,1.000000) -- (axis cs:3.400000,1.000000);
\draw (axis cs:3.600000,5.000000) -- (axis cs:4.400000,5.000000);
\draw (axis cs:4.600000,30.000000) -- (axis cs:5.400000,30.000000);
\draw (axis cs:5.600000,204.000000) -- (axis cs:6.400000,204.000000);
\draw (axis cs:6.600000,1497.000000) -- (axis cs:7.400000,1497.000000);
\draw (axis cs:7.600000,10922.000000) -- (axis cs:8.400000,10922.000000);
\end{axis}
\end{tikzpicture}
\begin{tikzpicture}[thick,scale=0.7, every node/.style={scale=0.75}]
\begin{axis}[
   title={Global coarsening},
   bar width=10pt,
   width=0.5\textwidth,
   height=0.35\textwidth,
   y tick label style={/pgf/number format/fixed, /pgf/number format/fixed zerofill,/pgf/number format/precision=1},
   xlabel={Level}, 
   xmajorgrids,xminorgrids,
   ymajorgrids,yminorgrids,
   xmin=-0.5,xmax=8.500000,
   ymin=0.0, ymax=12540.000000,
   xtick distance=1,legend columns=3, legend pos={north west}
]
\addlegendimage{ybar,ybar legend,blue, fill=blue!40!white}
\addlegendentry{min};
\addlegendimage{ybar,ybar legend,red, fill=red!40!white}
\addlegendentry{max};
\addlegendimage{}
\addlegendentry{avg};
\addplot [ybar stacked, blue, fill=blue!40!white] coordinates {
(0, 0.000000)
(1, 0.000000)
(2, 0.000000)
(3, 0.000000)
(4, 0.000000)
(5, 10.000000)
(6, 105.000000)
(7, 864.000000)
(8, 7351.000000)
};
\addplot [ybar stacked, red, fill=red!40!white]  coordinates {
(0, 1.000000)
(1, 8.000000)
(2, 8.000000)
(3, 8.000000)
(4, 10.000000)
(5, 22.000000)
(6, 87.000000)
(7, 568.000000)
(8, 4049.000000)
};
\draw (axis cs:-0.400000,0.000000) -- (axis cs:0.400000,0.000000);
\draw (axis cs:0.600000,0.000000) -- (axis cs:1.400000,0.000000);
\draw (axis cs:1.600000,0.000000) -- (axis cs:2.400000,0.000000);
\draw (axis cs:2.600000,0.000000) -- (axis cs:3.400000,0.000000);
\draw (axis cs:3.600000,3.000000) -- (axis cs:4.400000,3.000000);
\draw (axis cs:4.600000,24.000000) -- (axis cs:5.400000,24.000000);
\draw (axis cs:5.600000,181.000000) -- (axis cs:6.400000,181.000000);
\draw (axis cs:6.600000,1406.000000) -- (axis cs:7.400000,1406.000000);
\draw (axis cs:7.600000,11079.000000) -- (axis cs:8.400000,11079.000000);
\end{axis}
\end{tikzpicture}

\caption{Workload of an octant simulation with 192 processes.}

%% file: chapters/performance.tex
\section{Performance analysis: h-multigrid}\label{sec:performance_h}

In the following, we solve a {3D Poisson problem} with homogeneous Dirichlet boundary conditions and
{a constant right-hand side} as well as
 compare the performance of local-smoothing and
global-coarsening algorithms for the \texttt{octant} and \texttt{shell} test {cases}, as
introduced in Section~\ref{sec:modelling}. {The results for a more complex
setup with non-homogeneous Dirichlet boundary conditions are shown in Table~\ref{tab:gaussian} in Appendix~\ref{sec:app:performance_h}.}

We will use {a continuous Lagrange finite element,
which is defined as the tensor products of 1D  finite elements with
degree $p$}. For quadrature, we consider the consistent
Gauss--Legendre quadrature rule {with $(p+1)^3$ points}.

We start with investigating the serial performance, proceed with parallel execution
with moderate numbers of processes, and finally analyze the parallel behavior on the large scale (150k processes). {We conclude this section with
investigation of an alternative partitioning scheme for global coarsening.}

In order to obtain the best performance, the experiments are configured in the following way:
\begin{itemize}
\item Cells with hanging-node constraints are weighted by the factor of 2.
\item The conjugate-gradient solver is run until a reduction of the $l_2$-norm of the {unpreconditioned residual by $10^{4}$} is
obtained. {We choose a rather coarse tolerance, since} this is a common value for the solution of time-dependent {\color{\myred}problems, such as the
Navier--Stokes equations,} where good initial guesses can be obtained by projection and extrapolation
 without the
need to {converge multigrid to many digits}. {\color{\myred}Similarly, coarse tolerances also indicate the costs of solving in a full-multigrid scenario with the finest level only correcting against the next coarser one.}
\item The conjugate-gradient solver is preconditioned by a single V-cycle of either a local-smoothing or
a global-coarsening multigrid algorithm.
\item All operations in the multigrid V-cycle are run with single-precision floating-point numbers, while the conjugate-gradient solver is run in {double precision~\cite{ljungkvist2014matrix}}.
\item We use a Chebyshev smoother of degree 3 on all levels.
\item As coarse-grid solver, we use two {V-cycles} of AMG (double-precision, ML~\cite{gee2006ml} {\color{\myred}with parameters shown in Appendix~\ref{sec:app:amg}}).
\end{itemize}
The results of performance studies leading to the decision on the
configuration described above are presented in
{Tables~\ref{fig:appendix:hmg:smoother}-\ref{fig:appendix:hmg:weights:large} in} Appendix~\ref{sec:app:performance_h}. All experiments
have been conducted on the SuperMUC-NG supercomputer. Its compute nodes have 2 sockets (each with 24 cores of Intel Xeon
Skylake)
and the AVX-512 ISA
extension so that 8 doubles or 16 floats can be processed per instruction.  A detailed
specification of the hardware is given in Table~\ref{tab:systems}.
The parallel network is organized into islands of 792 compute nodes each.
The maximum network bandwidth per node
within an island is  100GBit/s=12.5GB/s\footnote{\url{https://doku.lrz.de/
display/PUBLIC/SuperMUC-NG}, retrieved on February 26, 2022.} via a fat-tree network topology. Islands are
connected via a pruned-tree network architecture (pruning factor 1:4).

\begin{table}

\centering

  \caption{Specification of the hardware system used for evaluation. Memory
    bandwidth is according to the STREAM triad benchmark (optimized variant
    without
read for ownership transfer involving two reads and one write), and GFLOP/s are
based
on the theoretical maximum at the AVX-512 frequency. The \texttt{dgemm}
performance is measured for $m=n=k=12{,}000$ with Intel MKL 18.0.2.
We {measured} a frequency of 2.5 GHz with AVX-512 dense code for
the current experiments. The empirical machine balance is computed as the ratio
of measured
\texttt{dgemm} performance and STREAM bandwidth from RAM memory.}
    \label{tab:systems}
  {
	\scriptsize
    \begin{tabular}{lc}
      \hline
      & Intel Skylake   Xeon Platinum 8174 \\
      \hline
      cores       & $2\times 24$ \\
      frequency base (max AVX-512 frequency) & 2.7 GHz \\
      SIMD width    & 512 bit \\
      arithmetic peak (\texttt{dgemm} performance)  & 4147 GFLOP/s (3318 GFLOP/s)\\
      memory interface & DDR4-2666, 12 channels \\
      STREAM memory bandwidth & 205 GB/s\\
      empirical machine balance & 14.3 FLOP/Byte  \\
      L1-/L2-/L3-/MEM size & 32kB (per core)/1MB (per core)/66MB (shared)/96GB(shared) \\
      \hline
      compiler + compiler flags & \texttt{g++}, version 9.1.0, \texttt{-O3 -funroll-loops -march=skylake-avx512}\\
      \hline
    \end{tabular}
  }
\end{table}

For local smoothing and global coarsening, we use different implementations
of the transfer operator from \texttt{deal.II} (see Subsection~\ref{sec:impl:transfer}). In order to demonstrate that they are equivalent
and results shown in the following are indeed related to the defintion of
multigrid levels {\color{\myred}and the resulting
different algorithms}, we present in Table~\ref{fig:appendix:hmg:hypercube} in Appendix~\ref{sec:app:performance_h} a performance comparison {of these implementations for
uniformly refined meshes of a cube}, for which both algorithms {are} equivalent.

{\color{\myred}Note: For global coarsening, we have investigated the possibility
to decrease the number of participating processes and
to {\color{\myblue}switch to the coarse-grid solver earlier}. For our
test problems, we could not see any obvious benefits for the
time to solution so that we do not use these
features of global coarsening, but defer their investigation
to future work.}

\subsection{Serial runs: overview}

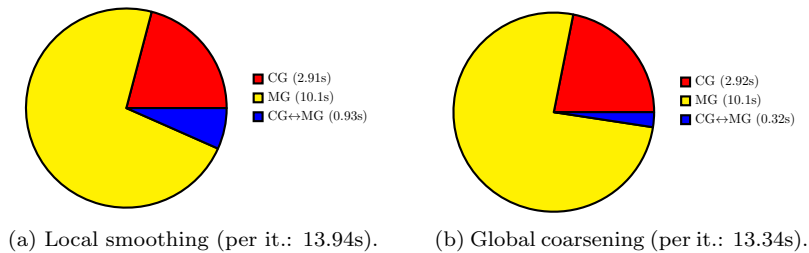
\begin{figure}[!t]
\centering
\subfloat[Local smoothing (per it.: 13.94s).]{
\begin{tikzpicture}[scale=0.44, every node/.style={scale=0.5}]
\pie[sum=auto, hide number, text = legend,color={red, yellow, blue}]
{
2.9139/CG (2.91s),
10.08/MG (10.1s),
0.928/CG$\leftrightarrow$MG (0.93s)
}
\end{tikzpicture}
}
\qquad
\subfloat[Global coarsening (per it.: 13.34s).]{
\begin{tikzpicture}[scale=0.44, every node/.style={scale=0.5}]
\pie[sum=auto, hide number, text = legend, color={red, yellow, blue}]
{
2.9152/CG (2.92s),
10.08/MG (10.1s),
0.322/CG$\leftrightarrow$MG (0.32s)
}
\end{tikzpicture}
}

\caption{\color{\myred}Time per iteration spent for conjugate-gradient solver (CG), multigrid preconditioner (MG), as well as
transfer between solver and preconditioner (CG$\leftrightarrow$MG) in a serial \texttt{shell} simulation with $k=4$ and $L=9$.}\label{fig:timing:overview}

\end{figure}

{\color{\myred}Figure~\ref{fig:timing:overview} gives an overview of the time shares
during the solution process in a serial \texttt{shell} simulation for local smoothing and global coarsening. Without going
into details of the actual numbers, one can see that most of the time is spent in
the multigrid preconditioner in the case both of local smoothing and of
global coarsening (72\%/76\%). It is followed by {\color{\myblue}the other operations in the outer conjugate gradient solver} (21\%/22\%). The
least time is spent for transferring data between preconditioner and
solver (7\%/2\%). It is well visible that more time (factor of approximate 3) is spent for
the transfer in the local-smoothing case. This is not surprising, since the {\color{\myblue}transfer
involves} all multigrid levels sharing cells with the active level.
Since the overwhelming share of solution time is taken by the multigrid preconditioner, {\color{\myblue}all detailed analysis in the remainder of this work concentrates on the multigrid V-cycle}.

One should note that spending 72\%/76\% of the solution
time within the multigrid preconditioner is already low,
particularly taking into account that the outer conjugate-gradient
solver also performs its operator evaluations (1 matrix-vector multiplication per iteration) in an
efficient matrix-free fashion. The low costs are related to
the usage of single-precision floating point numbers and
to the low number of pre-/postsmoothing steps, which results
in a total of 6-7 matrix-vector multiplications per level and iteration.
}

\subsection{Serial run}

{Tables}~\ref{tab:performance_h:quadrant} and~\ref{tab:performance_h:annulus} show the number of
iterations and the time to solution for the \texttt{octant} and the \texttt{shell} test cases run {\color{\myred}serially}
with local smoothing and
global coarsening as well as with the polynomial degrees $p=1$ and $p=4$.

\begin{figure}[!t]

\centering

\captionsetup{type=table}
\caption{Number of iterations and time to solution for local smoothing (LS) and global coarsening (GC) with 1 process
and 192 processes for the \texttt{octant} simulation case.}\label{tab:performance_h:quadrant}

\begin{scriptsize}
\begin{tabular}{c|cc|cc|cc|cc|cc|cc|cc|cccc}
\toprule
& \multicolumn{8}{c|}{1 process} &  \multicolumn{8}{c}{192 processes (4 nodes)} \\
& \multicolumn{4}{c|}{$p=1$} &  \multicolumn{4}{c|}{$p=4$} & \multicolumn{4}{c|}{$p=1$} &  \multicolumn{4}{c}{$p=4$} \\
& \multicolumn{2}{c|}{LS} & \multicolumn{2}{c|}{GC} & \multicolumn{2}{c|}{LS} & \multicolumn{2}{c|}{GC} & \multicolumn{2}{c|}{LS} & \multicolumn{2}{c|}{GC} & \multicolumn{2}{c|}{LS} & \multicolumn{2}{c}{GC} \\ \midrule
$L$ &  \#$i$ & $t$[s] &  \#$i$ & $t$[s] &  \#$i$ & $t$[s] &  \#$i$ & $t$[s] &  \#$i$ & $t$[s] &  \#$i$ & $t$[s] &  \#$i$ & $t$[s] &  \#$i$ & $t$[s] &  \\
\midrule
\input{./global-coarsening-paper-data/small-scaling/quadrant_times.tex}
\bottomrule
\end{tabular}
\end{scriptsize}

\vspace{0.4cm}

%

\captionsetup{type=table}
\caption{Number of iterations and time to solution for local smoothing (LS) and global coarsening (GC) with 1 process
and 192 processes for the \texttt{shell} simulation case.}\label{tab:performance_h:annulus}

\begin{scriptsize}
\begin{tabular}{c|cc|cc|cc|cc|cc|cc|cc|cccc}
\toprule
& \multicolumn{8}{c|}{1 process} &  \multicolumn{8}{c}{192 processes (4 nodes)} \\
& \multicolumn{4}{c|}{$p=1$} &  \multicolumn{4}{c|}{$p=4$} & \multicolumn{4}{c|}{$p=1$} &  \multicolumn{4}{c}{$p=4$} \\
& \multicolumn{2}{c|}{LS} & \multicolumn{2}{c|}{GC} & \multicolumn{2}{c|}{LS} & \multicolumn{2}{c|}{GC} & \multicolumn{2}{c|}{LS} & \multicolumn{2}{c|}{GC} & \multicolumn{2}{c|}{LS} & \multicolumn{2}{c}{GC} \\ \midrule
$L$ &  \#$i$ & $t$[s] &  \#$i$ & $t$[s] &  \#$i$ & $t$[s] &  \#$i$ & $t$[s] &  \#$i$ & $t$[s] &  \#$i$ & $t$[s] &  \#$i$ & $t$[s] &  \#$i$ & $t$[s] &  \\
\midrule
\input{./global-coarsening-paper-data/small-scaling/annulus_times.tex}
\bottomrule
\end{tabular}
\end{scriptsize}

\end{figure}

\begin{figure}[!t]

\centering

\input{./global-coarsening-paper-data/small-scaling/profile-0001.tex}\label{fig:performance_h:profile:serial}

\vspace{0.0cm}

\input{./global-coarsening-paper-data/small-scaling/profile-0192.tex}\label{fig:performance_h:profile:parallel}

\end{figure}

It is well visible that local smoothing has to perform at least as many iterations as global coarsening, with
the difference in iterations limited to 1 in the examples considered. This difference
is not surprising, since global coarsening does at least as much work as local smoothing, by
smoothing over the whole computational domain.
{\color{\myred}Note that global coarsening also benefits from the simple setup of a smooth solution with artificial refinement, see Table~\ref{tab:gaussian} in the appendix for a somewhat more realistic test case.}
As already seen in Tables~\ref{fig:tab:geom:stat:quadrant} and~\ref{fig:tab:geom:stat:annulus}, the serial workload is higher in the case of global coarsening; this is
also visible in the times of a single V-cycle (not shown). Nevertheless, the fewer number of iterations leads
to a smaller time to solution in the case of global coarsening in some instances. Generally,
the costs of a global-coarsening V-cycle are relatively more expensive in the case of the \texttt{shell} simulation with linear elements:
This is not surprising due to
the higher number of cells with hanging-node constraints (see also Figure~\ref{fig:metric:geom}) in the
\texttt{shell} case and {\color{\myred}the higher overhead of linear elements
for application of hanging-node constraints}, as analyzed by~\cite{munch2021hn}.

Figure~\ref{fig:performance_h:profile:serial} shows the distribution of times spent on each multigrid level
and in each multigrid stage for the \texttt{octant} case with $L=8$ and $p=4$.
While the times spent on each multigrid level show similar trends in the
case of local smoothing and global coarsening, distinct (and expected) differences are visible for the
stages:
{\color{\myred}The \texttt{restriction} and \texttt{prolongation} steps
take about the same time in both cases. The \texttt{presmoothing},
\texttt{residual}, and \texttt{postsmoothing} steps are slightly
more expensive in the global-coarsening case, which is
related to the observation that the evaluation of the level operator $\bm A^{(l)}$ (2/1/3-times) is the dominating factor and the presence of
hanging-node constraints makes its evaluation more expensive.
The observation that the residual evaluation is not more expensive
in the local-smoothing case is related to the fact that the
computation of $\bm A_{ES}^{(l)}\bm x_S^{(l)}$ is a side product of the application of $\bm A^{(l)}$. The only visible additional cost of
local smoothing is $\bm A_{SE}^{(l)} \bm x_E^{(l)}$ (\texttt{edge} step). However, its
evaluation is less expensive than the application of $\bm A^{(l)}$,
since only DoFs in proximity to the interface are updated.}

\subsection{Moderately parallel runs}\label{sec:performance_h:moderately}

For the discussion of moderately parallel runs, we have run simulations with 192 processes (on 4 nodes).
Table~\ref{fig:tab:geom:stat:quadrant} and~\ref{fig:tab:geom:stat:annulus} have shown
that imbalance in the
workload leads to parallel workload efficiencies of 40--50\% in the case of local smoothing. The imbalance is lower in the case of global coarsening at the price
of a permutation during prolongation and restriction.

Tables~\ref{tab:performance_h:quadrant} and~\ref{tab:performance_h:annulus} confirm the significance of
{a good workbalance} for the time to solution. The number of processes does not influence the number of iterations due to the
{chosen smoother}.  Speedups---compared to local smoothing---are reached: {\color{\myred}up to 2.3/1.7 for the \texttt{octant} ($p=1$/$p=4$) and up to 2.1/1.6 for the \texttt{shell} case if the different iteration numbers are considered. {\color{\myblue}Normalized per solver iteration, the advantage of global coarsening in these four cases is 2.0, 1.3, 1.7 and 1.6, respectively.}}
Figure~\ref{fig:performance_h:profile:parallel} gives an indication on this behavior by showing the
distribution of the minimum/maximum/average times spent on each multigrid level
and each multigrid stage for the \texttt{octant} case with $L=8$. In the case of global coarsening, it
is well visible that the load is equally distributed and the time spent on the levels is significantly reduced
level by level for the higher levels. For the finest ones, local smoothing shows a completely different picture.
On the finest level, there are processes with hardly any work, but nevertheless the average work is close
to the maximum value, indicating that the load is well-balanced among the processes with work. However, on
the second finest level, a significant workload imbalance---of factor 2.8---is visible also among
the processes with work. This leads to the situation that the maximum time spent on the second finest
level is just slightly less than the one on the finest level, contradicting our expectation of a geometric
{\color{\mygreen}series and leading to the observed increase in the total runtime}.

\subsection{Large-scale parallel run}

\begin{figure}

\centering

\hspace{0.6cm}

\caption{Strong-scaling comparison of local smoothing (LS) and global coarsening (GC) for \texttt{octant} for $p=1$ and $p=4$.
}\label{fig:performance:strongscaling:quadrant}

%

\hspace{0.6cm}%

\caption{Strong-scaling comparison of local smoothing (LS) and global coarsening (GC) for \texttt{shell} for $p=1$ and $p=4$.}\label{fig:performance:strongscaling:annulus}

\end{figure}


{Figures~\ref{fig:performance:strongscaling:quadrant} and~\ref{fig:performance:strongscaling:annulus}
show results of scaling experiments starting with 1 compute node (48 processes) up to 3,072 nodes (147,456 processes).
Besides the times of a single V-cycle, we plot the \textit{normalized throughput} (DoFs per process and time per iteration) against the \textit{time per iteration}.
{\color{\mygreen}The throughput in case of $p=4$ is significantly---by approximatly a factor of 3---higher than for $p=1$. This is expected and related to the
used matrix-free algorithms and their node-level performance, which improves with the polynomial order~\cite{kronbichler2018performance}.}
Just as in the moderately parallel case (see Subsection~\ref{sec:performance_h:moderately}),
we can observe better timings
in the case of global coarsening for a large range of configurations
({\color{\myred}max. speedup: \texttt{octant} 1.9/1.4 for $p=1$/$p=4$, \texttt{shell}: 2.4/2.4}). The number of iterations of local smoothing is 4 for both cases and
all refinement numbers. The number of iterations of global coarsening is 4 for the
\texttt{sphere} case and decreases from 4 to 2 with increasing number of
refinements in the case of \texttt{quadrant} so that the actual speedups reported
above are even higher. {\color{\myred}The high speedup numbers
of global coarsening in the  \texttt{shell} simulation case are particularly
related to its high workload and vertical efficiency, as shown in Table~\ref{tab:performance_h:annulus}.}

The normalized plots give additional insights. Apart from the obvious observations
that, with increasing number of refinements, the minimal time to solution increases (left bottom corner of the plots {\color{\myblue}in Figures~\ref{fig:performance:strongscaling:quadrant} and~\ref{fig:performance:strongscaling:annulus}}) and global coarsening starts with higher throughputs,
one can see that the decrease of parallel efficiency is more moderate in
the case of global coarsening. This is quite astonishing and means, e.g.,
for the \texttt{octant} case with $p=1$/$L=10$, that one can increase the number
of processes by a factor of 16 and still have a throughput per process that
is higher than the one in the case of single-node computations
of local smoothing, which normally shows a
kink in efficiency at early stages (particularly visible in the \texttt{shell} case{\color{\myred}, which matches the findings made in~\cite{clevenger2020flexible}}).
A further observation in the \texttt{shell} case is that
the lines for global coarsening overlap far from the scaling limit, i.e,
the throughput is independent of the number of processes and the number of refinements. This is not the case for local smoothing, where
the throughput deterioriates
with the number of levels, indicating load-balance problems.
The simulations with $p=4$ show similar trends, but the lines are not as smooth, possibly due to the decreased granularity for higher orders.}

%

\subsection{\color{\myred}First-child policy as alternative partitioning strategy for global coarsening}

{In Subsection~\ref{sec:performance_h:moderately}, we have discussed that local smoothing with first-child policy might suffer
from deteriorated reduction rates of the maximimum number of cells on each level; in particular,
there might be processes without any cells, i.e., any work, increasing the critical
path, although the vertical efficiency is optimal.
In this section, we consider the first-child policy, which we use in the context of local smoothing as an alternative for partitioning of the levels
for global coarsening.

Figure~\ref{fig:gmg:fcp:quadrant} shows the timings of large-scale \texttt{octant} simulations for 1) local smoothing, 2) global coarsening with default partitioning, and 3)
global coarsening with first-child policy for $p=1$/$p=4$. The timings of the latter
show similar trends as global coarsening with default partitioning and are
lower than the ones of local smoothing. To explain this counterintuitive observation,
Figure~\ref{fig:gmg:fcp:quadrant} shows the maximum number of cells on each multigrid level of the three approaches for $L=11$  on 256 nodes.
Since global coarsening with first-child policy does not perform any
repartitioning, {\color{\myred}better reduction trends as in the
local-smoothing case on all levels} can not be expected, however, one can observe that, for the local
section of the refinement tree, the number of cells on the finest
levels can be reduced nearly as well as in the case of the default partitioner{\color{\myred}, i.e., the parallel workload and the parallel efficiency are not much {\color{\myblue}worse, nonetheless, with} higher
vertical efficiency}.
This is not possible for the lower levels, but the behavior of the finest levels dominates the overall trends,
since they take the largest time of the
computation.
One should note that local smoothing has access to the same cells, but
simply skips them during smoothing of a given level, missing the opportunity
to reduce the problem size locally. In our experiments, it is not clear whether {\color{\myred}an optimal reduction of cells
is crucial for all configurations}: for $p=1$, global coarsening
with default partitioning is faster for most configurations (up to 20\%) and for $p=4$, global coarsening with first-child policy is faster (up to 10\%).  We could trace this difference back to the different costs of the transfer{\color{\myblue}, where
for $p=4$ a reduced data transfer (both within the compute node and across the network) and for $p=1$ a better load balance is
beneficial.}

In the \texttt{shell} case (Figure~\ref{fig:gmg:fcp:annulus}), we
observe that both global coarsening partitioning strategies result in
very similar reduction rates, which can be traced back to the fact that both strategies
lead to comparable partitionings of the levels {\color{\myred}(see also Table~\ref{tab:performance_h:annulus},
which shows a very high vertical efficiency of the default
partitioning) and---as a consequence---to very similar solution times ($\pm$10\%)}. Local smoothing only reduces
the maximum number of cells per level optimally once all locally refined cells
have been processed and the levels that have been constructed via global
refinement have been reached.
{\color{\myred}This case stresses the issue of load imbalances
related to reduction rates significantly differing between
processes.}

}

\begin{figure}

%

\caption{\color{\myred}Left: Maximum number of cells per level ($L=11$, 256 nodes) and right: strong scaling of global coarsening with first-child
policy (fcp) in comparison to local smoothing (LS) and global coarsening (GC)
with default policy for \texttt{octant}.}\label{fig:gmg:fcp:quadrant}

\hspace{-1.5cm}
%

\caption{
\color{\myred}Left: Maximum number of cells per level ($L=12$, 256 nodes) and right: strong scaling of global coarsening with first-child
policy (fcp) in comparison to local smoothing (LS) and global coarsening (GC)
with default policy for \texttt{shell}.}\label{fig:gmg:fcp:annulus}

\end{figure}

\section{Performance analysis: p-multigrid}\label{sec:performance_hp}

\begin{figure}[!t]
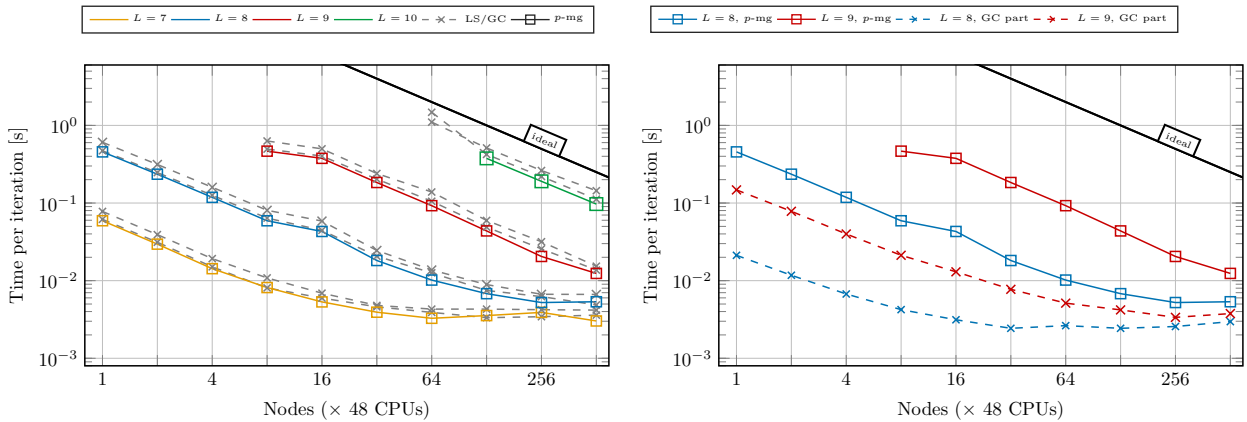


\centering

\hspace{0.3cm}%

\caption{\color{\myred}Strong scaling of $p$-multigrid for \texttt{octant} for $p=4$. $p$-multigrid switches to a global-coarsening
coarse-grid solver immediately once linear elements have been reached. Left: Comparison with local smoothing (LS) and
global coarsening (GC) for $p=4$. Right: Comparison with
global coarsening (GC) for $p=1$ (coarse-grid problem).
}\label{fig:performance_hp:scaling}

\end{figure}

In this section, we consider $p$-multigrid.
The settings {are as} described in Section~\ref{sec:performance_h}.
On the coarsest level ($p=1$, {\color{\myred}fine mesh with
hanging nodes}), we run a single V-cycle of either a
local-smoothing or a global-coarsening geometric multigrid solver in an $hp$ context.
In Appendix~\ref{sec:app:performance_hp}, we present a comparison
with state-of-the-art AMG solvers~\cite{gee2006ml,falgout2006design}
as coarse-grid solvers of $p$-multigrid on 16 nodes. The results show better timings in favor of
geometric multigrid, which also turned out to be more robust with a single V-cycle.

Figure~\ref{fig:performance_hp:scaling} presents a strong-scaling comparison of $h$- and $p$-multigrid versions of the global-coarsening algorithm for the \texttt{octant} case for $p = 4$. The measurements for local smoothing as coarse-grid solver
are skipped here, {since the trends are similar and the values for local smoothing are only {\color{\myred}a few percentage} higher than the ones for global coarsening}.\footnote{The raw data for $L=9$ is
provided in Appendix~\ref{sec:app:performance_hp}.} As we use
a bisection strategy in the context of $p$-multigrid, the overall multigrid algorithm has two additional levels
compared to pure $h$-multigrid (the same fine mesh, but with $p=1$ and $p=2$). In our experiments,
we have observed that $p$-multigrid needs at least as many iterations as the pure (global-coarsening)
$h$-multigrid algorithm (with {\color{\myblue}small differences}---max 1).
One can see that {one cycle of} $p$-multigrid is {\color{\myred}10-15\%} faster than
$h$-multigrid for moderate numbers of processes. {The reason for this is that the smoother application on the finest level is equally expensive, however, the transfer} between the two
finest levels is cheaper: due to the same partitioning of these levels, the data is mainly transferred between
cells that are locally owned on both the coarse and the fine level. At the scaling limit (not shown), $p$-multigrid falls
behind $h$-multigrid regarding performance. This is related to the increased number of levels.

{\color{\myred}For the sake of completeness, Figure~\ref{fig:performance_hp:scaling} also shows global-coarsening timings for $p=1$ for $L=8$/$L=9$ exemplarily,
since it is the coarse-grid problem of the $p$-multigrid
solver. It is well visible that the coarse-grid problem
is negligible  for a wide range of nodes, but becomes noticeable
at the scaling limit.}

%% file: global-coarsening-paper-data/small-scaling/quadrant_times.tex
3  & 4 & 5.0e-4 & 4 & 6.0e-4 & 4 & 3.8e-3 & 4 & 4.3e-3 & 4 & 1.2e-3 & 4 & 1.0e-3 & 4 & 2.5e-3 & 4 & 2.5e-3 \\
4  & 4 & 1.8e-3 & 4 & 2.5e-3 & 4 & 1.8e-2 & 4 & 2.0e-2 & 4 & 2.6e-3 & 4 & 1.7e-3 & 4 & 4.3e-3 & 4 & 4.1e-3 \\
5  & 4 & 8.6e-3 & 4 & 1.2e-2 & 4 & 1.1e-1 & 3 & 8.7e-2 & 4 & 5.3e-3 & 4 & 2.6e-3 & 4 & 6.9e-3 & 3 & 5.0e-3 \\
6  & 4 & 5.1e-2 & 4 & 6.5e-2 & 4 & 8.8e-1 & 3 & 6.5e-1 & 4 & 5.6e-3 & 4 & 4.1e-3 & 4 & 1.7e-2 & 3 & 1.0e-2 \\
7  & 4 & 3.6e-1 & 3 & 3.2e-1 & 4 & 6.9e+0 & 3 & 5.0e+0 & 4 & 1.3e-2 & 3 & 5.7e-3 & 4 & 8.4e-2 & 3 & 5.0e-2 \\
8  & 4 & 2.8e+0 & 3 & 2.3e+0 & 4 & 5.3e+1 & 3 & 3.8e+1 & 4 & 3.9e-2 & 3 & 2.1e-2 & 4 & 7.2e-1 & 3 & 4.3e-1 \\
9  & 4 & 2.2e+1 & 3 & 1.8e+1 & - & - & - & - & 4 & 2.3e-1 & 3 & 1.3e-1 & - & - & - & - \\

%% file: global-coarsening-paper-data/small-scaling/annulus_times.tex
5  & 5 & 3.6e-3 & 4 & 7.1e-3 & 4 & 3.1e-2 & 4 & 4.5e-2 & 5 & 6.6e-3 & 4 & 3.1e-3 & 4 & 6.4e-3 & 4 & 6.3e-3 \\
6  & 5 & 1.5e-2 & 4 & 2.8e-2 & 4 & 1.9e-1 & 4 & 2.2e-1 & 5 & 6.8e-3 & 4 & 4.1e-3 & 4 & 1.0e-2 & 4 & 9.2e-3 \\
7  & 5 & 7.3e-2 & 4 & 1.3e-1 & 4 & 1.0e+0 & 4 & 1.2e+0 & 5 & 9.2e-3 & 4 & 5.6e-3 & 4 & 2.3e-2 & 4 & 1.8e-2 \\
8  & 5 & 5.1e-1 & 4 & 7.3e-1 & 4 & 7.8e+0 & 4 & 7.9e+0 & 5 & 1.7e-2 & 4 & 1.2e-2 & 4 & 1.1e-1 & 4 & 6.9e-2 \\
9  & 5 & 3.7e+0 & 4 & 4.3e+0 & 4 & 5.8e+1 & 4 & 5.6e+1 & 5 & 6.3e-2 & 4 & 3.5e-2 & 4 & 8.5e-1 & 4 & 5.2e-1 \\
10  & - & - & - & - & - & - & - & - & 5 & 3.9e-1 & 4 & 1.9e-1 & - & - & - & - \\

%% file: global-coarsening-paper-data/small-scaling/profile-0001.tex
\begin{tikzpicture}[thick,scale=0.7, every node/.style={scale=0.75}]
\begin{axis}[
   title={Local smoothing},
   bar width=10pt,
   width=0.5\textwidth,
   height=0.35\textwidth,
   ylabel={Exclusive time [s]},
   xlabel={Level}, 
   xmajorgrids,xminorgrids,
   ymajorgrids,yminorgrids,
   xmin=-0.5,xmax=8.500000,
   ymin=0.0, ymax=9.139018,
   y tick label style={/pgf/number format/fixed, /pgf/number format/fixed zerofill,/pgf/number format/precision=2},scaled y ticks=false,
   xtick distance=1,legend columns=3, legend pos={north west}
]
\addplot [ybar stacked, blue, fill=blue!40!white] coordinates {
(0, 0.000198)
(1, 0.000097)
(2, 0.000330)
(3, 0.001006)
(4, 0.004213)
(5, 0.023247)
(6, 0.160459)
(7, 1.194836)
(8, 8.308198)
};
\end{axis}
\end{tikzpicture}
\begin{tikzpicture}[thick,scale=0.7, every node/.style={scale=0.75}]
\begin{axis}[
   title={Global coarsening},
   bar width=10pt,
   width=0.5\textwidth,
   height=0.35\textwidth,
   xlabel={Level}, 
   xmajorgrids,xminorgrids,
   ymajorgrids,yminorgrids,
   xmin=-0.5,xmax=8.500000,
   ymin=0.0, ymax=9.139018,
   y tick label style={/pgf/number format/fixed, /pgf/number format/fixed zerofill,/pgf/number format/precision=2},scaled y ticks=false,
   xtick distance=1,legend columns=3, legend pos={north west}
]
\addplot [ybar stacked, blue, fill=blue!40!white] coordinates {
(0, 0.000197)
(1, 0.000103)
(2, 0.000153)
(3, 0.000655)
(4, 0.003287)
(5, 0.019754)
(6, 0.141565)
(7, 1.108810)
(8, 8.293954)
};
\end{axis}
\end{tikzpicture}

\begin{tikzpicture}[thick,scale=0.7, every node/.style={scale=0.75}]
\begin{axis}[
   bar width=10pt,
   width=0.5\textwidth,
   height=0.35\textwidth,
   y tick label style={/pgf/number format/fixed, /pgf/number format/fixed zerofill,/pgf/number format/precision=2},scaled y ticks=false,
   ylabel={Time [s]},
   xmajorgrids,xminorgrids,
   ymajorgrids,yminorgrids,
   ymin=0.0, ymax=4.699626,
   xtick = {-1, 0, 1,2,3,4,5,6, 7},
   x tick label style={rotate=45,anchor=east},legend columns=3, legend pos={north west},
   xticklabels={,presmoothing, residual, restriction, coarse-grid solver, prolongation, edge, postsmoothing}
]
\addplot [ybar stacked, blue, fill=blue!40!white] coordinates {
(0, 2.875589)
(1, 1.236280)
(2, 0.530589)
(3, 0.000198)
(4, 0.424069)
(5, 0.522821)
(6, 4.103039)
};
\end{axis}
\end{tikzpicture}
\begin{tikzpicture}[thick,scale=0.7, every node/.style={scale=0.75}]
\begin{axis}[
   bar width=10pt,
   width=0.5\textwidth,
   height=0.35\textwidth,
   y tick label style={/pgf/number format/fixed, /pgf/number format/fixed zerofill,/pgf/number format/precision=2},scaled y ticks=false,
   xmajorgrids,xminorgrids,
   ymajorgrids,yminorgrids,
   ymin=0.0, ymax=4.699626,
   xtick = {-1, 0, 1,2,3,4,5,6, 7},
   x tick label style={rotate=45,anchor=east},legend columns=3, legend pos={north west},
   xticklabels={,presmoothing, residual, restriction, coarse-grid solver, prolongation, edge, postsmoothing}
]
\addplot [ybar stacked, blue, fill=blue!40!white] coordinates {
(0, 2.979526)
(1, 1.298900)
(2, 0.540270)
(3, 0.000197)
(4, 0.477192)
(5, 0.000005)
(6, 4.272387)
};
\end{axis}
\end{tikzpicture}

\vspace{-0.5cm}

\caption{Profile of a V-cycle of an octant simulation with a single process for $L=8$ and $p=4$}

%% file: global-coarsening-paper-data/small-scaling/profile-0192.tex
\begin{tikzpicture}[thick,scale=0.7, every node/.style={scale=0.75}]
\begin{axis}[
   title={Local smoothing},
   bar width=10pt,
   width=0.5\textwidth,
   height=0.35\textwidth,
   ylabel={Exclusive time [s]},
   xlabel={Level}, 
   xmajorgrids,xminorgrids,
   ymajorgrids,yminorgrids,
   xmin=-0.5,xmax=8.500000,
   ymin=0.0, ymax=0.100780,
   y tick label style={/pgf/number format/fixed, /pgf/number format/fixed zerofill,/pgf/number format/precision=2},scaled y ticks=false,
   xtick distance=1,legend columns=3, legend pos={north west}
]
\addlegendimage{ybar,ybar legend,blue, fill=blue!40!white}
\addlegendentry{min};
\addlegendimage{ybar,ybar legend,red, fill=red!40!white}
\addlegendentry{max};
\addlegendimage{}
\addlegendentry{avg};
\addplot [ybar stacked, blue, fill=blue!40!white] coordinates {
(0, 0.000101)
(1, 0.000014)
(2, 0.000009)
(3, 0.000012)
(4, 0.000935)
(5, 0.002216)
(6, 0.007881)
(7, 0.025080)
(8, 0.007599)
};
\addplot [ybar stacked, red, fill=red!40!white] coordinates {
(0, 0.000690)
(1, 0.000204)
(2, 0.000379)
(3, 0.000740)
(4, 0.000581)
(5, 0.000759)
(6, 0.002613)
(7, 0.042367)
(8, 0.078828)
};
\draw (axis cs:-0.400000,0.000508) -- (axis cs:0.400000,0.000508);
\draw (axis cs:0.600000,0.000022) -- (axis cs:1.400000,0.000022);
\draw (axis cs:1.600000,0.000032) -- (axis cs:2.400000,0.000032);
\draw (axis cs:2.600000,0.000248) -- (axis cs:3.400000,0.000248);
\draw (axis cs:3.600000,0.001278) -- (axis cs:4.400000,0.001278);
\draw (axis cs:4.600000,0.002415) -- (axis cs:5.400000,0.002415);
\draw (axis cs:5.600000,0.008296) -- (axis cs:6.400000,0.008296);
\draw (axis cs:6.600000,0.033938) -- (axis cs:7.400000,0.033938);
\draw (axis cs:7.600000,0.077684) -- (axis cs:8.400000,0.077684);
\end{axis}
\end{tikzpicture}
\begin{tikzpicture}[thick,scale=0.7, every node/.style={scale=0.75}]
\begin{axis}[
   title={Global coarsening},
   bar width=10pt,
   width=0.5\textwidth,
   height=0.35\textwidth,
   xlabel={Level}, 
   xmajorgrids,xminorgrids,
   ymajorgrids,yminorgrids,
   xmin=-0.5,xmax=8.500000,
   ymin=0.0, ymax=0.100780,
   y tick label style={/pgf/number format/fixed, /pgf/number format/fixed zerofill,/pgf/number format/precision=2},scaled y ticks=false,
   xtick distance=1,legend columns=3, legend pos={north west}
]
\addlegendimage{ybar,ybar legend,blue, fill=blue!40!white}
\addlegendentry{min};
\addlegendimage{ybar,ybar legend,red, fill=red!40!white}
\addlegendentry{max};
\addlegendimage{}
\addlegendentry{avg};
\addplot [ybar stacked, blue, fill=blue!40!white] coordinates {
(0, 0.000111)
(1, 0.000010)
(2, 0.000009)
(3, 0.000008)
(4, 0.000009)
(5, 0.000551)
(6, 0.001301)
(7, 0.008596)
(8, 0.075363)
};
\addplot [ybar stacked, red, fill=red!40!white] coordinates {
(0, 0.000436)
(1, 0.000089)
(2, 0.000333)
(3, 0.000434)
(4, 0.000700)
(5, 0.000607)
(6, 0.000771)
(7, 0.006197)
(8, 0.016255)
};
\draw (axis cs:-0.400000,0.000411) -- (axis cs:0.400000,0.000411);
\draw (axis cs:0.600000,0.000012) -- (axis cs:1.400000,0.000012);
\draw (axis cs:1.600000,0.000023) -- (axis cs:2.400000,0.000023);
\draw (axis cs:2.600000,0.000034) -- (axis cs:3.400000,0.000034);
\draw (axis cs:3.600000,0.000331) -- (axis cs:4.400000,0.000331);
\draw (axis cs:4.600000,0.000830) -- (axis cs:5.400000,0.000830);
\draw (axis cs:5.600000,0.001514) -- (axis cs:6.400000,0.001514);
\draw (axis cs:6.600000,0.009943) -- (axis cs:7.400000,0.009943);
\draw (axis cs:7.600000,0.087681) -- (axis cs:8.400000,0.087681);
\end{axis}
\end{tikzpicture}

\begin{tikzpicture}[thick,scale=0.7, every node/.style={scale=0.75}]
\begin{axis}[
   bar width=10pt,
   width=0.5\textwidth,
   height=0.35\textwidth,
   y tick label style={/pgf/number format/fixed, /pgf/number format/fixed zerofill,/pgf/number format/precision=2},scaled y ticks=false,
   ylabel={Time [s]},
   xmajorgrids,xminorgrids,
   ymajorgrids,yminorgrids,
   ymin=0.0, ymax=0.057936,
   xtick = {-1, 0, 1,2,3,4,5,6, 7},
   x tick label style={rotate=45,anchor=east},legend columns=3, legend pos={north west},
   xticklabels={,presmoothing, residual, restriction, coarse-grid solver, prolongation, edge, postsmoothing}
]
\addlegendimage{ybar,ybar legend,blue, fill=blue!40!white}
\addlegendentry{min};
\addlegendimage{ybar,ybar legend,red, fill=red!40!white}
\addlegendentry{max};
\addlegendimage{}
\addlegendentry{avg};
\addplot [ybar stacked, blue, fill=blue!40!white] coordinates {
(0, 0.028942)
(1, 0.004884)
(2, 0.006802)
(3, 0.000101)
(4, 0.003989)
(5, 0.005364)
(6, 0.018673)
};
\addplot [ybar stacked, red, fill=red!40!white] coordinates {
(0, 0.016093)
(1, 0.011642)
(2, 0.006372)
(3, 0.000690)
(4, 0.007633)
(5, 0.008480)
(6, 0.033996)
};
\draw (axis cs:-0.300000,0.033958) -- (axis cs:0.300000,0.033958);
\draw (axis cs:0.700000,0.013899) -- (axis cs:1.300000,0.013899);
\draw (axis cs:1.700000,0.009543) -- (axis cs:2.300000,0.009543);
\draw (axis cs:2.700000,0.000508) -- (axis cs:3.300000,0.000508);
\draw (axis cs:3.700000,0.008984) -- (axis cs:4.300000,0.008984);
\draw (axis cs:4.700000,0.010407) -- (axis cs:5.300000,0.010407);
\draw (axis cs:5.700000,0.047123) -- (axis cs:6.300000,0.047123);
\end{axis}
\end{tikzpicture}
\begin{tikzpicture}[thick,scale=0.7, every node/.style={scale=0.75}]
\begin{axis}[
   bar width=10pt,
   width=0.5\textwidth,
   height=0.35\textwidth,
   y tick label style={/pgf/number format/fixed, /pgf/number format/fixed zerofill,/pgf/number format/precision=2},scaled y ticks=false,
   xmajorgrids,xminorgrids,
   ymajorgrids,yminorgrids,
   ymin=0.0, ymax=0.057936,
   xtick = {-1, 0, 1,2,3,4,5,6, 7},
   x tick label style={rotate=45,anchor=east},legend columns=3, legend pos={north west},
   xticklabels={,presmoothing, residual, restriction, coarse-grid solver, prolongation, edge, postsmoothing}
]
\addlegendimage{ybar,ybar legend,blue, fill=blue!40!white}
\addlegendentry{min};
\addlegendimage{ybar,ybar legend,red, fill=red!40!white}
\addlegendentry{max};
\addlegendimage{}
\addlegendentry{avg};
\addplot [ybar stacked, blue, fill=blue!40!white] coordinates {
(0, 0.023702)
(1, 0.008242)
(2, 0.008568)
(3, 0.000111)
(4, 0.005723)
(5, 0.000003)
(6, 0.030516)
};
\addplot [ybar stacked, red, fill=red!40!white] coordinates {
(0, 0.005979)
(1, 0.004050)
(2, 0.006219)
(3, 0.000436)
(4, 0.011637)
(5, 0.000003)
(6, 0.014925)
};
\draw (axis cs:-0.300000,0.027279) -- (axis cs:0.300000,0.027279);
\draw (axis cs:0.700000,0.010422) -- (axis cs:1.300000,0.010422);
\draw (axis cs:1.700000,0.011425) -- (axis cs:2.300000,0.011425);
\draw (axis cs:2.700000,0.000411) -- (axis cs:3.300000,0.000411);
\draw (axis cs:3.700000,0.013039) -- (axis cs:4.300000,0.013039);
\draw (axis cs:4.700000,0.000005) -- (axis cs:5.300000,0.000005);
\draw (axis cs:5.700000,0.038199) -- (axis cs:6.300000,0.038199);
\end{axis}
\end{tikzpicture}

\vspace{-0.5cm}

\caption{Profile of a V-cycle of an octant simulation with 192 processes for $L=8$ and $p=4$}

%% file: chapters/stokes.tex
\section{Application: variable viscosity Stokes flow}\label{sec:stokes}

We conclude this publication by presenting preliminary results of a practical application from Geosciences
by integrating the global-coarsening framework into the mantle convection code \texttt{ASPECT}~\cite{KHB12,heister_aspect_methods2} and compare it against the existing local-smoothing implementation~\cite{clevenger2021comparison}.

We consider the variable-viscosity Stokes problem
\begin{align*}
   -\nabla \cdot \left(2\eta \varepsilon(\mathbf u)
                \right)
                + \nabla p &=
  f
  \\
  \nabla \cdot \mathbf u &= 0
\end{align*}
with a Q2-Q1 Taylor-Hood discretization of velocity $\vec{u}$, pressure $p$, and viscosity $\eta(x)$. The resulting linear system
\[
  \begin{pmatrix}
        A & B^T \\ B  & 0
   \end{pmatrix}
   \begin{pmatrix}
        U \\ P
   \end{pmatrix}
   =
   \begin{pmatrix}
        F \\ 0
   \end{pmatrix}
\]
is solved with a Krylov method (in our tests using IDR(2), see \cite{clevenger2021comparison}) preconditioned using a block preconditioner
\[
P^{-1} =
   \begin{pmatrix}
        A & B^T \\
        0 & -S
   \end{pmatrix} ^{-1}
\]
where the Schur complement $S=B A^{-1} B^T$ is approximated using a mass matrix weighted by the viscosity.
The inverses of the diagonal blocks of $A$ and the Schur complement approximation $\hat{S}$ are each approximated by applying a single V-cycle of geometric multigrid using a Chebyshev smoother of degree 4 and implemented in a matrix-free fashion. Each IDR(2) iteration consists of 3 matrix-vector products and 3 preconditioner applications. For a detailed description see \cite{clevenger2021comparison}.

We consider a 3D spherical shell benchmark problem called ``nsinker\_spherical\_shell'' that is part of \texttt{ASPECT}.
A set of 7 heavy sinkers are placed in a spherical shell with inner radius 0.54 and outer radius 1.0. The flow is driven by the density difference of the sinkers and the gravity of magnitude 1.
The viscosity is evaluated in the quadrature points of each cell on the finest level, averaged using harmonic averaging on each cell, and then interpolated to the coarser multigrid levels using the multigrid transfer operators (see also the function \texttt{interpolate\_to\_mg()}
in Figure~\ref{fig:transfer:uml}) for use in the multigrid preconditioner.

\begin{figure}
\centering
\begin{minipage}{0.32\textwidth}
\includegraphics[width=1.0\textwidth]{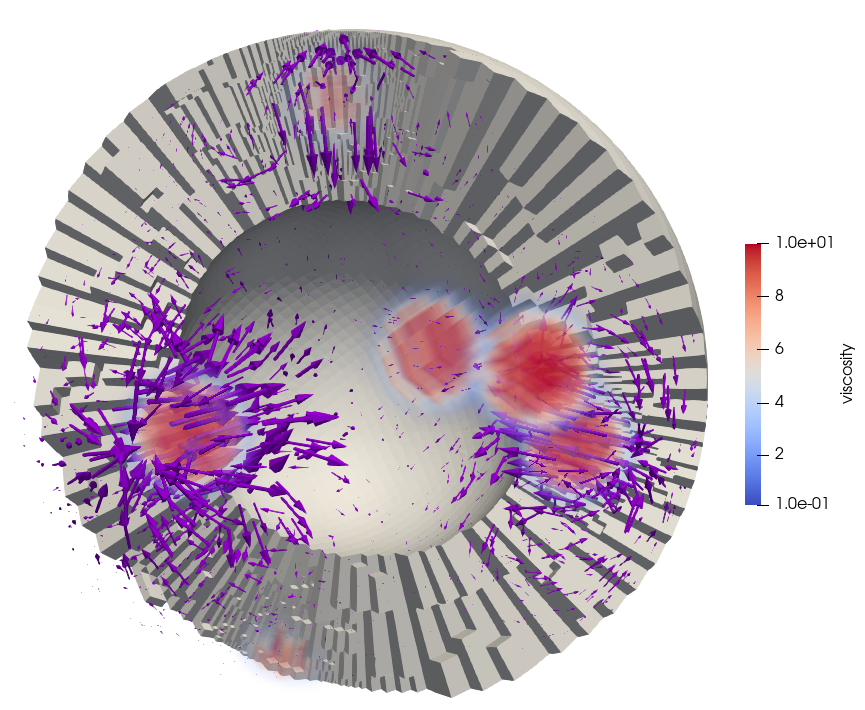}
\end{minipage}
\quad
\begin{minipage}{0.62\textwidth}
\begin{scriptsize}

\begin{tabular}{c|r|crc|crc}
\toprule
  &            & \multicolumn{3}{c|}{global coarsening} & \multicolumn{3}{c}{local smoothing} \\
$L$ & \#DoFs [1e6] & \#it     & solve [s]     & V-cycle [s]     & \#it    & solve [s]    & V-cycle [s]    \\
\midrule
5   & 10.0        & 25       & 1.38        & 0.003        & 25      & 1.15       & 0.006        \\
6   & 20.7        & 25       & 1.61        & 0.006        & 25      & 1.84       & 0.008        \\
7   & 43.2        & 25       & 2.26        & 0.009        & 26      & 2.63       & 0.013        \\
8   & 88.0        & 22       & 2.39        & 0.013        & 27      & 4.15       & 0.022        \\
9   & 178.0       & 25       & 4.62        & 0.020        & 28      & 7.87       & 0.047        \\
10   & 355.0       & 26       & 6.47        & 0.039        & 28      & 13.66      & 0.109        \\
11   & 715.7       & 27       & 10.65       & 0.069        & 27      & 23.67      & 0.214        \\
12   & 1441.4      & 28       & 22.17       & 0.141        & 29      & 50.85      & 0.436        \\
13   & 2896.7      & 28       & 37.38       & 0.266        & 26      & 99.07      & 0.971        \\
14   & 5861.9      & 29       & 77.94       & 0.515        & 26      & 193.19     & 2.060        \\
\bottomrule
\end{tabular}

\end{scriptsize}
\end{minipage}

\caption{Stokes flow in a spherical shell. Left: Visualization of the solution with the high viscosity sinkers (red), velocity vector field in purple, and adaptive mesh in the background. Right: Performance comparison with 7168 processes.}\label{fig:stokes}

\end{figure}

The initial mesh consisting of 96 coarse cells, 4 initial refinement steps and
a high-order manifold description
is refined adaptively using a gradient jump error estimator of the velocity field roughly doubling the number of unkowns in each step, see Figure~\ref{fig:stokes} (left).

We compare number of iterations, time to solution, and time for a single V-cycle of the $A$ block in Figure~\ref{fig:stokes} (right).
Computations are done on TACC Frontera\footnote{Using deal.II master f07477f502 with 64bit indices, Intel 19.1.0.20200306 with -O3 and -march=native} on 7168 processes (128 nodes with 56 cores each).

While iteration numbers are very similar, the simulation with
global-coarsening is about twice as fast in total. Each V-cycle is up to four times
faster, which is not surprising, since the mesh has---for high number of
refinements---a workload efficiency of approximately 20\%. The results suggest that the findings regarding
iteration numbers and performance obtained
in Section~\ref{sec:performance_h} for a simple Poisson problem are also applicable
for this non-trivial problem.

%% file: chapters/conclusions.tex
\section{Summary and outlook}\label{sec:outlook}

{\color{\myred}We have compared geometric local-smoothing and geometric global-coarsening multigrid implementations
for locally refined meshes.
They {\color{\myblue}are based on optimal node-level performance through matrix-free operator evaluation and have been integrated}
into the same finite-element library (\texttt{deal.II}) in order to enable a fair
comparison regarding implementation complexity  and performance.

From the implementation point of view, {the two multigrid versions are---except for some
subtle differences---similar, requiring special treatment} either of refinement edges or of hanging-node constraints during the application of the
smoother and the transfer operator. For the latter, one can usually
rely on existing infrastructure for the application
of constraints, already available from the context
of matrix-free operator evaluations~\cite{Kronbichler2012,munch2021hn}. 
During transfer, global coarsening needs to transfer 
between cells that are refined or not. To be able to
vectorize the cell-local transfer, we categorize
cells within the transfer operator and process cells with
the same category in one go. In the case of local smoothing, it
is possible but not common to repartition the multigrid 
levels; instead, one uses local strategies
for partitioning. {\color{\myblue}For global coarsening, we
investigated 
two partitioning strategies: {\color{\mygreen}one that optimally balances} workload during
smoothing via repartitioning each level and one that minimizes the data to be 
communicated during the transfer phase.}

In a large number of experiments, we have made the observation
that,
for serial simulations, geometric local smoothing is faster than geometric
global coarsening (if the 
number of iterations is the same), since 
the total number of cells to perform smoothing on is less and
the
need for evaluating hanging-node constraints on each level
might be noticeable, particularly for linear shape functions. {For parallel simulations, an equally distributed reduction
of cells is beneficial. Is this not given, load imbalance
is introduced and the critical path of cells, i.e., the time to
solution, is increased. In the case of local smoothing, there might be a non-negligible
number of processes without cells, i.e. work, naturally
introducing load imbalance already on the finest---computationally most expensive---levels. Global coarsening with
repartitioning alleviates this problem and reaches optimal parallel workload, however, comes with the disadvantage of expensive 
transfer steps due to permutation of the data. We made the observation that
global coarsening {\color{\myblue}even if levels are not partitioned optimally
for smoothing} allows---for the
examples considered---to reduce the number of cells on the first levels surprisingly well, introducing only a small
load imbalance, but allowing to perform transfers locally. We
could not make a definite statement 
{\color{\myblue}on the choice of partitioning strategies for global coarsening,}
since
it very much depends on whether the transfer or
the load imbalance is the bottleneck for the given problem.}}

We have also considered polynomial global coarsening ($p$-multigrid). 
Its implementation is conceptually similar to the one of geometric global coarsening so 
that the data structures of geometric global coarsening can be reused and only the 
setup differs. Far from the scaling limit, $p$-multigrid shows better timings per iteration, which can be contributed to the cheaper intergrid transfer.
At the scaling limit, the introduction of additional multigrid levels (when combining $h$ and $p$ multigrid) is noticeable and leads to slower times due to latency.

Global-coarsening algorithms also give the possibility
to simply remove ranks from MPI communicators, allowing to increase
the granularity of the problem on each level and once the problem size
can not be reduced by a sufficient degree anymore, one can simply switch
to the coarse-grid solver even on a level with hanging nodes. Furthermore,
$h$ and $p$ global coarsening could be done in one go: while this might lead to
an overly aggressive coarsening in many cases and, as a result, to a deterioration of the convergence
rate, in other cases, the reduced number of levels could result in  
a reduced latency and, as a consequence, in an improved strong-scaling
behavior. The investigation of these topics is deferred to future work.

In this publication, we focused on
geometrically refined meshes
consisting only of hexahedral shaped cells, where all cells have the same polynomial degree $p$. However, the
presented algorithms work both for $p$- and for $hp$-adaptive problems, where
the polynomial degree varies for each cell, as well as
for simplex or mixed meshes{\color{\myred}, as the implementation in \texttt{deal.II}~\cite{dealII93} shows.
Moreover, they are---as demonstrated
in \cite{Kronbichler2021}---applicable also
for auxiliary-space approximation for DG}.

%% file: chapters/appendix.tex
\begin{acks}

The authors acknowledge collaboration with
Maximilian Bergbauer,
Thomas C. Clevenger,
Ivo Dravins,
Niklas Fehn,
Marc Fehling,
and Magdalena Schreter
as well as the \texttt{deal.II} community.

This work was supported by the Bayerisches Kompetenznetzwerk f\"ur 
Technisch-Wissenschaftliches Hoch- und H\"ochstleistungsrechnen 
(KONWIHR) through the projects ``Performance tuning of
high-order discontinuous Galerkin solvers for SuperMUC-NG''
and ``High-order matrix-free finite element implementations with
hybrid parallelization and improved data locality''. The authors
gratefully acknowledge the Gauss Centre for Supercomputing e.V.
(\url{www.gauss-centre.eu}) for funding this project by providing
computing time on the GCS Supercomputer SuperMUC-NG at
Leibniz Supercomputing Centre (LRZ, \url{www.lrz.de}) through
project id pr83te.

Timo Heister was partially supported by the National Science Foundation (NSF)
Award DMS-2028346, OAC-2015848, EAR-1925575, by the Computational
Infrastructure in Geodynamics initiative (CIG), through the NSF under Award
EAR-0949446 and EAR-1550901 and The University of California -- Davis, and by
Technical Data Analysis, Inc. through US Navy STTR Contract N68335-18-C-0011.
Clemson University is acknowledged for generous allotment of compute time on Palmetto cluster.

\end{acks}

\appendix

\section{APPENDIX TO SECTION~\ref{sec:performance_h}}\label{sec:app:performance_h}
\setcounter{section}{1}


\begin{table}[H]

\centering

\caption{Number of iterations and time to solution for local smoothing (LS) and global coarsening (GC) for the \texttt{octant} simulation case with given analytical solution with 768 processes (16 nodes):
$u(\vec{x})=\left(\frac{1}{\alpha \sqrt{2\pi}}\right)^3\exp\left(-||\vec{x} - \vec{x}_0 || / \alpha^2 \right)$ with $\vec{x}_0 = (-0.5, -0.5, -0.5)^T$ and $\alpha=0.1$. Right-hand-side function $f(\vec{x})$ and inhomogeneous 
Dirichlet boundary conditions have been selected appropriately. }\label{tab:gaussian}

\begin{scriptsize}
\begin{tabular}{c|cc|cc|cc|cc}
\toprule
& \multicolumn{4}{c|}{$p=1$} &  \multicolumn{4}{c}{$p=4$} \\
& \multicolumn{2}{c|}{LS} & \multicolumn{2}{c|}{GC} & \multicolumn{2}{c|}{LS} & \multicolumn{2}{c}{GC}  \\ \midrule
$L$ &  \#$i$ & $t$[s] &  \#$i$ & $t$[s] &  \#$i$ & $t$[s] &  \#$i$ & $t$[s]  \\
\midrule
\input{./global-coarsening-paper-data/small-scaling-gaussian/data.tex}
\bottomrule
\end{tabular}
\end{scriptsize}

\end{table}

\begin{table}[H]

\centering

\caption{Number of iterations and time to solution for local smoothing, global coarsening, and AMG with
768 processes (16 nodes) for the \texttt{octant} simulation case. For local smoothing and global
coarsening, results are shown for \textbf{Chebyshev smoothing degrees} of $k=3$ and $k=6$.}\label{fig:appendix:hmg:smoother}

\begin{scriptsize}
\begin{tabular}{c|cc|cc|cc|cc|cc|cc|cc|cc|cc}
\toprule
& \multicolumn{10}{c|}{$p=1$} &  \multicolumn{8}{c}{$p=4$} \\
& \multicolumn{4}{c|}{local smoothing} &  \multicolumn{4}{c|}{global coarsening} &  \multicolumn{2}{c|}{\multirow{2}{*}{AMG}} & \multicolumn{4}{c|}{local smoothing} &  \multicolumn{4}{c}{global coarsening} \\
& \multicolumn{2}{c|}{$k=3$} & \multicolumn{2}{c|}{$k=6$} & \multicolumn{2}{c|}{$k=3$} & \multicolumn{2}{c|}{$k=6$} & \multicolumn{2}{c|}{} & \multicolumn{2}{c|}{$k=3$} & \multicolumn{2}{c|}{$k=6$} & \multicolumn{2}{c|}{$k=3$} & \multicolumn{2}{c}{$k=6$} \\ \midrule
$L$ &  \#$i$ & $t$[s] &  \#$i$ & $t$[s] &  \#$i$ & $t$[s] &  \#$i$ & $t$[s] &  \#$i$ & $t$[s] &  \#$i$ & $t$[s] &  \#$i$ & $t$[s] &  \#$i$ & $t$[s] &  \#$i$ & $t$[s] \\
\midrule
\input{./global-coarsening-paper-data/parameters-smoother-degree/data.tex}
\bottomrule
\end{tabular}
\end{scriptsize}

\end{table}



\begin{table}[H]

\centering

\caption{Number of iterations and time to solution for local smoothing (LS) and global coarsening (GC) for the \texttt{octant} simulation case with 768 processes (16 nodes). 
All operations {in the outer CG solver are run with double-precision floating-point numbers and} in the multigrid V-cycle are run with the following \textbf{multigrid number types}: single- or double-precision floating-point numbers.}\label{fig:appendix:hmg:mgnumber}

\begin{scriptsize}
\begin{tabular}{c|cc|cc|cc|cc|cc|cc|cc|cc}
\toprule
& \multicolumn{8}{c|}{double} &  \multicolumn{8}{c}{float} \\
& \multicolumn{4}{c|}{$p=1$} &  \multicolumn{4}{c|}{$p=4$} & \multicolumn{4}{c|}{$p=1$} &  \multicolumn{4}{c}{$p=4$} \\
& \multicolumn{2}{c|}{LS} & \multicolumn{2}{c|}{GC} & \multicolumn{2}{c|}{LS} & \multicolumn{2}{c|}{GC} & \multicolumn{2}{c|}{LS} & \multicolumn{2}{c|}{GC} & \multicolumn{2}{c|}{LS} & \multicolumn{2}{c}{GC} \\ \midrule
$L$ &  \#$i$ & $t$[s] &  \#$i$ & $t$[s] &  \#$i$ & $t$[s] &  \#$i$ & $t$[s] &  \#$i$ & $t$[s] &  \#$i$ & $t$[s] &  \#$i$ & $t$[s] &  \#$i$ & $t$[s] \\
\midrule
\input{./global-coarsening-paper-data/parameters-mgnumber/data.tex}
\bottomrule
\end{tabular}
\end{scriptsize}

\end{table}



\begin{table}[H]

\centering

\caption{Number of iterations for local smoothing (LS) and global coarsening (GC) for the \texttt{octant} simulation case and for different \textbf{global relative solver tolerances} with
768 processes (16 nodes).}\label{fig:appendix:hmg:tolerance}

\begin{scriptsize}
\begin{tabular}{c|cc|cc|cc|cc|cc|cc|cc|cc}
\toprule
& \multicolumn{4}{c|}{$10^{-4}$} &  \multicolumn{4}{c|}{$10^{-6}$} &  \multicolumn{4}{c|}{$10^{-8}$} & \multicolumn{4}{c}{$10^{-10}$}  \\
& \multicolumn{2}{c|}{$p=1$} & \multicolumn{2}{c|}{$p=4$} & \multicolumn{2}{c|}{$p=1$} & \multicolumn{2}{c|}{$p=4$} & \multicolumn{2}{c|}{$p=1$} & \multicolumn{2}{c|}{$p=4$} & \multicolumn{2}{c|}{$p=1$} & \multicolumn{2}{c}{$p=4$} \\ \midrule
$L$ &  LS & GC &  LS & GC &  LS & GC &  LS & GC &  LS & GC &  LS & GC &  LS & GC &  LS & GC \\
\midrule
\input{./global-coarsening-paper-data/parameters-tolerance/data.tex}
\bottomrule
\end{tabular}
\end{scriptsize}

\end{table}



\begin{table}[H]

\centering

\caption{Time to solution for global coarsening for the \texttt{octant} simulation case and \textbf{different cell weights} {for cells near hanging nodes compared to regular cells} \textbf{with
768 processes} (16 nodes).}\label{fig:appendix:hmg:weights:small}

\begin{scriptsize}
\begin{tabular}{c|ccccc|ccccc}
\toprule
& \multicolumn{5}{c|}{$p=1$} &  \multicolumn{5}{c}{$p=4$}\\
$L/w$ &  1.0 & 1.5 &  2.0 & 2.5 &  3.0 & 1.0 & 1.5 &  2.0 & 2.5 &  3.0 \\
\midrule
\input{./global-coarsening-paper-data/parameters-weights/data.tex}
\bottomrule
\end{tabular}
\end{scriptsize}

\end{table}

\begin{table}[H]

\centering

\caption{Time to solution for global coarsening for the \texttt{octant} simulation case and \textbf{different cell weights } {for cells near hanging nodes compared to regular cells} with \textbf{24,576 processes} (512 nodes).}\label{fig:appendix:hmg:weights:large}

\begin{scriptsize}
\begin{tabular}{c|ccccc|ccccc}
\toprule
& \multicolumn{5}{c|}{$p=1$} &  \multicolumn{5}{c}{$p=4$}\\
$L/w$ &  1.0 & 1.5 &  2.0 & 2.5 &  3.0 & 1.0 & 1.5 &  2.0 & 2.5 &  3.0 \\
\midrule
\input{./global-coarsening-paper-data/parameters-weights/data-512.tex}
\bottomrule
\end{tabular}
\end{scriptsize}

\end{table}

\begin{table}[H]

\centering

\caption{Number of iterations and time to solution for local smoothing (LS) and global coarsening (GC) for {a uniformly refined mesh of a cube} (\textbf{without hanging nodes}) with
768 processes (16 nodes).}\label{fig:appendix:hmg:hypercube}

\begin{scriptsize}
\begin{tabular}{c|cc|cc|cc|cc}
\toprule
& \multicolumn{4}{c|}{$p=1$} &  \multicolumn{4}{c}{$p=4$} \\
& \multicolumn{2}{c|}{LS} & \multicolumn{2}{c|}{GC} & \multicolumn{2}{c|}{LS} & \multicolumn{2}{c}{GC}  \\ \midrule
$L$ &  \#$i$ & $t$[s] &  \#$i$ & $t$[s] &  \#$i$ & $t$[s] &  \#$i$ & $t$[s]  \\
\midrule
\input{./global-coarsening-paper-data/small-scaling-hypercube/data.tex}
\bottomrule
\end{tabular}
\end{scriptsize}

\end{table}

\section{APPENDIX TO SECTION~\ref{sec:performance_hp}}\label{sec:app:performance_hp}


\begin{table}[H]

\centering

\caption{Number of iterations and time to solution for local smoothing (LS), global coarsening (GC), and AMG as \textbf{coarse-grid solver} of $p$-multigrid with
768 processes (16 nodes) for the \texttt{octant} simulation case for $p=4$. For AMG, different numbers of V-cycles \#$v$ are investigated. AMG parameters used
are shown in Appendix~\ref{sec:app:amg}.}

\begin{scriptsize}
\begin{tabular}{c|cc|cc|cc|cc|cc|cc|cc}
\toprule
& \multicolumn{2}{c|}{\multirow{2}{*}{LS}} &  \multicolumn{2}{c|}{\multirow{2}{*}{GC}} &  \multicolumn{8}{c|}{AMG (ML)} &  \multicolumn{2}{c}{AMG (BoomerAMG)}  \\
 & \multicolumn{2}{c|}{} & \multicolumn{2}{c|}{} & \multicolumn{2}{c|}{\#$v$=1} & \multicolumn{2}{c|}{\#$v$=2} & \multicolumn{2}{c|}{\#$v$=3} & \multicolumn{2}{c|}{\#$v$=4} & \multicolumn{2}{c}{\#$v$=1}\\ \midrule
$L$ &  \#$i$ & $t$[s] &  \#$i$ & $t$[s] &  \#$i$ & $t$[s] &  \#$i$ & $t$[s] &  \#$i$ & $t$[s] &  \#$i$ & $t$[s] &  \#$i$ & $t$[s]   \\
\midrule
\input{./global-coarsening-paper-data/parameters-hp-amg/data.tex}
\bottomrule
\end{tabular}
\end{scriptsize}

\end{table}


\begin{table}[H]

\centering

\caption{Number of iterations and time of a single iteration for $h$-multigrid (local smoothing (LS), global coarsening (GC)) and $p$-multigrid (local smoothing or global coarsening as coarse-grid solver) for $L=9$ and $p=4$.}

\begin{scriptsize}
\begin{tabular}{c|cc|cc|cc|cc}
\toprule
& \multicolumn{4}{c|}{$h$-mg} &  \multicolumn{4}{c}{$p$-mg}  \\
 & \multicolumn{2}{c|}{LS} & \multicolumn{2}{c|}{GC} & \multicolumn{2}{c|}{LS} & \multicolumn{2}{c}{GC}\\ \midrule
\#nodes &  \#$i$ & $t$[s] &  \#$i$ & $t$[s] &  \#$i$ & $t$[s] &  \#$i$ & $t$[s] \\
\midrule
\input{./global-coarsening-paper-data/large-scaling-new/comp.tex}
\bottomrule
\end{tabular}
\end{scriptsize}

\end{table}

\section{AMG Parameters}\label{sec:app:amg}


\vspace{-0.7cm}

\begin{lstlisting}[language=C++, caption={ML~\cite{gee2006ml} (Trilinos 12.12.1)}]
Teuchos::ParameterList parameter_list;
ML_Epetra::SetDefaults("SA", parameter_list);

parameter_list.set("smoother: type", "ILU");
parameter_list.set("coarse: type", coarse_type);
parameter_list.set("initialize random seed", true);
parameter_list.set("smoother: sweeps", 1);
parameter_list.set("cycle applications", 2);
parameter_list.set("prec type", "MGV");
parameter_list.set("smoother: Chebyshev alpha", 10.);
parameter_list.set("smoother: ifpack overlap", 0);
parameter_list.set("aggregation: threshold", 1e-4);
parameter_list.set("coarse: max size", 2000);
\end{lstlisting}

\vspace{-0.7cm}

\begin{lstlisting}[language=C++, caption={BoomerAMG~\cite{falgout2006design} (PETSc 3.14.5)}]
PCHYPRESetType(pc, "boomeramg");

set_option_value("-pc_hypre_boomeramg_agg_nl", "2");
set_option_value("-pc_hypre_boomeramg_max_row_sum", "0.9");
set_option_value("-pc_hypre_boomeramg_strong_threshold", "0.5");
set_option_value("-pc_hypre_boomeramg_relax_type_up", "SOR/Jacobi");
set_option_value("-pc_hypre_boomeramg_relax_type_down", "SOR/Jacobi");
set_option_value("-pc_hypre_boomeramg_relax_type_coarse", "Gaussian-elimination");
set_option_value("-pc_hypre_boomeramg_grid_sweeps_coarse", "1");
set_option_value("-pc_hypre_boomeramg_tol", "0.0");
set_option_value("-pc_hypre_boomeramg_max_iter", "2"));
\end{lstlisting}

%% file: global-coarsening-paper-data/small-scaling-gaussian/data.tex
3  & 4 & 1.4e-3 & 4 & 1.0e-3 & 4 & 2.9e-3 & 4 & 2.9e-3 \\
4  & 4 & 2.8e-3 & 4 & 1.7e-3 & 4 & 4.7e-3 & 4 & 4.2e-3 \\
5  & 4 & 4.0e-3 & 4 & 2.7e-3 & 4 & 6.9e-3 & 4 & 6.5e-3 \\
6  & 4 & 5.7e-3 & 4 & 4.2e-3 & 4 & 1.2e-2 & 4 & 1.0e-2 \\
7  & 4 & 8.2e-3 & 4 & 6.1e-3 & 4 & 2.8e-2 & 4 & 2.4e-2 \\
8  & 5 & 2.2e-2 & 4 & 1.2e-2 & 5 & 2.8e-1 & 4 & 1.7e-1 \\
9  & 5 & 8.8e-2 & 4 & 5.2e-2 & 5 & 2.5e+0 & 4 & 1.6e+0 \\
10  & 5 & 6.1e-1 & 4 & 3.8e-1 & - & - & - & - \\

%% file: global-coarsening-paper-data/parameters-smoother-degree/data.tex
3  & 4 & 1.3e-3 & 3 & 1.4e-3 & 4 & 1.0e-3 & 2 & 7.0e-4 & 1 & 6.0e-3 & 4 & 2.6e-3 & 3 & 2.5e-3 & 4 & 2.5e-3 & 2 & 1.9e-3 \\
4  & 4 & 2.8e-3 & 3 & 2.9e-3 & 4 & 1.7e-3 & 2 & 1.3e-3 & 1 & 2.5e-2 & 4 & 4.3e-3 & 3 & 4.6e-3 & 4 & 4.0e-3 & 2 & 3.0e-3 \\
5  & 4 & 4.3e-3 & 3 & 4.6e-3 & 4 & 2.6e-3 & 2 & 2.3e-3 & 5 & 4.2e-3 & 4 & 6.1e-3 & 3 & 6.6e-3 & 3 & 4.9e-3 & 2 & 5.1e-3 \\
6  & 4 & 5.8e-3 & 3 & 6.4e-3 & 4 & 4.3e-3 & 2 & 3.4e-3 & 8 & 1.5e-2 & 4 & 9.9e-3 & 3 & 1.1e-2 & 3 & 7.2e-3 & 2 & 7.1e-3 \\
7  & 4 & 8.2e-3 & 3 & 9.3e-3 & 3 & 4.5e-3 & 2 & 4.9e-3 & 6 & 1.2e-2 & 4 & 2.7e-2 & 3 & 3.0e-2 & 3 & 1.8e-2 & 2 & 1.7e-2 \\
8  & 4 & 1.7e-2 & 3 & 2.1e-2 & 3 & 9.1e-3 & 2 & 1.0e-2 & 6 & 3.0e-2 & 4 & 2.3e-1 & 3 & 2.6e-1 & 3 & 1.4e-1 & 2 & 1.4e-1 \\
9  & 4 & 7.1e-2 & 3 & 8.6e-2 & 3 & 3.8e-2 & 2 & 4.3e-2 & 7 & 2.1e-1 & 4 & 2.0e+0 & 3 & 2.3e+0 & 3 & 1.2e+0 & 2 & 1.3e+0 \\
10  & 4 & 5.0e-1 & 3 & 6.1e-1 & 2 & 1.9e-1 & 2 & 3.0e-1 & 7 & 1.6e+0 & - & - & - & - & - & - & - & - \\

%% file: global-coarsening-paper-data/parameters-mgnumber/data.tex
3  & 4 & 1.3e-3 & 4 & 9.0e-4 & 4 & 2.6e-3 & 4 & 2.7e-3 & 4 & 1.3e-3 & 4 & 1.0e-3 & 4 & 2.4e-3 & 4 & 2.4e-3 \\
4  & 4 & 2.8e-3 & 4 & 1.6e-3 & 4 & 4.5e-3 & 4 & 4.1e-3 & 4 & 2.7e-3 & 4 & 1.7e-3 & 4 & 4.1e-3 & 4 & 4.1e-3 \\
5  & 4 & 4.0e-3 & 4 & 2.6e-3 & 4 & 6.6e-3 & 3 & 5.1e-3 & 4 & 4.0e-3 & 4 & 2.6e-3 & 4 & 6.1e-3 & 3 & 4.8e-3 \\
6  & 4 & 7.8e-3 & 4 & 4.0e-3 & 4 & 1.1e-2 & 3 & 9.6e-3 & 4 & 6.0e-3 & 4 & 4.1e-3 & 4 & 1.1e-2 & 3 & 7.3e-3 \\
7  & 4 & 8.9e-3 & 3 & 4.9e-3 & 4 & 3.8e-2 & 3 & 2.7e-2 & 4 & 8.3e-3 & 3 & 4.8e-3 & 4 & 2.7e-2 & 3 & 1.8e-2 \\
8  & 4 & 2.1e-2 & 3 & 1.0e-2 & 4 & 4.0e-1 & 3 & 2.4e-1 & 4 & 1.7e-2 & 3 & 9.2e-3 & 4 & 2.3e-1 & 3 & 1.4e-1 \\
9  & 4 & 9.9e-2 & 3 & 5.1e-2 & 4 & 3.2e+0 & 3 & 1.9e+0 & 4 & 7.1e-2 & 3 & 3.9e-2 & 4 & 2.0e+0 & 3 & 1.2e+0 \\
10  & 4 & 7.7e-1 & 2 & 2.8e-1 & - & - & - & - & 4 & 5.0e-1 & 2 & 1.9e-1 & - & - & - & - \\

%% file: global-coarsening-paper-data/parameters-tolerance/data.tex
3  & 4 & 4 & 4 & 4 & 6 & 5 & 6 & 6 & 7 & 6 & 8 & 7 & 9 & 8 & 10 & 9 \\
4  & 4 & 4 & 4 & 4 & 6 & 6 & 6 & 5 & 8 & 7 & 8 & 7 & 10 & 9 & 10 & 9 \\
5  & 4 & 4 & 4 & 3 & 6 & 5 & 6 & 5 & 8 & 7 & 8 & 7 & 10 & 8 & 10 & 8 \\
6  & 4 & 4 & 4 & 3 & 6 & 5 & 6 & 5 & 8 & 7 & 8 & 6 & 10 & 8 & 10 & 8 \\
7  & 4 & 3 & 4 & 3 & 6 & 5 & 6 & 4 & 8 & 6 & 8 & 6 & 10 & 8 & 10 & 8 \\
8  & 4 & 3 & 4 & 3 & 6 & 4 & 6 & 4 & 8 & 6 & 8 & 6 & 9 & 8 & 10 & 7 \\
9  & 4 & 3 & 4 & 3 & 6 & 4 & 6 & 4 & 8 & 6 & 8 & 6 & 9 & 7 & 10 & 7 \\
10  & 4 & 2 & - & - & 6 & 4 & - & - & 8 & 6 & - & - & 9 & 7 & - & - \\

%% file: global-coarsening-paper-data/parameters-weights/data.tex
3  & 1.1e-3 & 1.0e-3 & 1.0e-3 & 1.0e-3 & 1.0e-3 & 2.7e-3 & 2.6e-3 & 2.6e-3 & 2.5e-3 & 2.5e-3 \\
4  & 1.8e-3 & 1.8e-3 & 1.8e-3 & 1.7e-3 & 1.7e-3 & 3.9e-3 & 3.9e-3 & 3.9e-3 & 3.9e-3 & 3.9e-3 \\
5  & 2.5e-3 & 2.7e-3 & 2.7e-3 & 2.6e-3 & 2.8e-3 & 4.4e-3 & 4.7e-3 & 4.6e-3 & 4.5e-3 & 4.6e-3 \\
6  & 4.0e-3 & 4.1e-3 & 4.1e-3 & 4.2e-3 & 4.1e-3 & 7.3e-3 & 7.1e-3 & 7.0e-3 & 7.1e-3 & 7.0e-3 \\
7  & 4.8e-3 & 4.8e-3 & 4.7e-3 & 4.4e-3 & 5.4e-3 & 1.8e-2 & 1.7e-2 & 1.6e-2 & 1.7e-2 & 1.8e-2 \\
8  & 1.1e-2 & 1.0e-2 & 9.1e-3 & 8.5e-3 & 8.2e-3 & 1.1e-1 & 1.0e-1 & 1.0e-1 & 1.0e-1 & 1.0e-1 \\
9  & 5.3e-2 & 4.4e-2 & 3.8e-2 & 3.5e-2 & 3.5e-2 & 9.0e-1 & 8.2e-1 & 7.9e-1 & 7.8e-1 & 7.9e-1 \\
10  & 2.3e-1 & 1.9e-1 & 1.7e-1 & 1.5e-1 & 1.5e-1 & - & - & - & - & - \\

%% file: global-coarsening-paper-data/parameters-weights/data-512.tex
3  & 1.4e-3 & 1.3e-3 & 1.4e-3 & 1.6e-3 & 1.4e-3 & 3.1e-3 & 3.0e-3 & 3.0e-3 & 2.9e-3 & 2.9e-3 \\
4  & 2.1e-3 & 2.2e-3 & 2.1e-3 & 2.0e-3 & 2.1e-3 & 4.3e-3 & 4.4e-3 & 4.3e-3 & 4.3e-3 & 4.3e-3 \\
5  & 2.9e-3 & 2.9e-3 & 2.9e-3 & 2.9e-3 & 2.9e-3 & 4.6e-3 & 5.0e-3 & 4.6e-3 & 5.1e-3 & 4.7e-3 \\
6  & 4.5e-3 & 4.1e-3 & 4.1e-3 & 4.2e-3 & 4.1e-3 & 7.0e-3 & 7.0e-3 & 7.1e-3 & 6.9e-3 & 6.6e-3 \\
7  & 5.9e-3 & 4.8e-3 & 5.7e-3 & 5.1e-3 & 5.6e-3 & 1.1e-2 & 1.1e-2 & 1.1e-2 & 1.1e-2 & 1.0e-2 \\
8  & 7.6e-3 & 9.0e-3 & 7.9e-3 & 9.2e-3 & 8.9e-3 & 1.5e-2 & 1.5e-2 & 1.6e-2 & 1.5e-2 & 1.4e-2 \\
9  & 1.2e-2 & 1.1e-2 & 1.1e-2 & 1.1e-2 & 1.0e-2 & 4.9e-2 & 4.2e-2 & 4.3e-2 & 4.2e-2 & 4.2e-2 \\
10  & 1.8e-2 & 1.5e-2 & 1.4e-2 & 1.3e-2 & 1.3e-2 & 3.7e-1 & 3.4e-1 & 3.4e-1 & 3.4e-1 & 3.4e-1 \\
11  & 7.6e-2 & 6.8e-2 & 5.7e-2 & 5.4e-2 & 5.3e-2 & - & - & - & - & - \\

%% file: global-coarsening-paper-data/small-scaling-hypercube/data.tex
3  & 4 & 1.80e-3 & 4 & 1.80e-3 & 4 & 2.90e-3 & 4 & 3.20e-3 \\
4  & 4 & 2.50e-3 & 4 & 2.40e-3 & 4 & 3.90e-3 & 4 & 4.00e-3 \\
5  & 4 & 3.30e-3 & 4 & 3.20e-3 & 4 & 5.90e-3 & 4 & 6.40e-3 \\
6  & 4 & 4.60e-3 & 4 & 5.40e-3 & 4 & 1.61e-2 & 4 & 1.69e-2 \\
7  & 4 & 9.00e-3 & 4 & 9.00e-3 & 4 & 1.29e-1 & 4 & 1.29e-1 \\
8  & 4 & 4.18e-2 & 4 & 4.21e-2 & 4 & 1.11e+0 & 4 & 1.10e+0 \\
9  & 4 & 3.00e-1 & 4 & 2.98e-1 & - & - & - & - \\

%% file: global-coarsening-paper-data/parameters-hp-amg/data.tex
3  & 4 & 2.60e-3 & 4 & 2.20e-3 & 4 & 1.70e-3 & 4 & 2.00e-3 & 4 & 2.20e-3 & 4 & 2.50e-3 & 4 & 1.70e-3 \\
4  & 4 & 5.10e-3 & 4 & 3.90e-3 & 4 & 5.40e-3 & 4 & 8.30e-3 & 4 & 1.11e-2 & 4 & 1.39e-2 & 4 & 4.90e-3 \\
5  & 4 & 7.60e-3 & 4 & 6.10e-3 & 4 & 1.17e-2 & 4 & 1.95e-2 & 4 & 2.74e-2 & 4 & 3.52e-2 & 4 & 1.17e-2 \\
6  & 4 & 1.11e-2 & 4 & 9.20e-3 & 5 & 2.37e-2 & 4 & 3.24e-2 & 4 & 4.45e-2 & 4 & 5.96e-2 & 5 & 3.60e-2 \\
7  & 4 & 2.35e-2 & 4 & 2.11e-2 & 6 & 4.70e-2 & 5 & 5.47e-2 & 4 & 5.59e-2 & 4 & 6.72e-2 & 6 & 1.04e-1 \\
8  & 4 & 1.81e-1 & 4 & 1.72e-1 & 7 & 3.26e-1 & 5 & 2.65e-1 & 4 & 2.37e-1 & 4 & 2.63e-1 & 7 & 4.76e-1 \\
9  & 4 & 1.53e+0 & 4 & 1.51e+0 & 9 & 3.54e+0 & 7 & 2.97e+0 & 6 & 2.75e+0 & 5 & 2.46e+0 & 8 & 3.83e+0 \\

%% file: global-coarsening-paper-data/large-scaling-new/comp.tex
8 & 4 & 6.30e-1 & 3 & 4.95e-1 & 4 & 4.75e-1 & 4 & 4.66e-1\\
16 & 4 & 4.99e-1 & 3 & 4.03e-1 & 4 & 3.83e-1 & 4 & 3.77e-1\\
32 & 4 & 2.40e-1 & 3 & 2.03e-1 & 4 & 1.87e-1 & 4 & 1.84e-1\\
64 & 4 & 1.39e-1 & 3 & 1.05e-1 & 4 & 9.57e-2 & 4 & 9.25e-2\\
128 & 4 & 5.97e-2 & 3 & 4.90e-2 & 4 & 4.49e-2 & 4 & 4.37e-2\\
256 & 4 & 3.17e-2 & 3 & 2.57e-2 & 4 & 2.35e-2 & 4 & 2.05e-2\\
512 & 4 & 1.54e-2 & 3 & 1.37e-2 & 4 & 1.38e-2 & 4 & 1.24e-2\\